\documentclass[a4paper]{article}

\usepackage{dsfont}
\usepackage{amssymb}
\usepackage{amsthm}
\usepackage{latexsym}
\usepackage{amsmath}
\usepackage{color}
\usepackage{comment}
\usepackage{enumerate}

\newif\ifrs
\rstrue
\ifrs \usepackage{mathrsfs} \fi

\newif\ifcol
\colfalse
\ifcol
\newcommand{\colorr}{\color[rgb]{0.8,0,0}}

\newcommand{\colorb}{\color[rgb]{0,0,0.8}}

\newcommand{\colorn}{\color[rgb]{1,1,1}}

\newcommand{\colory}{\color{yellow}}
\newcommand{\coloro}{\color[rgb]{1,0.4,0}}
\newcommand{\coloroy}{\color[rgb]{1,0.95,0}}
\newcommand{\colorsb}{\color[rgb]{0,0.95,1}}
\newcommand{\colorro}{\color[rgb]{0.851,0.255,0.467}} 
\else
\newcommand{\colr}{\color{black}}
\newcommand{\colorr}{\color{black}}
\newcommand{\colb}{\color{black}}
\newcommand{\colorn}{\color{black}}
\newcommand{\colory}{\color{black}}
\newcommand{\coloro}{\color{black}}
\newcommand{\colo}{\color[rgb]{1,0.851,0}}
\newcommand{\coloroy}{\color{black}}
\newcommand{\colorsb}{\color{black}}
\newcommand{\colorro}{\color{black}}
\newcommand{\coloraka}{\color{black}}
\fi

\excludecomment{en-text}
\includecomment{jp-text}

\excludecomment{comment}

\newtheorem{lemma}{Lemma}

\newtheorem{theorem}{Theorem}
\newtheorem{remark}{Remark}

\def\B{{\bf B}}

\def\D{{\bf D}}
\def\E{{\bf E}}
\def\F{{\bf F}}

\def\I{{\bf I}}

\def\calc{{\cal C}}

\def\cale{{\cal E}}
\def\calf{{\cal F}}

\def\cali{{\cal I}}
\def\calj{{\cal J}}
\def\calk{{\cal K}}
\def\call{{\cal L}}

\def\calo{{\cal O}}

\def\cals{{\cal S}}

\def\sskip{\hspace{0.5cm}}
\def\simleq{ \raisebox{-.7ex}{\em $\stackrel{{\textstyle <}}{\sim}$} }

\def\ep{\epsilon}
\def\half{\frac{1}{2}}
\def\iku{\rightarrow}

\def\up{\uparrow}

\def\y{\vspace*{3mm}\\}
\def\halflineskip{\vspace*{3mm}}
\def\nn{\nonumber}
\def\be{\begin{equation}}
\def\ee{\end{equation}}
\def\bea{\begin{eqnarray}}
\def\eea{\end{eqnarray}}
\def\beas{\begin{eqnarray*}}
\def\eeas{\end{eqnarray*}}
\def\l{\left}
\def\r{\right}

\newcommand{\bbA}{{\mathbb A}}

\newcommand{\bbC}{{\mathbb C}}
\newcommand{\bbD}{{\mathbb D}}
\newcommand{\bbE}{{\mathbb E}}
\newcommand{\bbF}{{\mathbb F}}
\newcommand{\bbG}{{\mathbb G}}
\newcommand{\bbH}{{\mathbb H}}
\newcommand{\bbI}{{\mathbb I}}

\newcommand{\bbN}{{\mathbb N}}

\newcommand{\bbP}{{\mathbb P}}

\newcommand{\bbR}{{\mathbb R}}

\newcommand{\bbW}{{\mathbb W}}

\newcommand{\bbZ}{{\mathbb Z}}

\def\koko{{\colo{koko}}}
\def\bd{\begin{description}}
\def\ed{\end{description}}
\def\partialbs{\backslash\!\!\!\partial}
\def\D2{\bbD_{2,\infty-}}

\def\dotc{\stackrel{\circ}{C}}
\def\dotf{\stackrel{\circ}{F}}
\def\dotw{\stackrel{\circ}{W}}
\def\dote{\stackrel{\circ}{e}}

\begin{document}

\title{
Martingale Expansion in Mixed Normal Limit 
\footnote{
This work was in part supported by 
Grants-in-Aid for Scientific Research No. 19340021, 
No. 24340015 (Scientific Research), 
No. 24650148 (Challenging Exploratory Research); 
the Global COE program ``The Research and Training Center for New Development in Mathematics'' 
of the Graduate School of Mathematical Sciences, University of Tokyo; 
JST Basic Research Programs PRESTO; 
and by a Cooperative Research Program of the Institute of Statistical Mathematics. 
The author thanks NS Solutions Corporation for its support. 
The main parts in this paper were presented at 
Workshop on ``Finance and Related Mathematical and Statistical Issues'',  
September 3-6, 2008, Kyoto Research Park, Kyoto, 
International conference ``Statistique Asymptotique 
des Processus Stochastiques VII'', Universit\'e du Maine, Le Mans, March 16-19, 2009, 
MSJ Spring Meeting 2010, March 24-27, 2010, Keio University, Mathematical Society of Japan, 
and 
International conference ``DYNSTOCH Meeting 2010'', 
Angers, June 16-19, 2010. 
The author thanks to the organizers of the meetings for opportunities of the talks. 
}
\author{Nakahiro YOSHIDA\\
\begin{small}
University of Tokyo
\end{small}\\
\begin{small}
Graduate School of Mathematical Sciences, 
3-8-1 Komaba, Meguro-ku, Tokyo 153, Japan. 
\end{small}\\
\begin{small}
e-mail: nakahiro@ms.u-tokyo.ac.jp
\end{small}
}}
\date{September 14, 2010, \\
Revised version 1: March 14, 2012,\\
{\coloraka Revised version 2: October 1, 2012},\\
Revised version 3: December 31, 2012\\
}
\maketitle
\ \\
{\it Summary} 
The quasi-likelihood estimator and the Bayesian type estimator of the volatility parameter 
are in general asymptotically mixed normal. In case the limit is normal, the asymptotic expansion 
was derived in \cite{Yoshida1997} as an application of the martingale expansion. 
The expansion for the asymptotically mixed normal distribution is then indispensable 
to develop the higher-order approximation and inference for the volatility. 
The classical approaches in limit theorems, where the limit is a process with independent increments 
or a {\colb simple} mixture, 
do not work.  
We present asymptotic expansion of a martingale with asymptotically mixed normal distribution. 
The expansion formula is expressed by the adjoint of a random symbol with coefficients 
described by the Malliavin calculus, 
differently from the standard invariance principle. 
Applications to a quadratic form of a diffusion process (``realized volatility'') are discussed. 
\ \\
\ \\
{\it Keywords and phrases} 
Asymptotic expansion, martingale, mixed normal distribution, 
Malliavin calculus, random symbol, 
double It\^o integral, quadratic form. 
\ \\




\section{Introduction}

The asymptotic expansion is a tool to give a precise approximation to the probability distribution.  
As commonly well known, it is the basement of the contemporary fields in theoretical statistics 
such as the asymptotic decision theory, 
prediction, information criterion, bootstrap and resampling methods, 
information geometry, and stochastic numerical analysis with applications to finance, 
as well as the higher-order approximation of distributions. 
The methodology of the asymptotic expansion has been well established, 
and historically well developed  
especially for independent observations (\cite{BhattacharyaRanga1986}). 

For stochastic processes, there are two principles of asymptotic expansion. 
The first one is the mixing approach. It is also called the {\it local approach} 
because it takes advantage of the factorization of the characteristic function 
by the Markovian property more or less,  
and corresponds to the classical asymptotic expansion for independent models. 
A compilation of the studies for Markovian or near Markovian chains 
is G\"otze and Hipp \cite{GotzeHipp1983}. 
The local nondegeneracy, that is, the decay 
of each factor of the characteristic function becomes an essential problem, 
and G\"otze and Hipp \cite{GotzeHipp1994} showed it for time series models. 
For continuous time, Kusuoka and Yoshida \cite{KusuokaYoshida2000} and 
Yoshida \cite{Yoshida2001b} treated the $\ep$-Markovian process to give asymptotic expansion 
under the mixing condition. 
Then the local nondegeneracy of the Malliavin covariance works 
for objects expressed by the stochastic analysis. 
Though this method is relatively new, 
there are already not a few applications: 
expansion for a functional of the $\ep$-Markov processes, statistical estimators, 
information criterion, stochastic volatility model, empirical distribution 
 (\cite{SakamotoYoshida2004}, \cite{UchidaYoshida2001}, \cite{MasudaYoshida2005}, 
\cite{KutoyantsYoshida2007}  among others).  

Although the local approach is more efficient when one treats mixing processes, 
the asymptotic expansion for the ergodic diffusion was first derived 
through the martingale expansion (\cite{Yoshida1997}, \cite{SakamotoYoshida1998a}, \cite{Yoshida2001c}). 
When the strong mixing coefficient decays sufficiently fast, 
under suitable moments conditions, 
the functionals appearing as the higher-order terms in the stochastic expansion of the functional 
bear a jointly normal limit distribution and it yields a classical expansion formula 
each of whose terms is the Hermite polynomial times the normal density. 
On the other hand, if we consider the sum of quadratics of the increments of a diffusion process, 
it is observed that in the higher-order, the variables have a non-Gaussian limit even when 
the sum is asymptotically normal in the first order. 
Such phenomenon is observed in the estimation of the statistical parameter 
in the essentially linearly parametrized diffusion coefficient. 
Due to the non-Gaussianity, we cannot apply the mixing approach in this case. 
Nonetheless, the martingale expansion can apply to obtain the expansion 
for the estimator of the volatility parameter; see \cite{Yoshida1997}. 
This example shows that the martingale expansion is not inferior to the mixing method 
but even superior in some situations. 

Estimation for the diffusion coefficient 
has been attracting statisticians' interests. 
There are many studies in theoretical statistics such as 
\cite{Dohnal1987}, \cite{PrakasaRao1983}, \cite{PrakasaRao1988}, 
\cite{Yoshida1992b}, \cite{Kessler1997}, 
\cite{SoerensenUchida2003}, 
\cite{IacusUchidaYoshida2009}, 
\cite{Uchida2010}, \cite{ShimizuYoshida2006} 
among others. 
The estimators, including the so-called the realized volatility,  
are in general asymptotically mixed normal; see for example 
\cite{Dohnal1987}, \cite{Genon-CatalotJacod1993}, 
\cite{HayashiYoshida2008a}, 
\cite{UchidaYoshida2008}.  
Today vast literature about this topic is available 
around financial data analysis. We refer the reader to 
\cite{BGJPS2006} and \cite{PodolskijVetter2009} 
and references therein 
for recent advances in the first-order asymptotic theory 
and access to related papers. 
The statistical theory for mixed limits is called the {\it non-ergodic} statistics since 
the Fisher information and the observed information are random even in the limit, differently from 
the classical cases where those are a constant. 
In this sense, it may be said that the mixing approach belongs to the classical theory.
As for the asymptotic expansion in non-ergodic statistics, it seems that 
there is 
room for study. 

In order to explain the technical difficulties in this question, let us recall the method in the proof 
of the martingale central limit theorem and the classical martingale expansion. 
Let $M^n=(M^n_t)_{t\in[0,1]}$ be a continuous martingale; it is possible to consider 
more general local martingale but the existence of jumps for example does not change the situation essentially. 
Under the condition that $\langle M^n\rangle_1\to^p1$ as $n\to\infty$ ,
it suffices to show that 
\beas 
E[e^{iuM^n_1}]e^{\half u^2}=E[e^{iuM^n_1+\half u^2}]\to1
\eeas
as $n\to\infty$ for every $u\in\bbR$.  
Since for large $n$, 
\beas 
E[e^{iuM^n_1+\half u^2}]
&=&
E[e^{iuM^n_1+\half u^2\langle M^n\rangle_1} e^{-\half u^2(\langle M^n\rangle_1-1)}]
\\&\sim&
E[e^{iuM^n_1+\half u^2\langle M^n\rangle_1} ]
\\&=&
1 
\eeas
by the martingale property of the exponential functional if necessary by suitable localization.  
Thus we obtain the central limit theorem. 
Roughly speaking, if evaluating the gap $\langle M^n\rangle_1-1$ more precisely, 
we can obtain the asymptotic expansion. 
In this standard proof in the classical theory of the limit theorems 
for semimartingales converging to a process with independent increments, 
the commutativity of the expectation and $\cale(M^\infty_1)(u)^{-1}$, where 
$\cale(M^\infty_1)(u)$ is the characteristic function of the limit of $M^n_1$, 
was essential. {\coloraka However,} when it is random, this method does not work. 
It explains the difficulty with the asymptotic expansion 
in the {\colb limit with random characteristics.} 

As already mentioned, the asymptotic expansion is inevitable to form 
{\colb modern theories} 
in the non-ergodic statistics, and so it seems natural to try to extend the classical theory of 
asymptotic expansion in the ergodic statistics to a theory that is applicable to the non-ergodic statistics.  
The aim of this article is to present a new martingale expansion to answer this question. 

As a prototype problem, in Section \ref{220915-3}, for a diffusion process 
satisfying the It\^o integral equation
\bea\label{240221-1} 
X_t &=& X_0+\int_0^t b(X_s)ds+\int_0^t\sigma(X_s)dw_s 
\eea
we will consider a quadratic form of the increments of $X$ 
with a strongly predictable kernel:  
\bea\label{240221-10} 
U_n &=&
\sum_{j=1}^n c(X_{t_{j-1}}) (\Delta_j X)^2, 
\eea
where 
$c$ is a function, 
$\Delta_j X=X_{t_j}-X_{t_{j-1}}$ and 
$t_j=j/n$. 
Under suitable conditions, $U_n$ has the in-probability limit 
\bea\label{240221-12} 
U_\infty &=& 
\int_0^1 c(X_s)\sigma(X_s)^2ds. 
\eea
Then the problem is to derive second-order approximation of the distribution of 
$\zeta_n=\sqrt{n}(U_n-U_\infty)$. 
When $c\equiv1$, this gives asymptotic expansion of the realized volatility. 
It turns out that $\zeta_n$ admits a stochastic expansion with a double stochastic integral perturbed by higher-order terms. 
In Section \ref{220915-2}, we will present an expansion formula for a general perturbed double stochastic integral. 
When the function $\sigma$ of (\ref{240221-1}) involves an unknown parameter $\theta$ and $b(X_t)$ is unobservable, 
this gives a semiparametric estimation problem. The error of the quasi maximum likelihood estimator $\hat{\theta}_n$ of $\theta$ 
has a representation by a perturbed double stochastic integral,   
the kernel of which is different from that of {\colb the realized volatility.}
Our result provides asymptotic expansion for 
$\hat{\theta}_n$ though we do not go into this question here. 
Access Section \ref{221003-1}, {\colb Section \ref{240324-2} or} Section \ref{220915-3} first 
if the reader 
wants to know quickly the results in typical problems of a quadratic form, 
while the results presented in the preceding sections are more general. 

The organization of this paper is as follows. 
In Section \ref{220912-1}, we define our object and prepare some notation. 
We introduce the notion of random symbols and 
define an adjoint operation of random symbols in Section \ref{210909-2}. 
The second-order part in the asymptotic expansion consists of two terms and 
they are represented by certain random symbols. 
The first one (adaptive symbol) corresponds to 
the second-order correction term of the classical martingale expansion in \cite{Yoshida1997}, 
while the second one (anticipative symbol) is new. 
Since the characteristics of the targeted distribution are random, it is natural to consider symbols with randomness. 
An error bound of the asymptotic expansion will be presented in Section \ref{240324-1}. 
At this stage, the expression of the second-order term 
associated with the anticipative random symbol 
is not explicit yet. In Section \ref{220915-2}, we treat a variable whose principal part is given 
by a sum of double stochastic integrals. 
It is natural to consider such a variable because it appears 
in the context of the statistical inference for diffusion coefficients and  
realized volatility. 
It turns out that the anticipative random symbol involves objects from the infinite dimensional calculus. 
Section \ref{220826-1} gives expansion for a simple quadratic form of the increments of a Brownian motion. 
Studentization procedure is also discussed there. 
The arguments on the nondegeneracy will be applied in Section \ref{220915-3}, where 
we derive the asymptotic expansion for 
a quadratic form of the increments of a diffusion process. 
Precise approximation to the realized volatility in finance is one of applications of our result  
though our aim is to develop a new methodology in the theory of limit theorems.

\section{Functionals}\label{220912-1}

Let 
$(\Omega,\calf,\F=(\calf_t)_{t\in[0,1]},P)$ be a stochastic basis 
with $\calf=\calf_1$. 
On $\Omega$, 
we will consider a sequence of $d$-dimensional functionals 
each of which admits the decomposition
\beas 
Z_n &=& M_n + W_n+r_n N_n. 
\eeas
For every $n\in\bbN$, $M^n=(M^n_t)_{t\in[0,1]}$ denotes a $d$-dimensional martingale 
with respect to $\F$ and 
$M_n$ denotes the terminal variable of $M^n$, i.e., $M_n = M^n_{1}$. 
In the above decomposition, 
$W_n$, $N_n\in\calf(\Omega;\bbR^d)$\footnote{The set of $d$-dimensional 
measurable mappings.} 
and $(r_n)_{n\in\bbN}$ is a sequence of positive numbers tending to zero 
as $n\to\infty$. 

We will essentially treat a conditional expectation 
given a certain functional 
$F_n\in\calf(\Omega;\bbR^{d_1})$, $n\in\bbN$. 
Later we will specify those functionals more precisely 
to validate computations involving conditional expectations.

{\colorsb{

We write 
\beas 
M_\infty= M^\infty_1, \hspace{5mm}
C^n_t = \langle M^n \rangle_t, \hspace{5mm}
C_n = \langle M^n \rangle_1. 
\eeas
Process $M^\infty$ will be specified later. 
Let $C_\infty\in\calf(\Omega;\bbR^d\otimes\bbR^d)$, 
$W_\infty\in\calf(\Omega;\bbR^d)$ 
and 
$F_\infty\in\calf(\Omega;\bbR^{d_1})$. 
The tangent random vectors are given by  
\beas
\dotc_n&=& r_n^{-1}(C_n-C_\infty),
\\
\dotw_n&=& r_n^{-1}(W_n-W_\infty),
\\
\dotf_n &=& r_n^{-1}(F_n-F_\infty).
\eeas


Consider an extension 
\beas 
(\bar{\Omega},\bar{\calf},\bar{P})=
(\Omega\times{\stackrel{\circ}{\Omega}},\calf\times{\stackrel{\circ}{\calf}},P\times{\stackrel{\circ}{P}})
\eeas
of $(\Omega,\calf,P)$ by a probability space 
$({\stackrel{\circ}{\Omega}},{\stackrel{\circ}{\calf}},{\stackrel{\circ}{P}})$. 
\begin{comment}
{\coloro{
\ \\
-------------------------
?so that 
$\zeta$ defined on $\bar{\Omega}$ depends only on $\omega'\in\Omega'$. 

(or $(M^\infty_\cdot,
N_\infty,\dotc_\infty,\dotw_\infty,\dotf_\infty,
C^\infty_\cdot,W_\infty,F_\infty)$) 
\\ 
------------------------
\ \\
}}
\end{comment}
%
Let 
$M^\infty\in\calf(\bar{\Omega};C([0,1];\bbR^d))$, 
$N_\infty\in\calf(\bar{\Omega};\bbR^d)$, 
$\dotc_\infty\in\calf(\bar{\Omega};\bbR^d\otimes\bbR^d)$, 
$\dotw_\infty\in\calf(\bar{\Omega};\bbR^d)$ 
and 
$\dotf_\infty\in\calf(\bar{\Omega};\bbR^{d_1})$.

We assume
\footnote{[B1](i)=[R1](ii), [B1](ii)=[R1](iii) in \cite{Yoshida2010}.} 
\colr{
\bd
\item[[B1\!\!]] 
{\bf (i)} $({\coloro{M^n}},N_n,\dotc_n,\dotw_n,\dotf_n)
\iku^{d_s(\calf)}
(M^\infty,N_\infty,\dotc_\infty,\dotw_\infty,\dotf_\infty)$. 
\bd
\item[(ii)] $\call\{M^\infty_t|\calf\} = N_d(0,C^\infty_t)$. 
\ed
\ed
}
The convergence in (i) is $\calf$-stable convergence. 
In particular, $(C_n,W_n,F_n)\iku^p(C_\infty,W_\infty,F_\infty)$, 
the variables $C_\infty$, $W_\infty$ and $F_\infty$ are 
$\calf$-measurable, 
and 
\beas
(M_n,N_n,\dotc_n,\dotw_n,\dotf_n,C_n,W_n,F_n)
\iku^d
(M_\infty,N_\infty,\dotc_\infty,\dotw_\infty,\dotf_\infty,
C_\infty,W_\infty,F_\infty). 
\eeas
\vspace{5mm}

Let $\check{\calf}=\calf\vee\sigma[M^\infty_1]$. 
\begin{en-text}
, that is, 
the $\sigma$-field in $\bar{\calf}$ generated by 
the mapping $\check{\Omega}\ni(\omega,\omega')\mapsto(\omega,\zeta(\omega'))
\in\Omega\times\bbR^d$. 
Then there exists a measurable function 
$\check{C}^o_\infty:\Omega\times\bbR^d\iku\bbR^d\otimes\bbR^d$ such that 
\beas 
\check{C}^o_\infty(\omega,\zeta)
&=& 
\bbE[\dotc_\infty|\check{\calf}]. 
\eeas 
If we set 
\beas 
\check{C}_\infty(\omega,z) &=& \check{C}^o_\infty(\omega,C_\infty^{-\half}z), 
\eeas 
then 
\end{en-text}
Then there exists a measurable mapping 
$\check{C}_\infty:\Omega\times\bbR^d\iku\bbR^d\otimes\bbR^d$ 
such that 
\beas 
\check{C}_\infty(\omega,M_\infty)&=&E[\dotc_\infty|\check{\calf}]
\eeas
and we simply write it as $\check{C}_\infty(M_\infty)$. 
%
The uniqueness of the mapping $\check{C}_\infty$ 
is not necessary in what follows. 
In the same way, we define 
$\check{W}_\infty(\omega,z)$, $\check{F}_\infty(\omega,z)$ and 
$\check{N}_\infty(\omega,z)$ so that 
\beas 
\check{W}_\infty(\omega,M_\infty)&=&E[\dotw_\infty|\check{\calf}]
\\
\check{F}_\infty(\omega,M_\infty)
&=&E[\dotf_\infty|\check{\calf}]
\\
\check{N}_\infty(\omega,M_\infty)
&=&E[N_\infty|\check{\calf}].
\eeas
Further, we introduce the notation 
\beas 
\tilde{C}_\infty(z) \equiv  
\tilde{C}_\infty(\omega,z) &:=& \check{C}_\infty(\omega,z-W_\infty)
\eeas
and similarly 
$\tilde{W}_\infty(\omega,z)$, $\tilde{F}_\infty(\omega,z)$ and 
$\tilde{N}_\infty(\omega,z)$.

{\colr{

\section{Random symbol}\label{210909-2}

We will need a notion of random symbols to express the asymptotic expansion formula. 

Given an ${\sf r}$-dimensional Wiener space $(\bbW,\bbP)$ over time interval $[0,1]$, 
the Cameron-Martin subspace is denoted by $H$. 
We will assume that 
the probability space $(\Omega,\calf,P)$ admits the structure such that 
$\Omega=\Omega'\times\bbW$, $\calf=\calf'\otimes\B(\bbW)$ and 
$P=P'\times\bbP$ for some probability space $(\Omega',\calf',P')$, and consider 
the partial Malliavin calculus based on the shifts in $H$. 
See Ikeda and Watanabe \cite{IkedaWatanabe1989}, Nualart \cite{Nualart2006} for the Malliavin calculus. 
$\bbD_{\ell,p}$ denotes the Sobolev space on $\Omega$ 
with indices $\ell\in\bbR$ and $p\in(1,\infty)$. 
Let $\bbD_{\ell,\infty}=\cap_{p\geq2}\bbD_{\ell,p}$.
 
Let $\ell,{\mathfrak m}\in\bbZ_+$. 
Denote by $\calc(\ell,{\mathfrak m})$ the set of 
functions $c:\bbR^d\to\bbD_{\ell,\infty}$ 
satisfying the following conditions: 
\begin{description}
\item[(i)] $c\in C^{\mathfrak m}(\bbR^d;\bbD_{\ell,\infty})$, that is, 
for every $p\geq2$, 
$c:\bbR^d\to\bbD_{\ell,p}$ is Fr\'echet differentiable 
${\mathfrak m}$ times in $z$ and all the derivatives 
(each element taking values in $\bbD_{\ell,\infty}$) 
up to order ${\mathfrak m}$ are continuous. 

\item[(ii)] For every $p\geq2$, 
\beas 
\|\partial_z^\mu c(z) \|_{\ell,p} 
&\leq& 
C_{\ell,p}(1+|z|)^{C_{\ell,p}}
\sskip(z\in\bbR^d,\  \mu\in\bbZ_+^d 
\mbox{with }|\mu|\leq {\mathfrak m})
\eeas
for some constant $C_{\ell,p}$ depending on $\ell$ and $p$. 
\end{description}
Note that if $c=c(z)$ does not depend on $z$, then this condition is reduced to 
that $\|c\|_{\ell,p}<\infty$. 
The sum of the elements of $\alpha\in\bbZ_+^d$, i.e., the length of $\alpha$, is denoted by $|\alpha|$. 

Let $\ell,{\mathfrak m},{\mathfrak n}\in\bbZ_+$. 
A function 
$\varsigma:\Omega\times\bbR^d\times (i\>\bbR^d)\times(i\>\bbR^{d_1})
\to\bbC$ 
is called a {\it random symbol}. 
We say that a random symbol $\varsigma$ is of {\it class} $(\ell,{\mathfrak m},{\mathfrak n})$ 
if $\varsigma$ admits a representation 
\bea\label{210815-1} 
\varsigma(z,iu,iv)&=&\sum_j c_j(z)(iu)^{m_j}(iv)^{n_j}
\sskip(\mbox{finite sum})
\eea
for some $c_j\in\calc(\ell,{\coloraka |}m_j{\coloraka |})$, 
$m_j\in\bbZ_+^d$ ($|m_j|\leq {\mathfrak m}$) and 
$n_j\in\bbZ_+^{d_1}$ ($|n_j|\leq {\mathfrak n}$). 
We denote by $\cals(\ell,{\mathfrak m},{\mathfrak n})$ 
the set of random symbols of class $(\ell,{\mathfrak m},{\mathfrak n})$.

The Malliavin covariance matrix of a multi-dimensional functional $F$ is denoted by $\sigma_F$,  
and we write $\Delta_F=\det\sigma_F$. 
Suppose that {\coloraka $\ell\geq \ell_0:=2([d_1/2]+1+[({\mathfrak n}+1)/2])$}
\footnote{{\coloraka The number 
$\ell_0=2[({\mathfrak n}+d_1+2)/2]$ is possible to use for this $\ell_0$ if we regard $\partial^{n_j}\delta_x$ as a Schwartz distribution.}} 
and the following conditions are satisfied: 
\begin{description}
\item[(i)] 
$F_\infty\in\bbD_{\ell+1,\infty}(\bbR^{d_1})$, 
$W_\infty\in\bbD_{\ell,\infty}(\bbR^d)$ and 
$C_\infty\in\bbD_{\ell,\infty}(\bbR^d\otimes\bbR^d)$; 

\item[(ii)] 
$\Delta_{F_\infty}^{-1},\>\det C_\infty^{-1}\in\cap_{p\geq2} L^p$; 

\item[(iii)] 
$\varsigma\in\cals(\ell,{\mathfrak m},{\mathfrak n})$. 
\end{description}
For the random symbol $\varsigma\in\cals(\ell,{\mathfrak m},{\mathfrak n})$ 
taking the form of (\ref{210815-1}), 
we define the adjoint operator $\varsigma^*$ of $\varsigma$ as follows. 
It applies 
to $\phi(z;W_\infty,C_\infty)\delta_x(F_\infty)$ as 
\bea\label{210815-2} 
\varsigma(z,\partial_z,\partial_x)^*
\Big\{
\phi(z;W_\infty,C_\infty)\delta_x(F_\infty)\Big\}
&=&
\sum_j (-\partial_z)^{m_j}(-\partial_x)^{n_j} 
\Big(
c_j(z)
\phi(z;W_\infty,C_\infty)\delta_x(F_\infty)\Big). 
\eea
Here 
$\phi(z;\mu,C)$ is the normal density with mean vector $\mu$ and covariance matrix $C$, and 
the derivatives are interpreted as the Fr\'echet derivatives 
in the space of the generalized Wiener functionals of Watanabe (\cite{Watanabe1983}, \cite{IkedaWatanabe1989}). 
By assumptions, 
$(-\partial_z)^{m_j}\phi(z;W_\infty,C_\infty)\in\bbD_{\ell,\infty}$ 
and 
$(-\partial_x)^{n_j}\delta_x(F_\infty)=(\partial^{n_j}\delta_x)(F_\infty)
\in\bbD_{-\ell,\infty}$. 
$p^{F_\infty}$ denotes the density of $F_\infty$. 
The generalized expectation of the formula (\ref{210815-2}) 
gives a formula with the usual expectation: 
\beas &&
E\Big[
\varsigma(z,\partial_z,\partial_x)^*
\Big\{
\phi(z;W_\infty,C_\infty)\delta_x(F_\infty)\Big\}
\Big]
\\&=& 
\sum_j (-\partial_z)^{m_j}(-\partial_x)^{n_j} 
E\Big[c_j(z)
\phi(z;W_\infty,C_\infty)\delta_x(F_\infty)
\Big]
\\&=& 
\sum_j (-\partial_z)^{m_j}(-\partial_x)^{n_j} 
\bigg(E\big[c_j(z)
\phi(z;W_\infty,C_\infty)\big|F_\infty=x\big]p^{F_\infty}(x)\bigg)
\\&=:&
E\bigg[\varsigma(z,\partial_z,\partial_x)^*
\bigg\{
\phi(z;W_\infty,C_\infty)\bigg|F_\infty=x\bigg]p^{F_\infty}(x)
\bigg\}. 
\eeas

For $m'\in\bbZ_+^d$ and $n'\in\bbZ_+^{d_1}$, 
\bea\label{210923-1} 
&&
\sup_{z,x}
\Big| 
z^{m'}x^{n'}
E\Big[(-\partial_z)^{\mu}c_j(z)\cdot
(-\partial_z)^{m_j-\mu}\phi(z;W_\infty,C_\infty)
(-\partial_x)^{n_j}\delta_x(F_\infty)\Big]
\Big|
%
%
\nn\\&\leq& 
\sup_{z,x}
\|(-\partial_z)^{\mu}c_j(z)\|_{\ell,p}
\ 
\|z^{m'}(-\partial_z)^{m_j-\mu}\phi(z;W_\infty,C_\infty)\|_{\ell,p_1}
\ 
\|F_\infty^{n'}\|_{\ell,p_2}
\ 
\|(-\partial_x)^{n_j}\delta_x(F_\infty)\|_{-\ell,p_3}
\nn\\&\simleq&
\sup_{z,x}\{(1+|z|)^{C_{\ell,p}}\cdot (1+|z|)^{-C_{\ell,p}} \cdot 1\}
\nn\\&<&
\infty, 
\eea
where $p,p_1,p_2,p_3>1$ with $p^{-1}+p_1^{-1}+p_2^{-1}+p_3^{-1}=1$. 
Consequently, we can apply the Fourier transform to obtain 
\begin{en-text}
if a random symbol $\varsigma$ has the form 
\beas 
\varsigma(z,iu,iv)&=&\sum_j c_j(z)(iu)^{m_j}(iv)^{n_j}
\eeas
for random mesurable functions 
$c_j:\Omega\times\bbR^d\iku\bbR$, 
where $m_j\in\bbZ_+^d$ and $n_j\in\bbZ_+^{d_1}$. 
%
%
$\varsigma(x,iu,iv)=\sum_j c_j(x)(iu)^{m_j}(iv)^{n_j}$, 
\end{en-text}
%
%
\begin{en-text}
\beas 
&&
\calf_{(z,x)}
\bigg[\!\!\bigg[
E\bigg[ \varsigma(z,x,\partial_z,\partial_x)^* 
\bigg\{\phi(z;W_\infty,C_\infty) \big| 
F_\infty=x \bigg] 
p^{F_\infty}(x)\bigg\} 
\bigg]\!\!\bigg]
(u,v)
\\&=&
\int \exp(iu\cdot z+iv\cdot x) 
E\bigg[ \varsigma(z,x,\partial_z,\partial_x)^* 
\bigg\{\phi(z;W_\infty,C_\infty) \big| 
F_\infty=x \bigg] 
p^{F_\infty}(x)\bigg\} dzdx
\\&=&
E\bigg[ \int_{\bbR^d}
\exp(iu\cdot z+iF_\infty[v]) \phi(z;W_\infty,C_\infty) 
\varsigma(z,F_\infty,iu,iv) dz \bigg]
\eeas 
\end{en-text}
\bea\label{210922-5}
&&
\calf_{(z,x)}
\bigg[\!\!\bigg[
E\bigg[ \varsigma(z,\partial_z,\partial_x)^* 
\bigg\{
{\colb \phi(z;W_\infty,C_\infty) \delta_x(F_\infty)}
\bigg\} \bigg]
\bigg]\!\!\bigg]
(u,v)
\nn\\&=&
\int \exp(iu\cdot z+iv\cdot x) 
E\bigg[ \varsigma(z,\partial_z,\partial_x)^* 
\bigg\{
{\colb \phi(z;W_\infty,C_\infty) \delta_x(F_\infty)}
\bigg\}\bigg]
dzdx
\nn\\&=&
E\bigg[ \int_{\bbR^d}
\exp(iu\cdot z+iF_\infty[v]) \phi(z;W_\infty,C_\infty) 
\varsigma(z,iu,iv) dz \bigg].
\eea
{\colb{
More generally, duality argument also yields 
\beas
\calf_{(z,x)}
\bigg[\!\!\bigg[
\varsigma(z,\partial_z,\partial_x)^* 
\bigg\{
{\colb \phi(z;W_\infty,C_\infty) \delta_x(F_\infty)}
\bigg\} 
\bigg]\!\!\bigg]
(u,v)
&=&
\calf_{z}
\bigg[\!\!\bigg[
\varsigma(z,iu,iv)\phi(z;W_\infty,C_\infty) e^{iu\cdot F_\infty}
\bigg]\!\!\bigg](u). 
\eeas
}}


\section{Asymptotic expansion formula}\label{240324-1}

Nondegeneracy of the targeted distribution is indispensable for asymptotic expansion. 
However, the complete nondegeneracy, that implies absolute continuity, is not necessary, nor can we assume 
in statistical inference. For example, the maximum likelihood estimator does not admit a density in general;  
even existence of itself is not ensured on the whole probability space. 
Besides, the complete nondegeneracy is often hard to prove although it would be possible, and this restricts applications of the result. 
On the other hand, partial nondegeneracy is easier to work with. There the localization method plays an essential role. 
The localization is realized through a sequence of $\calf$-measurable truncation functionals $\xi_n$. 
In applications, we construct a suitable $\xi_n$ to validate necessary nondegeneracy of the functional in question.  

We will assume the following conditions.\footnote{[B2]$_\ell$(i),(ii) are [S]$_{\ell,{\mathfrak m},{\mathfrak n}}$(i),(ii) of \cite{Yoshida2010}, respectively. 
[B3](i),(ii),(iii) are [S]$_{\ell,{\mathfrak m},{\mathfrak n}}$(iii),(iv),(v) of \cite{Yoshida2010}, respectively.}

\bd
\item[[B2\!\!]]$_ {\ell}$
{\bf (i)}  
$F_\infty\in\bbD_{\ell+1,\infty}(\bbR^{d_1})$, 
$W_\infty\in\bbD_{{\coloro{\ell+1}},\infty}(\bbR^d)$ 
and $C_\infty\in\bbD_{\ell,\infty}(\bbR^d\otimes\bbR^d)$. 

\bd
\item[(ii)] 
{\colb $M_n\in\bbD_{{\ell+1},\infty}(\bbR^d)$, } 
$F_n\in\bbD_{\ell+1,\infty}(\bbR^{d_1})$, 
$W_n\in\bbD_{{\coloro{\ell+1}},\infty}(\bbR^d)$, 
$C_n\in\bbD_{\ell,\infty}(\bbR^d\otimes\bbR^d)$, 
$N_n\in\bbD_{{\coloro{\ell+1}},\infty}(\bbR^d)$  
and $\xi_n\in\bbD_{\ell,\infty}(\bbR)$. 
Moreover, 
\beas 
\sup_{{\coloraka n\in\bbN}}\Big\{ 
\|M_n\|_{\ell+1,p}+\|\dotc_n\|_{\ell,p}+\|\dotw_n\|_{{\coloro{\ell+1}},p}
+\|\dotf_n\|_{\ell+1,p}+\|N_n\|_{{\coloro{\ell+1}},p}+\|\xi_n\|_{\ell,p} 
\Big\} <\infty
\eeas
for every $p\geq2$. 
\ed
\ed

\bd
\item[[B3\!\!]] 
{\bf (i)}  
{\coloro{$\lim_{n\to\infty}P\l[|\xi_n|\leq\half\r]=1$. }}
\bd
\item[(ii)] 
$|C_n-C_\infty|>r_n^{1-a}$ 
implies $|\xi_n|\geq 1$, 
where $a\in(0,{\colorr{1/3}})$ is a constant. 

\item[(iii)] For every $p\geq2$, 
\beas 
\limsup_{n\to\infty}E\Big[1_{\{|\xi_n|\leq1\}} \Delta_{(M_n+W_\infty,F_\infty)}^{-p}\Big] 
<\infty
\eeas
and moreover 
$\det C_\infty^{-1}\in\cap_{p\geq2} L^p$. 

\end{description}
\ed

{\colb We have $\Delta_{F_\infty}^{-1}\in\cap_{p\geq2} L^p$ 
by Fisher's inequality and Fatou's lemma. }

{\colb{We}} will give an asymptotic expansion formula 
to approximate the 
joint 
distribution of $(Z_n,F_n)$. 
Let 
\bea\label{240218-1} 
\underline{\sigma}(z,iu,iv)
&=&
\half\tilde{C}_\infty(z)^{j,k}(iu_j)(iu_k)
+
\tilde{W}_\infty(z)^j (iu_j)
\nn\\&&
+
\tilde{N}_\infty(z)^j (iu_j)
+
\tilde{F}_\infty(z)^l (iv_l)
\eea
for $u\in\bbR^d$ and $v\in\bbR^{d_1}$.\footnote{Einstein's rule for repeated indices. } 

The symbol (\ref{240218-1}) is the adaptive random random symbol, which corresponds to the second-order correction term of  
the asymptotic expansion in the normal limit (\cite{Yoshida1997}). 
In the mixed normal limit case, we need an additional one referred to as the anticipative random symbol 
and denoted by $\overline{\sigma}(z,iu,iv)$. 
As will be seen in applications in this paper, the anticipative random symbol is given by 
the Malliavin derivatives, which shows non-classical nature of the present asymptotic expansion 
beyond the standard invariance principle. 

We shall define the anticipative random symbol. 
Let 
\beas
\Psi_\infty(u,v)
&=&
\exp\l\{ {\colorr iW_\infty[u]}-\half C_\infty[u^{\otimes2}]+iF_\infty[v] \r\}.
\eeas
Moreover, set 
$
L^n_t(u) 
=
e^n_t(u)-1
$
with 
\beas 
e^n_t(u)
&=&
\exp\l(iM^n_t[u]+\half C^n_t[u^{\otimes2}]\r). 
\eeas
We write 
$\partialbs^\alpha=i^{-|\alpha|}\partial^\alpha$, 
$\partialbs^{\alpha_1}_u=i^{-|\alpha_1|}\partial_u^{\alpha_1}$, 
$\partialbs^{\alpha_2}_v=i^{-|\alpha_2|}\partial_v^{\alpha_2}$ 
and 
$
\partialbs^\alpha=\partialbs_u^{\alpha_1}\partialbs_v^{\alpha_2}
$ 
for {\colb the multi-index} 
$\alpha=(\alpha_1,\alpha_2)\in\bbZ_+^{\check{d}}$, $\check{d}=d+d_1$. 
Let $\psi\in C^\infty(\bbR;[0,1])$ such that $\psi(x)=1$ if $|x|\leq1/2$ and $\psi(x)=0$ if $|x|\geq1$ 
and let $\psi_n=\psi(\xi_n)$ for $n\in\bbN$. 
Furthermore, let 
\beas 
\Phi^{2,\alpha}_n(u,v)&=&\partialbs^\alpha E[L^n_1(u)\Psi_\infty(u,v)\psi_n].
\eeas
If $\psi_n$ is equal to one, in fact it is usually so without large deviation probability, 
and $W_\infty$, $C_\infty$ and $F_\infty$ are constants, as this is the case in the normal limit case, 
then $\Phi^{2,\alpha}_n=0$ since $L^n(u)$ becomes a mean zero martingale. 
However $\Phi^{2,\alpha}_n$ does not vanish in general. 
We may intuitively say that $\Phi^{2,\alpha}_n$ measures the torsion of martingales under the shift of 
measure $P$ by $\Psi_\infty(u,v)$. 
The effect of this torsion appears quite differently, depending on the cases. 
{\colb Thus, we will treat the effect in a slightly abstract shape for a while.}
The adaptive random symbol describes it as (ii) of the following condition.
\footnote{
[B4]$_{\ell,{\mathfrak m},{\mathfrak n}}$(i),(ii) are [S]$_{\ell,{\mathfrak m},{\mathfrak n}}$(vi),(vii) of \cite{Yoshida2010}.}
{\coloraka Let $\ell_*=2[d_1/2]+4$.} 

\bd
\item[[B4\!\!]]$_{\ell,{\mathfrak m},{\mathfrak n}}$
{\bf (i)} $\underline{\sigma}\in{\coloraka \cals(\ell_*,2,1)}$.
\footnote{{\coloraka The number $\ell_*=2[(d_1+3)/2]$ is possible as the previous footnote. }} 
\item[(ii)] 
There exists a random symbol $\overline{\sigma}\in\cals(\ell,{\mathfrak m},{\mathfrak n})$ 
admitting a representation 
\beas 
\overline{\sigma}(iu,iv)
&=& 
\sum_jc_j(iu)^{m_j}(iv)^{n_j}\sskip(\mbox{finite sum})
\eeas
for some random variables $c_j$ and satisfying 
\bea\label{220131-1}
\lim_{n\to\infty}r_n^{-1}\Phi^{2,\alpha}_n(u,v)
&=&
\partialbs^\alpha E\bigg[\int_{\bbR^d}
\exp\big(iu\cdot z+iF_\infty[v]\big)
\phi(z;W_\infty,C_\infty)dz\>\overline{\sigma}(iu,iv)\bigg]
\nn\\&{\coloraka =}&
{\coloraka
\partialbs^\alpha
E\bigg[\exp\bigg(iW_\infty[u]-\half C_\infty[u^{\otimes2}]+iF_\infty[v]\bigg)\>\overline{\sigma}(iu,iv)\bigg]
}
\eea
for $u\in\bbR^d$, $v\in\bbR^{d_1}$ and $\alpha=(\alpha_1,\alpha_2)\in\bbZ^{\check{d}}_+$. 
\ed

\begin{en-text}
In particular, 
\begin{description}
\item[(i)] 
$(C^n_\cdot,W_n,F_n)\iku^p(C^\infty_\cdot,W_\infty,F_\infty)$. 

\item[(iii)] 
$\call\{M^\infty_t|\calf\} = N_d(0,C^\infty_t)$. 
\end{description}

In particular, 
\beas &&
(M^n_\cdot,N_n,\dotc_n,\dotw_n,\dotf_n,C^n_\cdot,W_n,F_n)
\iku^d
(M^\infty_\cdot,N_\infty,\dotc_\infty,\dotw_\infty,\dotf_\infty,
C^\infty_\cdot,W_\infty,F_\infty)
\eeas
\end{en-text}

\begin{comment}
{\coloro{------------------------------------

$M_\infty=M^\infty_1=^d C_\infty^\half[\zeta,\cdot]$ 
for $\zeta=B_1$. 

---------------------------------------------
}}
\end{comment}

\begin{en-text}
$\xi_n$: $\calf$-measurable smooth 
(or finitely many times differentialble) functional such that 
\begin{description}
\item[(i)] 
$
P\l[|\xi_n|>\half\r]=O(r_n^{1+\kappa})
$, where 
$\kappa$ is related to other parameters in practice. 

\item[(ii)] 
$|C_n-C_\infty|>r_n^{1-a}$ or ... 
implies $|\xi_n|\geq 1$, 
where $a\in(0,1)$ is a constant. 
\end{description}
\end{en-text}

\begin{en-text}
Let $\xi_n$ $(n\in\bbN)$ be $\calf$-measurable functionals 
satisfying the following conditions. 

{\coloro{
\bd
\item[[R2\!\!]] 
{\bf (i)} $\{\xi_n;n\in\bbN\}$ is bounded in $L^{\infty-}$. 
\begin{description}

\item[(ii)] 
{\coloro{$\lim_{n\to\infty}P\l[|\xi_n|\leq\half\r]=1$. }}

\item[(iii)] 
$|C_n-C_\infty|>r_n^{1-a}$ 
implies $|\xi_n|\geq 1$, 
where $a\in(0,{\colorr{1/3}})$ is a constant. 
\end{description}
\ed
}}
\ \\

We write $\psi_n=\psi(\xi_n)$. 
\end{en-text}

\begin{en-text}
$(M_\infty,N_\infty,\dotc_\infty,\dotw_\infty,\dotf_\infty)$ are 
defined on the extention 
$(\Omega\times\Omega',\calf\times\calf',\bbP)$. 

\beas 
\check{\calf}=\calf\vee\sigma[\zeta]
\eeas

\beas 
\check{C}_\infty(\omega,z)
&=&\bbE[\dotc_\infty|\calf\vee\sigma[\zeta]]
|_{(\omega,\zeta)=(\omega,z)}
\\
&=&
\int_\Omega\int_{\Omega'} \dotc_\infty(\omega,\omega')
\bbP(d\omega,d\omega'|\zeta=z)
\eeas
($\calf$ is extended over $\tilde{\Omega}$. 
$\zeta$ is defined on $\tilde{\Omega}$. )
Then 
\beas 
\bbE_{\check{\calf}}[\dotc_\infty]
=
\check{C}_\infty(\cdot,\zeta)
=\check{C}_\infty(\zeta)\sskip\mbox{simply write}
\eeas 
In the same way, we define
\beas 
\check{W}_\infty(\omega,z)
&=&\bbE[\dotw_\infty|\calf\vee\sigma[\zeta]]
|_{(\omega,\zeta)=(\omega,z)}
\\
\check{F}_\infty(\omega,z)
&=&\bbE[\dotf_\infty|\calf\vee\sigma[\zeta]]
|_{(\omega,\zeta)=(\omega,z)}
\\
\check{N}_\infty(\omega,z)
&=&\bbE[N_\infty|\calf\vee\sigma[\zeta]]
|_{(\omega,\zeta)=(\omega,z)}
\eeas
\end{en-text}
%

Set 
$\tilde{\Phi}^{2,\alpha}(u,v)
=
\lim_{n\iku\infty} r_n^{-1}
\Phi^{2,\alpha}_n(u,v)$. 
%

\begin{en-text}
\subsection{Second-order term of type II (Delete this subsection!?)}


Suppose that 
[B1] 
[B2]$_\ell$ for some $\ell\geq0$  
and [B3] (i), (ii) 
are satisfied. 
We assume that there exists the limit 
\bea\label{220131-1}
\tilde{\Phi}^{2,\alpha}(u,v)
&:=&
\lim_{n\iku\infty} r_n^{-1}
\partialbs^\alpha\Phi^{2,0}_n(u,v)
\nn\\&=&
\lim_{n\iku\infty} 
r_n^{-1}
\partialbs^\alpha
E\l[L^n_1(u)\Psi_\infty(u,v)\psi_n\r]
\nn\\&=&
\partialbs^\alpha
E\bigg[ \int_{\bbR^d}\exp\big(iu\cdot z+iF_\infty[v] \big) 
\phi(z;W_\infty,C_\infty)\>dz 
\>\bar{\sigma}(iu,iv)
\bigg]
\eea
for some random symbol $\bar{\sigma}(iu,iv)$ 
which admits a representation
\beas
\bar{\sigma}(iu,iv)&=&\sum_j c_j(iu)^{m_j}(iv)^{n_j}
\sskip(\mbox{finite sum}), 
\eeas
where $c_j$ are random variables satisfying 
$c_j\in\cap_{p>1}L^p$. 
\end{en-text}

\begin{remark}\rm 
It is also possible to consider a more general random symbol 
$\bar{\sigma}(z,x,iu,iv)$. 
\end{remark}

\begin{remark}\rm
{\colb As mentioned above,}
$\tilde{\Phi}^{2,\alpha}$ vanishes 
in the classical case of deterministic $\Psi_\infty(u,v)$ 
thanks to the local martingale property of $L^n_t(u)$. 
{\colb Non-vanishing case will appear later. }
\end{remark}

The full random symbol for the second-order terms is 
\bea\label{220829-1}
\sigma&:=& \underline{\sigma}+\bar{\sigma}. 
\eea
%
In order to approximate the joint local density of $(Z_n,F_n)$, 
we use the density function 
\beas 
p_n(z,x) 
&=& 
E\bigg[
\phi(z;W_\infty,C_\infty)\bigg|F_\infty=x\bigg]
p^{F_\infty}(x)
\\&&
+
r_n 
E\bigg[\sigma(z,\partial_z,\partial_x)^*
\bigg\{
\phi(z;W_\infty,C_\infty)\bigg|F_\infty=x\bigg]
p^{F_\infty}(x)
\bigg\}. 
\eeas
Note that $p_n(z,x)$ is well defined under 
{\coloraka [B2]$_\ell$, [B3] and [B4]$_{\ell,{\mathfrak m},{\mathfrak n}}$} when 
{\coloraka $\ell\geq\ell_0$}. 
We should remark that $p_n(z,x)$ is written in terms of the conditional expectation given $F_\infty$. 
This suggests that conditioning by $F_\infty$ or equivalently by $F_n$ under truncation is essentially used in validation of the formula.  
With Watanabe's delta functional, we can write 
\bea\label{220829-2} 
p_n(z,x) 
&=& 
E\bigg[
\phi(z;W_\infty,C_\infty)\delta_x(F_\infty)\bigg]
\nn\\&&
+
r_n 
E\bigg[\sigma(z,\partial_z,\partial_x)^*
\bigg\{
\phi(z;W_\infty,C_\infty)\delta_x(F_\infty)\bigg \}\bigg]. 
\eea

For $M,\gamma>0$, let 
$\cale(M,\gamma)$ denote the set of measurable functions $f:\bbR^{\check{d}}\to\bbR$ 
satisfying $|f(z,x)|\leq M(1+|z|+|x|)^\gamma$. 
Let 
\beas 
\Delta_n(f)
&=&
\bigg| E\bigl[f(Z_n,F_n)\bigr] 
- \int f(z,x)p_n(z,x)\>dzdx \bigg|. 
\eeas
{\coloraka{
Let $\Lambda^0_n(d,q)=\{u\in\bbR^d;|u|\leq r_n^{-q}\}$,  
where 
$q=(1-a)/2\in(1/3,1/2)$. 
}}
\noindent 
Let 
\beas 
\ep(k,n)
&=&
\max_{\alpha:|\alpha|\leq {\coloro{k}}}
\frac{1}{(2\pi)^{\check{d}}}
r_n
\int_{\Lambda^0_n(\check{d},q)
}
\big|r_n^{-1}\Phi^{2,\alpha}_n(u,v)-\tilde{\Phi}^{2,\alpha}(u,v)\big|dudv. 
\eeas 
The following theorem gives an extension of Theorem 4 of \cite{Yoshida1997}. 

\begin{theorem}\label{210215-1}
{\coloro{Let $\ell={\coloraka \ell_0}
\vee(\check{d}+3)$.}}
Suppose that {\rm [B1], [B2]$_\ell$, [B3] and [B4]$_{\ell,{\mathfrak m},{\mathfrak n}}$} are fulfilled. 
Let $M,\gamma\in(0,\infty)$ and $\theta\in(0,1)$ be arbitrary numbers. 
Then 
\bd
\item[(a)] 
there exist constants $C_1=C(M,\gamma,\theta)$ and $C_2=C(M,\gamma)$ such that 
\beas\label{210922-8} 
\sup_{f\in\cale(M,\gamma)}\Delta_n(f)
&\leq&
{\colorr{C_1P\Bigl[|\xi_n|>\half\Bigr]^\theta 
+C_2\ep([\gamma+\check{d}]+1,n)
+o(r_n)}}
\eeas
as $n\to\infty$. 

\item[(b)] 
If 
\bea\label{210314-1} 
\sup_n\sup_{(u,v)\in\Lambda^0_n(\check{d},q)}
r_n^{-1}|(u,v)|^{{\coloraka \check{d}+1-\ep}}
|\Phi^{2,\alpha}_n(u,v)|
<\infty
\eea
{\coloro{for every $\alpha\in\bbZ_+^{\check{d}}$}} 
{\coloroy{and some $\ep=\ep(\alpha)\in(0,1)$}}, 
%
then for some constant $C_1=C(M,\gamma,\theta)$, 
\bea\label{210926-5} 
\sup_{f\in\cale(M,\gamma)}\Delta_n(f)
&\leq&
{\colorr{C_1P\Bigl[|\xi_n|>\half\Bigr]^\theta 
+o(r_n)}}
\eea
as $n\to\infty$. \ed
\end{theorem}

In Theorem \ref{210215-1}, ${\mathfrak m}$ is some number, which puts restriction when (\ref{210314-1}) is verified. 
Next, we shall present a version of Theorem \ref{210215-1}. 
{\coloraka Let $s_n$ be a positive random variable on $\Omega$. } 
\bd
\item[[B2$'$\!\!]]$_\ell$ 
{\bf (i)}  
$F_\infty\in\bbD_{\ell+1,\infty}(\bbR^{d_1})$, 
$W_\infty\in\bbD_{{\coloro{\ell+1}},\infty}(\bbR^d)$ 
and $C_\infty\in\bbD_{\ell,\infty}(\bbR^d\otimes\bbR^d)$. 
\bd
\item[(ii)] 
{\colb $M_n\in\bbD_{{\ell+1},\infty}(\bbR^d)$,}
$F_n\in\bbD_{\ell+1,\infty}(\bbR^{d_1})$, 
$W_n\in\bbD_{{\coloro{\ell+1}},\infty}(\bbR^d)$, 
$C_n\in\bbD_{\ell,\infty}(\bbR^d\otimes\bbR^d)$, 
$N_n\in\bbD_{{\coloro{\ell+1}},\infty}(\bbR^d)$  
and 
$s_n\in\bbD_{\ell,\infty}(\bbR)$. 
Moreover, 
\beas 
\sup_{{\coloraka n\in\bbN}}\Big\{ 
\|M_n\|_{\ell+1,p}+\|\dotc_n\|_{\ell,p}+\|\dotw_n\|_{{\coloro{\ell+1}},p}
+\|\dotf_n\|_{\ell+1,p}+\|N_n\|_{{\coloro{\ell+1}},p}+\|s_n\|_{\ell,p} 
\Big\} <\infty
\eeas
for every $p\geq2$. 
\ed
\ed
\bd
\item[[B3$'$\!\!]] 
{\bf (i)}  
{\coloro{$P\l[
\Delta_{(M_n+W_\infty,F_\infty)}< s_n
\r]=O(r_n^{1+\kappa})$ as $n\to\infty$ 
for some $\kappa>0$. }}
\bd
\item[(ii)] For every $p\geq2$, 
\beas 
\limsup_{n\to\infty}E\big[s_n^{-p}\big] 
<\infty
\eeas
and moreover 
$\det C_\infty^{-1}\in\cap_{p\geq2} L^p$. 
\ed
\ed

\begin{en-text}
\bd
{\coloro{\item[[T\!\!]]$_{\ell,{\mathfrak m},{\mathfrak n}}$ }}
{\bf (i)}  
$F_\infty\in\bbD_{\ell+1,\infty}(\bbR^{d_1})$, 
$W_\infty\in\bbD_{{\coloro{\ell+1}},\infty}(\bbR^d)$ 
and $C_\infty\in\bbD_{\ell,\infty}(\bbR^d\otimes\bbR^d)$. 

\bd
\item[(ii)] 
$F_n\in\bbD_{\ell+1,\infty}(\bbR^{d_1})$, 
$W_n\in\bbD_{{\coloro{\ell+1}},\infty}(\bbR^d)$, 
$C_n\in\bbD_{\ell,\infty}(\bbR^d\otimes\bbR^d)$, 
$N_n\in\bbD_{{\coloro{\ell+1}},\infty}(\bbR^d)$  
and a positive functional $s_n\in\bbD_{\ell,\infty}(\bbR)$. 
Moreover, 
\beas 
\sup\Big\{ 
\|M_n\|_{\ell+1,p}+\|\dotc_n\|_{\ell,p}+\|\dotw_n\|_{{\coloro{\ell+1}},p}
+\|\dotf_n\|_{\ell+1,p}+\|N_n\|_{{\coloro{\ell+1}},p}+\|s_n\|_{\ell,p} 
\Big\} <\infty.
\eeas
for every $p\geq2$. 

\item[(iii)] 
{\coloro{$P\l[
\Delta_{(M_n+W_\infty,F_\infty)}< s_n
\r]=O(r_n^{1+\kappa})$ as $n\to\infty$ 
for some $\kappa>0$. }}

\item[(iv)] For every $p\geq2$, 
\beas 
\limsup_{n\to\infty}E\Big[s_n^{-p}\Big] 
<\infty
\eeas
and moreover 
$\Delta_{F_\infty}^{-1}, \>\det C_\infty^{-1}\in\cap_{p\geq2} L^p$. 

\item[(v)] 
$\underline{\sigma}\in\cals(2[(d_1+3)/2],2,1)$. 

\item[(vi)] 
$\bar{\sigma}\in\cals(\ell,{\mathfrak m},{\mathfrak n})$ and 
$\calf^{-1}_{(u,v)}[\tilde{\Phi}^{2,0}](z,x)$ 
has the random symbol $\bar{\sigma}$. 
\ed
\ed
\ \\
\end{en-text}

\halflineskip
\begin{theorem}\label{210926-6}
{\coloro{Let $\ell={\coloraka \ell_0}
\vee(\check{d}+3)$.}}
Suppose that {\rm [B1], [B2$'$]$_\ell$, [B3$'$] and [B4]$_{\ell,{\mathfrak m},{\mathfrak n}}$} are fulfilled. 
Then for any $M,\gamma\in(0,\infty)$, 
\bd
\item[(a)] 
there exists a constant $C$ such that 
\beas\label{210922-8} 
\sup_{f\in\cale(M,\gamma)}\Delta_n(f)
&\leq&
C\ep([\gamma+\check{d}]+1,n)
+o(r_n)
\eeas
as $n\to\infty$. 

\item[(b)] 
If $($\ref{210314-1} $)$ is satisfied, then 
\beas
\sup_{f\in\cale(M,\gamma)}\Delta_n(f)
&=&
o(r_n)
\eeas
as $n\to\infty$. \ed
\end{theorem}
}}

The above results can apply to the expansion of 
the distribution of functions of $(Z_n,F_n)$, 
in particular, $Z_n/F_n$ and $Z_n/\sqrt{F_n}$. 


\section{Proof of Theorems \ref{210215-1} and \ref{210926-6} 
}\label{220915-1}

\subsection{Decomposition of the joint characteristic function}\label{220912-2}

Let $\check{Z}_n=(Z_n,F_n)$. 
Then according to the notation of the multi-index, 
$\check{Z}_n^\alpha=Z_n^{\alpha_1}F_n^{\alpha_2}$ 
for 
$\alpha=(\alpha_1,\alpha_2)\in\bbZ_+^d\times\bbZ_+^{d_1}
=\bbZ_+^{\check{d}}$, $\check{d}=d+d_1$. 
%
%
%
%
Let 
\beas 
\hat{g}_n^\alpha(u,v) 
&=&
E \l[\psi_n\check{Z}_n^\alpha \exp\l(iZ_n[u]+iF_n[v]\r)\r]
\eeas
for $u\in\bbR^d$ and $v\in\bbR^{d_1}$. 
We denote  
\beas 
g_n^\alpha(z,x)
&=&
\frac{1}{(2\pi)^{\check{d}}} \int_{\bbR^{\check{d}}} 
e^{-iu\cdot z -iv\cdot x} \>
\hat{g}_n^\alpha(u,v) \>dudv
\eeas
if the integral on the right-hand side exists. 

\begin{en-text}
\beas 
\hat{h}^\alpha_n(u,v)
&=&
\int_{\bbR^d}
\eeas
\end{en-text}
%
%
%
%
%
\begin{comment} 
{\coloro{Comment. 
The local characteristic function of $(Z_n,F_n)$ is given by 
%
\beas 
\hat{g}_n^\alpha(u,v)
&=&
\partialbs^\alpha E\l[\psi_n\exp\l\{iZ_n[u]+iF_n[v]\r\} \r]
\\&=&
E\l[\psi_n \check{Z}_n^\alpha \exp\l\{iZ_n[u]+iF_n[v]\r\} \r]. 
\eeas
}}
\end{comment}
%
Let
\beas
\dote(x)&=&\int_0^1e^{sx}ds, 
\eeas
then 
$e^x=1+x\dote(x)$. 
For 
\beas
\Psi_n(u,v)
&=&
\exp\l\{ {\colorr iW_n[u]}-\half C_n[u^{\otimes2}]+iF_n[v] \r\}, 
\eeas
we have 
\beas 
\ep_n(u,v)
&=&
\log\l(
\Psi_n(u,v)\Psi_\infty(u,v)^{-1}
\r)
{\colorr +ir_nN_n[u]}
\\&=&
-\half (C_n-C_\infty)[u^{\otimes2}]
+i(F_n[v]-F_\infty[v])
\\&&
{\colorr +i(W_n[u]-W_\infty[u])+ir_nN_n[u]}. 
\eeas
Then we have
\bea\label{200322-1}
\hat{g}_n^\alpha(u,v)
&=&
\partialbs^\alpha E[\Psi_\infty(u,v)\psi_n]
+
\partialbs^\alpha 
E\l[e^n_1(u) \Psi_\infty(u,v) \ep_n(u,v) \dote(\ep_n(u,v))\psi_n \r]
\nn\\&&
+
\partialbs^\alpha E\l[L^n_1(u)\Psi_\infty(u,v)\psi_n\r]
\nn\\&=:&
\Phi^{0,\alpha}_n(u,v)+\Phi^{1,\alpha}_n(u,v)
+\Phi^{2,\alpha}_n(u,v)
.
\eea

\begin{comment}
{\coloro{
\begin{enumerate}
\item Check (\ref{200322-1}). OK
\item Do it for $\partial_{(u,v)}^\alpha$. OK
\end{enumerate}
}}
\end{comment}

Set 
\beas 
B^n_t(u,v) 
&=&
e^n_t(u)\Psi_\infty(u,v). 
\eeas
Then 
\beas 
B^n_t(u,v) &=&
\exp\l(iM^n_t[u]+\half C^n_t[u^{\otimes2}]\r)
\exp\l\{ {\colorr iW_\infty[u]}-\half C_\infty[u^{\otimes2}]+iF_\infty[v] \r\}
\\&=&
\exp\big(iM^n_t[u]\big)
\exp\big(iW_\infty[u]+iF_\infty[v]\big)
\exp\bigg(\half(C^n_t-C_\infty)[u^{\otimes2}]\bigg). 
\eeas
Condition 
[B3] (ii) 
implies 
\bea\label{210908-1} 
|B^n_t(u,v)|
&=&
\exp\bigg(\half(C^n_t-C^n_1)[u^{\otimes2}]\bigg)
\exp\bigg(\half(C_n-C_\infty)[u^{\otimes2}]\bigg)
\nn\\&\leq&
\exp\l(\half r_n^{1-a}|u|^2\r)
\eea
and 
\bea\label{210908-2}
{\mathfrak{Re\>}}\ep_n(u,v)
&\leq&
\half r_n^{1-a}|u|^2
\eea
whenever $\psi_n>0$.  

We shall consider the second-order term of type I. 


\begin{lemma}\label{210825-1}
{\coloro{Suppose that 
$[B1]$ and $[B2]_\ell$ for some $\ell\geq0$ and $[B3]$ $(i)$, $(ii)$ \footnote{It suffices to assume $[R1]$ and $[R2]$ of \cite{Yoshida2010}.} 
are fulfilled. }} Then 
the limit 
\beas 
\tilde{\Phi}^{1,\alpha}(u,v)
&=&
\lim_{n\iku\infty}r_n^{-1}\Phi^{1,\alpha}_n(u,v)
\eeas
exists and takes the form
\beas 
\tilde{\Phi}^{1,\alpha}(u,v)
&=&
\partialbs^\alpha
E \Bigl[ 
\int_{\bbR^d}
\exp\l(iu\cdot z +iF_\infty[v]\r)
\\&&
\>\cdot 
\Bigl( -\half\tilde{C}_\infty(z)[u^{\otimes2}]
+i\>\tilde{W}_\infty(z)[u]
+i\>\tilde{N}_\infty(z)[u]
	+i\>\tilde{F}_\infty(z)[v] \Bigr)
\>\phi(z;W_\infty,C_\infty)\>dz
\Bigr]
\eeas
for every $\alpha\in\bbZ_+^{\check{d}}$. 
\end{lemma}
\proof 
Since the family of functions $u\mapsto B^n_1(u,v)\psi_n$ $(n\in\bbN$) 
is uniformly (in $n$) locally (in $u$) bounded 
due to (\ref{210908-1}), we have, by using the weak convergence condition 
and  (\ref{210908-2}), 
\beas 
\tilde{\Phi}^{1,\alpha}(u,v)
&=&
\lim_{n\iku\infty}r_n^{-1}\Phi^{1,\alpha}_n(u,v)
\\&=&
\lim_{n\iku\infty}r_n^{-1}
\partialbs^\alpha
E\bigg[
\exp\big(iM^n_1[u]\big)
\exp\big(iW_\infty[u]+iF_\infty[v]\big)
\ep_n(u,v)
\exp\bigg(\half(C_n-C_\infty)[u^{\otimes2}]\bigg)
\dote(\ep_n(u,v))\psi_n
\bigg]
\\&=&
\lim_{n\iku\infty}
E\bigg[
\partialbs^\alpha\bigg\{
\exp\big(iM^n_1[u]\big)
\exp\big(iW_\infty[u]+iF_\infty[v]\big)
r_n^{-1}\ep_n(u,v)
\bigg\}
\exp\bigg(\half(C_n-C_\infty)[u^{\otimes2}]\bigg)
\dote(\ep_n(u,v))\psi_n
\bigg]
\\&=&
\partialbs^\alpha
E \Bigl[ 
\exp\Bigl(iM^\infty_1[u] +iW_\infty[u]+iF_\infty[v]\Bigr)
\\&&
\>\cdot 
\Bigl( -\half\dotc_\infty[u^{\otimes2}]+i\dotw_\infty[u]+iN_\infty[u]
+i\dotf_\infty[v] \Bigr)
\Bigr]. 
\eeas
The last expression is equal to
\beas 
&&
\partialbs^\alpha
E \Bigl[ 
\exp\Bigl(iM^\infty_1[u] +iW_\infty[u]+iF_\infty[v]\Bigr)
\\&&
\>\cdot 
\Bigl( -\half\bbE_{\check{\calf}}[\dotc_\infty][u^{\otimes2}]
+i\>\bbE_{\check{\calf}}[\dotw_\infty][u]
+i\>\bbE_{\check{\calf}}[N_\infty][u]
+i\>\bbE_{\check{\calf}}[\dotf_\infty][v] \Bigr)
\Bigr]
\\&=&
\partialbs^\alpha
E \Bigl[ 
\exp\l(iM_\infty[u] +iW_\infty[u]+iF_\infty[v]\r)
\\&&
\>\cdot 
\Bigl( -\half\check{C}_\infty(M_\infty)[u^{\otimes2}]
+i\>\check{W}_\infty(M_\infty)[u]
+i\>\check{N}_\infty(M_\infty)[u]
+i\>\check{F}_\infty(M_\infty)[v] \Bigr)
\Bigr]
\\&=&
\partialbs^\alpha
E \Bigl[ 
\int_{\bbR^d}
\exp\l(iu\cdot z +iW_\infty[u]+iF_\infty[v]\r)
\\&&
\>\cdot 
\Bigl( -\half\check{C}_\infty(z)[u^{\otimes2}]
+i\>\check{W}_\infty(z)[u]
+i\>\check{N}_\infty(z)[u]
+i\>\check{F}_\infty(z)[v] \Bigr)
\>\phi(z;0,C_\infty)\>dz
\Bigr]
\\&=&
\partialbs^\alpha
E \Bigl[ 
\int_{\bbR^d}
\exp\l(iu\cdot z +iF_\infty[v]\r)
\\&&
\>\cdot 
\Bigl( -\half\tilde{C}_\infty(z)[u^{\otimes2}]
+i\>\tilde{W}_\infty(z)[u]
+i\>\tilde{N}_\infty(z)[u]
+i\>\tilde{F}_\infty(z)[v] \Bigr)
\>\phi(z;W_\infty,C_\infty)\>dz
\Bigr]
\eeas
by $\calf$-conditional expectation, which is what we desired. 
\qed 

\ \\

In particular, 
\bea\label{221003-3}
\tilde{\Phi}^{1,\alpha}(u,v)
&=& \partialbs^\alpha 
\tilde{\Phi}^{1,0}(u,v). 
\eea

\begin{comment}
{\coloro{Comment. 
\noindent
$u$-partial Fourier inversion of $\tilde{\Phi}^{1,0}(u,v)$: 
\beas 
\calf^{-1}_u[\tilde{\Phi}^{1,0}](z)
&=&
E\Bigl[\Bigl\{
\half \partial_{z_j}\partial_{z_k} 
\l(\tilde{C}_\infty(z)^{jk} \phi(z;W_\infty,C_\infty)\r)
\\&&
-
\partial_{z_j}\l(\tilde{W}_\infty(z)^j \phi(z;W_\infty,C_\infty)\r)
\\&&
-
\partial_{z_j}\l(\tilde{N}_\infty(z)^j \phi(z;W_\infty,C_\infty)\r)
\\&&
+i\>\tilde{F}_\infty(z)[v]  \phi(z;W_\infty,C_\infty)
\Bigr\} \exp\l(iF_\infty[v]\r)
\Bigr];
\eeas
moreover, $v$-Fourier inversion of this: 
\beas 
\calf^{-1}_{(u,v)}[\tilde{\Phi}^{1,0}](z,x)
&=&
\half \partial_{z_j}\partial_{z_k} 
\l(E[\tilde{C}_\infty(z)^{j,k}\phi(z;W_\infty,C_\infty)|F_\infty=x] \r)
p^{F_\infty}(x)
\\&&
-
\partial_{z_j}\l(E[\tilde{W}_\infty(z)^j 
\phi(z;W_\infty,C_\infty)|F_\infty=x]\r)p^{F_\infty}(x)
\\&&
-
\partial_{z_j}\l(E[\tilde{N}_\infty(z)^j 
\phi(z;W_\infty,C_\infty)|F_\infty=x]\r)p^{F_\infty}(x)
\\&&
-\partial_{x_l}\l(E[\tilde{F}_\infty(z)^l  
\phi(z;W_\infty,C_\infty)|F_\infty=x] p^{F_\infty}(x)\r).
\eeas
}}
\end{comment}

\begin{en-text}
\beas 
\half\partial_z^2 
\Bigl\{ 
\bbE_\calf\l[\xi|\zeta=z\r] 
\phi(z;W_\infty,C_\infty) \Bigr\}
-
\partial?
\eeas
It is better to consider the joint density? 
\end{en-text}

\begin{comment}
{\coloro{
\begin{enumerate}
\item When $\Psi_\infty$ is deterministic, $\tilde{\Phi}^{2,0}$-term 
is zero, and the asymptotic expansion lives in this world. 

\item Even when $\Psi_\infty$ is deterministic, we still obtain 
conditional asymptotic expansion. 

\item Integrability should be mentioned. 
\end{enumerate}
\end{comment}

{\colb Under [B4]$_{\coloraka \ell,{\mathfrak m},{\mathfrak n}}$ (i),} the conditions just before (\ref{210815-2}) 
are satisfied for $\varsigma=\underline{\sigma}${\coloraka , and}   
%
\beas 
\calf^{-1}_{(u,v)}[\tilde{\Phi}^{1,0}](z,x)
&=&
\half \partial_{z_j}\partial_{z_k} 
\l(E[\tilde{C}_\infty(z)^{j,k}\phi(z;W_\infty,C_\infty)|F_\infty=x] \r)
p^{F_\infty}(x)
\\&&
-
\partial_{z_j}\l(E[\tilde{W}_\infty(z)^j 
\phi(z;W_\infty,C_\infty)|F_\infty=x]\r)p^{F_\infty}(x)
\\&&
-
\partial_{z_j}\l(E[\tilde{N}_\infty(z)^j 
\phi(z;W_\infty,C_\infty)|F_\infty=x]\r)p^{F_\infty}(x)
\\&&
-\partial_{x_l}\l(E[\tilde{F}_\infty(z)^l  
\phi(z;W_\infty,C_\infty)|F_\infty=x] p^{F_\infty}(x)\r)
\\&=&
E\bigg[\underline{\sigma}(z,\partial_z,\partial_x)^*
\bigg\{
\phi(z;W_\infty,C_\infty)\bigg|F_\infty=x\bigg]
p^{F_\infty}(x)
\bigg\}.
\eeas

We can say that the random symbol $\underline{\sigma}$ 
corresponds to the case where 
the martingale central limit theorem occurs 
without conditioning. 
On the other hand, we need another random symbol 
which reflects the deviation of the martingale  
in question from a {\it real martingale} under conditioning. 
%
The Fourier inversion $\calf^{-1}_{(u,v)}[\tilde{\Phi}^{2,0}]$ has a random symbol 
$\bar{\sigma}(\partial_z,\partial_x)$ in that 
\beas 
\calf^{-1}_{(u,v)}[\tilde{\Phi}^{2,0}](z,x)
&=&
E\bigg[\bar{\sigma}(\partial_z,\partial_x)^*
\bigg\{
\phi(z;W_\infty,C_\infty)\bigg|F_\infty=x\bigg]
p^{F_\infty}(x)
\bigg\}. \hspace{5cm}
\eeas
See (\ref{210922-5}) and (\ref{220131-1}) for this equality. 
}}
\begin{en-text}
Let $,{\mathfrak m}\in\bbZ_+$ and ${\mathfrak n}\in\bbN$. For the random symbols 
$\underline{\sigma}(z,iu,iv)$ and 
$\bar{\sigma}(\partial_z,\partial_x)$, we suppose 
\bd
{\coloro{
\item[[R\!\!]]$_{\ell,{\mathfrak m},{\mathfrak n}}$
$\ell\geq2[({\mathfrak n}+d_1+2)/2]$ and}}
\bd
\item[(i)] 
$F_\infty\in\bbD_{\ell+1,\infty}(\bbR^{d_1})$, 
$W_\infty\in\bbD_{\ell,\infty}(\bbR^d)$ and 
$C_\infty\in\bbD_{\ell,\infty}(\bbR^d\otimes\bbR^d)$. 

\item[(ii)] 
$\Delta_{F_\infty}^{-1},\>\det C_\infty^{-1}\in\cap_{p\geq2} L^p$. 

\item[(iii)] 
$\underline{\sigma}\in\cals(2[(d_1+3)/2],2,1)$. 

\item[(iv)] 
$\bar{\sigma}\in\cals(\ell,{\mathfrak m},{\mathfrak n})$. 

\item[(v)] 
$\calf^{-1}_{(u,v)}[\tilde{\Phi}^{2,0}]$ has a random symbol 
$\bar{\sigma}(\partial_z,\partial_x)$. 
\ed
\ed
\ \\
In the above, 
$\bar{\sigma}(\partial_z,\partial_x)^*$ is the adjoint operator 
of $\bar{\sigma}(\partial_z,\partial_x)$. 
\end{en-text}
\begin{en-text}
We need a sequence of functionals $\xi_n$ 
to avoid the difficulties in the degeneracy. 
Let $\psi\in C^\infty(\bbR;[0,1])$ such that $\psi(x)=1$ if $|x|\leq1/2$ 
and $\psi(x)=0$ if $|x|\geq1$. 
Let $\psi_n=\psi(\xi_n)$. 
\end{en-text}

\subsection{Estimates of error bounds}

{\coloro{
\begin{comment} 
Comment. 
To balance the treatments of two types of higher order terms, 
it seems better to write the first type in a more abstract way and 
put them before the "Second-order term of type I". 
%
Otherwise, the structure of $M^n$ should be restricted 
from the beginning as the double stochastic integral. 
\end{comment}
}}

{\colr{
We shall investigate the approximation error of $p_n(z,x)$ to 
the joint local density of $(Z_n,F_n)$. 
The error bounds depend on the smoothness of the 
distribution $\call\{(Z_n,F_n)\}$, and the arguments 
will require a tool to evaluate it quantitatively. 
The conditions are written in terms of the Malliavin calculus 
since it is convenience in practical uses while 
more primitive expression of them would be possible without it 
if we would admit more cumbersome descriptions. 
}}

\begin{en-text}
Let $m\in\bbZ_+$ and $n\in\bbN$.  
We consider the following conditions. 
\bd
{\coloro{
\item[[S\!\!]]$_{\ell,{\mathfrak m},{\mathfrak n}}$
}}
{\bf (i)}  
$F_\infty\in\bbD_{\ell+1,\infty}(\bbR^{d_1})$, 
$W_\infty\in\bbD_{{\coloro{\ell+1}},\infty}(\bbR^d)$ 
and $C_\infty\in\bbD_{\ell,\infty}(\bbR^d\otimes\bbR^d)$. 

\bd
\item[(ii)] 
$F_n\in\bbD_{\ell+1,\infty}(\bbR^{d_1})$, 
$W_n\in\bbD_{{\coloro{\ell+1}},\infty}(\bbR^d)$, 
$C_n\in\bbD_{\ell,\infty}(\bbR^d\otimes\bbR^d)$, 
$N_n\in\bbD_{{\coloro{\ell+1}},\infty}(\bbR^d)$  
and $\xi_n\in\bbD_{\ell,\infty}(\bbR)$. 
Moreover, 
\beas 
\sup\Big\{ 
\|M_n\|_{\ell+1,p}+\|\dotc_n\|_{\ell,p}+\|\dotw_n\|_{{\coloro{\ell+1}},p}
+\|\dotf_n\|_{\ell+1,p}+\|N_n\|_{{\coloro{\ell+1}},p}+\|\xi_n\|_{\ell,p} 
\Big\} <\infty.
\eeas
for every $p\geq2$. 

\item[(iii)] 
{\coloro{$\lim_{n\to\infty}P\l[|\xi_n|\leq\half\r]=1$. }}

\item[(iv)] 
$|C_n-C_\infty|>r_n^{1-a}$ implies $|\xi_n|\geq 1$, 
where $a\in(0,1/3)$ is a constant.

\item[(v)] For every $p\geq2$, 
\beas 
\limsup_{n\to\infty}E\Big[1_{\{|\xi_n|\leq1\}} \Delta_{(M_n+W_\infty,F_\infty)}^{-p}\Big] 
<\infty
\eeas
and moreover 
$\Delta_{F_\infty}^{-1}, \>\det C_\infty^{-1}\in\cap_{p\geq2} L^p$.

\item[(vi)] 
$\underline{\sigma}\in\cals(2[(d_1+3)/2],2,1)$. 

\item[(vii)] 
$\bar{\sigma}\in\cals(\ell,{\mathfrak m},{\mathfrak n})$ and 
$\calf^{-1}_{(u,v)}[\tilde{\Phi}^{2,0}](z,x)$ 
has the random symbol $\bar{\sigma}$. 
\ed
\ed
\ \\
\end{en-text}

\begin{comment}
{\coloro{Comment. 
The nondegeneracy of $\det\sigma_{(M_n+W_\infty,F_\infty)}$  
{\colorr{implies}} the nondegeneracy of $F_\infty$.}}
{\colorr{Because
\beas 
1/|A|=|D-.....|/|M|\leq |D|/|M|.
\eeas
and $D$ is $O(1)$. 
So we do not need to assume $\Delta_{F_\infty}^{-1}\in\cap_{p\geq2} L^p$. 
}}
\end{comment}

\colr{
Let 
\beas 
h^0_n(z,x) 
&=& 
E\bigg[
\psi_n
\phi(z;W_\infty,C_\infty)\bigg|F_\infty=x\bigg]
p^{F_\infty}(x)
\\&&
+
r_n 
E\bigg[\sigma(z,\partial_z,\partial_x)^*
\bigg\{
\phi(z;W_\infty,C_\infty)\bigg|F_\infty=x\bigg]
p^{F_\infty}(x)
\bigg\}. 
\eeas
and let 
\beas 
h^\alpha_n(z,x) = (z,x)^\alpha h^0_n(z,x). 
\eeas 

}

{\colorsb{

\begin{comment}
{\coloro{
ASS. 
\begin{itemize}
\item 
$\psi_n=0$ whenever 
($|C_\infty|> \half r_n^{(1-a)}$ or 
$|C_n-C_\infty|>\half r_n^{(1-a)}$)

\item 
$\dotc_n$, $\dotw_n$, $\dotf_n$ and $N_n$ $(n\in\bbN)$ 
are bounded in $\bbD_{\ell,\infty-}$. 
($\ell$ is an index of differentiability we will later specify. 
At least, $\ell>d+d_1+2$; however this bound will be replaced by 
a larger value later. )

\item 
$(M_n+W_\infty,F_\infty)$ $(n\in\bbN)$ is uniformly nondegenerate 
under truncation by $\psi_n$. 
Moreover, 
$M_n,W_\infty,F_\infty$ $(n\in\bbN)$ are bounded in $\bbD_{\ell+1,\infty-}$. 
\end{itemize}
}}
\end{comment}

%
\begin{en-text}
{\colr{
Let $\Lambda^0_n(d,q)=\{u\in\bbR^d;|u|\leq r_n^{-q}\}$,  
where 
$q=(1-a)/2\in(1/3,1/2)$. 
\end{en-text}

\begin{lemma}\label{210311-2}  
Let $\ell\geq2$. 
Suppose that 
$[B2]_\ell$ and $[B3]$ 
are fulfilled. Then  
\beas 
\sup_n \sup_{u\in\Lambda^0_n(d,q)}
\sup_{v\in\Lambda^0_n(d_1,2q) }
r_n^{-1}|(u,v)|^{\ell-2}|\Phi^{1,\alpha}_n(u,v)|
<\infty.
\eeas
Moreover, 
\beas 
\sup_{u\in\bbR^d}
\sup_{v\in\bbR^{d_1}}
|(u,v)|^{\ell-2} |\tilde{\Phi}^{1,\alpha}(u,v)|
<\infty.
\eeas
\end{lemma}
\proof 
By definition, 
\bea\label{210320-1} 
r_n^{-1}\Phi^{1,\alpha}_n(u,v)
&=&
\partialbs^\alpha
E\bigg[
\exp\bigg(i(M^n_1+W_\infty)[u]+iF_\infty[v]\bigg)
\nn\\&&\cdot 
r_n^{-1}\ep_n(u,v)
\exp\bigg(\half(C_n-C_\infty)[u^{\otimes2}]\bigg)
\dote(\ep_n(u,v))\psi_n
\bigg]
\nn\\&=&
E\bigg[
\exp\bigg(i(M^n_1+W_\infty)[u]+iF_\infty[v]\bigg)
\nn\\&&\cdot 
p_{2,1,\alpha}
\bigg(u,v;r_n,
{\colb M^n_1,W_\infty,F_\infty, }
N_n,\dotc_n,\dotw_n,\dotf_n\bigg)
\psi_n
\bigg],
\eea
where $p_{2,1,\alpha}$ is a smooth function such that 
\beas
|p_{2,1,\alpha}(u,v;{\sf x})| 
&\leq& 
C_\alpha(1+|u|^2+|v|) (1+|{\sf x}|)^{C_\alpha}
\sskip 
(u\in\bbR^d,v\in\bbR^{d_1}, \>{\sf x}\in \bbR^{{\colb d^2+4d+2d_1+1}})
\eeas
for some constant $C_\alpha$, and moreover 
the derivative of $p_{2,1,\alpha}$ of any order 
admits the same type estimate. 
It should be noted that in the expectation 
on the right-hand side of (\ref{210320-1}), 
we do not need an exponential like factor 
in $p_{2,1,\alpha}$, 
due to the truncation by $\psi_n$.

\begin{en-text}
\bea\label{210311-1} 
r_n^{-1}\Phi^{1,\alpha}_n(u,v)
&=&
r_n^{-1}
\partialbs^\alpha 
E\l[e^n_1(u) \Psi_\infty(u,v) \ep_n(u,v) \dote(\ep_n(u,v))\psi_n \r]
\nn\\&=&
\partialbs^\alpha 
E\bigg[
\exp\big(iM^n_1[u]+iW_\infty[u]+iF_\infty[v]\big)
\exp\big(\half (C^n_1-C_\infty)[u^{\otimes2}]\big)
\nn\\&&
\cdot 
r_n^{-1}\ep_n(u,v)
\dote(\ep_n(u,v))\psi_n 
\bigg]
\eea
\end{en-text}
We will apply the integration-by-parts formula $\ell$-times for the pull-back 
of the function 
$f(z,x) = e^{iu\cdot z+iv\cdot x}$ by 
taking the advantage of 
the uniform nondegeneracy of the functional 
{\coloro{$(M_n+W_\infty,F_\infty)$}} under $\psi_n$. 
Then, by truncation by $\psi_n$ and the restriction of 
the region of $(u,v)$, 
the functional in the expectation of (\ref{210320-1}) 
is essentially quadratic in $u$ and linear in $v$. 
It should be noted that all variables related to $u$ or $v$ 
are differentiated at the same time 
when applying the IBP-formula, therefore the index of differentiability 
should be common. 

The second inequality follows from the first one 
if one takes the limit in $n$. 
\qed\\
}}

\ \\
\colr{
Define the $\check{d}\times\check{d}$ random matrix $R_n'$ by
\beas 
R_n'
&=&
\sigma_{Q_n}^{-1}\bigl(r_n\langle DQ_n,DR_n\rangle 
+ r_n\langle DR_n,DQ_n\rangle +r_n^2\langle DR_n,DR_n\rangle\bigr), 
\eeas
where 
$Q_n=(M_n+W_\infty,F_\infty)$ and $R_n=(\dotw_n+N_n,\dotf_n)$. 
Obviously 
\bea\label{2109111212-1}
\sigma_{(Z_n,F_n)} &=& \sigma_{Q_n}(I_{\check{d}}+R_n'). 
\eea
Let $\xi_n'={\colb r_n^{-1}}|R_n'|^2$. 
We redefine $\psi_n$ by 
\bea\label{210922-9} 
\psi_n &=& \psi(\xi_n)\psi(\xi_n'). 
\eea
$\Phi^{j,\alpha}_n$ ($j=0,1,2$), 
$g^\alpha_n$ and $h^\alpha_n$ will be defined for 
$\psi_n$ given in (\ref{210922-9}). 
The following is an extension of the expansion of the local density 
in the normal limit case; see Theorems 3 and 5 of \cite{Yoshida1997}. 

\begin{lemma}\label{210215-2}
{\coloro{Let $\ell={\coloraka \ell_0}
\vee(\check{d}+3)$.}}
Suppose that 
$[B1]$, $[B2]_\ell$, $[B3]$ and $[B4]_{\ell,{\mathfrak m},{\mathfrak n}}$ 
are satisfied. Then 
for every $k\in\bbZ_+$, 
\bea\label{210314-5} 
\sup_{(z,x)\in\bbR^{d+d_1}}
\big||{\colorr{(z,x)}}|^k
\big(g_n^0(z,x)-h^0_n(z,x)\big)\big|
\leq
{\colorr{\ep(k,n)+o(r_n)}}. 
\eea
{\colorr{Furthermore, if 
(\ref{210314-1}) is satisfied 
{\coloro{for every $\alpha\in\bbZ_+^{\check{d}}$}} 
{\coloroy{and some $\ep=\ep(\alpha)\in(0,1)$}}, 
then $\ep(k,n)=o(r_n)$ {\coloro{for every $k\in\bbZ_+$}}. 
}}
\end{lemma}
}
\proof 
By definition, 
\beas 
\hat{h}^\alpha_n(u,v)
&=&
\partialbs^\alpha
E\big[\Psi_\infty(u,v)\psi_n\big] 
+r_n\partialbs^\alpha\tilde{\Phi}^{1,0}(u,v)
+r_n\partialbs^\alpha\tilde{\Phi}^{2,0}(u,v). 
\eeas 
%
%
{\colb Moreover,} we have (\ref{200322-1}) for $\hat{g}_n^\alpha(u,v)$. 
%
\begin{comment}
\\ 
{\coloro{Comment. 
\beas
\hat{g}_n^\alpha(u,v)
&=&
\partialbs^\alpha E[\Psi_\infty(u,v)\psi_n]
+
\partialbs^\alpha 
E\l[e^n_1(u) \Psi_\infty(u,v) \ep_n(u,v) \dote(\ep_n(u,v))\psi_n \r]
\nn\\&&
+
\partialbs^\alpha E\l[L^n_1(u)\Psi_\infty(u,v)\psi_n\r]
\nn\\&=&
\Phi^{0,\alpha}_n(u,v)+\Phi^{1,\alpha}_n(u,v)
+\Phi^{2,\alpha}_n(u,v).
\eeas
}}
\end{comment}
%
%
%
Applying the IBP formula $\ell$-times under truncation by $\psi_n$ 
{\colb to the definition of $\hat{g}^\alpha_n$}, 
we see that for every $\alpha\in\bbZ_+^{\check{d}}$, 
there exist $C>0$ and $n^*\in\bbN$ such that 
\bea\label{210922-6} 
|\hat{g}^\alpha_n(u,v)| &\leq& 
C (1+|u|+|v|)^{-\ell}
\eea
for all $(u,v)\in\bbR^{\check{d}}$ and $n\geq n^*$. 
In what follows, we will consider only sufficiently large $n$. 
Therefore the local densities $g^\alpha_n(z,x)$ are well defined. 
On the other hand, one can verify the integrability of 
$\tilde{\Phi}^{2,\alpha}(u,v)$ by using the nondegeneracy of $C_\infty$,  
the nondegeneracy of the Malliavin covariance of 
$F_\infty$ in 
[B3] 
and 
the representation {\coloraka (\ref{220131-1})}. 
{\coloraka [The aftereffect by the integration-by-parts for $F_\infty$ is absorbed by the nondegeneracy of $C_\infty$.] }
The integrability of $\Phi^{1,\alpha}_n(u,v)$ and $\tilde{\Phi}^{1,\alpha}(u,v)$ has been obtained 
in Lemma \ref{210311-2} . 
The decomposition (\ref{200322-1}) implies that 
$\Phi^{2,\alpha}_n$ is integrable for each $n$. 
Thus, the decomposition and the estimate of the following inequality 
based on 
the Fourier inversion formula are valid: 
\bea\label{210922-7}
\sup_{(z,x)\in\bbR^{\check{d}}}
\big|g^\alpha_n(z,x)-h^\alpha_n(z,x)\big|
&=&
\sup_{(z,x)\in\bbR^{\check{d}}}
\frac{1}{(2\pi)^{\check{d}}}
\bigg|\int_{\bbR^{\check{d}}}e^{-iu\cdot z-iv\cdot x}
\big(\hat{g}^\alpha_n(u,v)-\hat{h}^\alpha_n(u,v)\big)dudv\bigg|
\nn\\&\leq&
\frac{1}{(2\pi)^{\check{d}}}
\int_{\bbR^{\check{d}}\setminus\Lambda^0_n(\check{d},q)
}
\big(|\hat{g}^\alpha_n(u,v)|+|\hat{h}^\alpha_n(u,v)|\big)dudv
\nn\\&&
+\frac{1}{(2\pi)^{\check{d}}}
r_n
\int_{\Lambda^0_n(\check{d},q)
}
\big|r_n^{-1}\Phi^{1,\alpha}_n(u,v)-\tilde{\Phi}^{1,\alpha}(u,v)\big|dudv
\nn\\&&
+\frac{1}{(2\pi)^{\check{d}}}
r_n
\int_{\Lambda^0_n(\check{d},q)
}
\big|r_n^{-1}\Phi^{2,\alpha}_n(u,v)-\tilde{\Phi}^{2,\alpha}(u,v)\big|dudv. 
\eea
%
\begin{comment}
{\coloro{
Comment. 
In the classical setting, the last term was estimated 
by comparing it with $0$ (the expectation of a martingale), which 
gave $\|1-\psi_n\|_p$. 
}}
\end{comment}
%
We remark that $\tilde{\Phi}^{2,\alpha}=\partialbs^\alpha\tilde{\Phi}^{2,0}$ from (\ref{220131-1}) 
as well as $\tilde{\Phi}^{1,\alpha}=\partialbs^\alpha\tilde{\Phi}^{1,0}$ mentioned at (\ref{221003-3}). 

By (\ref{210922-6}), 
\beas 
\int_{\bbR^{\check{d}}\setminus\Lambda^0_n(\check{d},q)}
|\hat{g}^\alpha_n(u,v)|dudv
&=&
O(r_n^{q(\ell-\check{d})})
=O(r_n^{3q})=o(r_n)
\eeas
since $3q>1$ 
due to {\colorr{$a<1/3$}}.

As suggested above, we use the Gaussianity in $u$ and the IBP formula in $v$ 
to obtain
\beas 
\int_{\bbR^{\check{d}}\setminus\Lambda^0_n(\check{d},q)}
|\hat{h}^\alpha_n(u,v)|dudv
&\simleq&
\int_{\bbR^{\check{d}}\setminus\Lambda^0_n(\check{d},q)}|(u,v)|^{-\check{d}-3}dudv
\\&&
+r_n\int_{\bbR^{\check{d}}\setminus\Lambda^0_n(\check{d},q)}
|(u,v)|^{-\check{d}-3+2}dudv
\\&&
+r_n\int_{\bbR^{\check{d}}\setminus\Lambda^0_n(\check{d},q)}
(1+|u|)^{-k}
|v|^{-\{({\mathfrak n}+d_1+1)\}+{\mathfrak n}}dudv
\\&=&
O(r_n^{3q})+o(r_n)=o(r_n),
\eeas 
where $k\geq d+1$ is an arbitrary number. 
%
\begin{comment}
\\ \noindent
{\coloro{Comments. 
We can use the IBP formula in the direction of 
${\rm argmax}\{|u_i|,|v_j|\}$ to get $|(u,v)|$. 
The estimate in the direction of $u$ is no problem by using the Gaussianity 
but the IBP formula. \\
}}
\end{comment}
%
The second term on the right-hand side of (\ref{210922-7}) 
is estimated by Lemma \ref{210311-2}. 

{\colr{In order to conclude the second assertion of the lemma, 
to the third term on the right-hand side of (\ref{210922-7}), obviously, 
we apply (\ref{210314-1}) and the inequality which is obtained by its limit as $n\to\infty$. 
\qed\\

{\it Proof of Theorem \ref{210215-1}.} 
\begin{en-text}
Under the nondegeneracy of $(Z_n,F_n)$, we have 
\beas 
\sup_{(u,v)\in\bbR^{d+d_1}} 
|(u,v)|^j
\big| 
\E \l[\psi_n\check{Z}_n^\alpha \exp\l(iZ_n[u]+ivF_n[v]\r)\r]\big|
<\infty\sskip(\alpha\in\bbZ_+^{d+d_1},\ j\in\bbZ_+). 
\eeas
Applying this inequality to the case $j>d+d_1$, 
\end{en-text}
Let $f\in\cale(M,\gamma)$. 
By the estimate (\ref{210922-6}), 
we {\colb know} that 
the local density $g^0_n$ exists. 
In fact $g^0_n$ is a continuous version of 
$\{E[\psi_n|\check{Z}_n=(z,x)]dP^{\check{Z}_n}\}/dzdx$, 
it has moments of any order and 
\beas 
E[f(\check{Z}_n)\psi_n] 
&=&
\int_{\bbR^{\check{d}}}f(z,x)g_n^0(z,x)dzdx. 
\eeas
In particular, the integrability of $f(z,x)g_n^0(z,x)$ is obvious. 
We have 
\beas 
\big|E[f(\check{Z}_n)]-E[f(\check{Z}_n)\psi_n]\big|
&\leq&
\|f(\check{Z}_n)\|_{L^{p'}}\|1-\psi_n\|_{L^p}
\eeas
for any dual pair of positive numbers $(p,p')$ with $1/p+1/p'=1$.  
The integrability of 
$f(z,x)h^0_n(z,x)$ is known from (\ref{210923-1}) in Section \ref{210909-2} 
or in particular from Lemma \ref{210215-2}, 
and also we have 
\beas &&
\big|
\int_{\bbR^{\check{d}}}f(z,x)g_n^0(z,x)dzdx
- 
\int_{\bbR^{\check{d}}}f(z,x)h_n^0(z,x)dzdx
\big|
\\&\leq&
\int_{\bbR^{\check{d}}}|f(z,x)|(1+|(z,x)|^2)^{-k/2}dzdx
\cdot 
\sup_{(z,x)\in\bbR^{\check{d}}}
\big|(1+|(z,x)|^2)^{k/2}
\big(g_n^0(z,x)-h_n^0(z,x)\big)\big| 
\eeas
for any $k>\gamma+\check{d}$. 
Apply a similar estimate as in (\ref{210923-1}) for 
$c_j=1-\psi_n$ and $p{\coloraka \in(1,1/\theta)}$ to obtain 
\beas 
&&
\Big| \int_{\bbR^{\check{d}}}f(z,x)
E\Big[ (1-\psi_n) \phi(z;W_\infty,C_\infty)\delta_x(F_\infty) \Big] \>dzdx
\Big| 
\\&\leq&  
{\sf C}\|1-\psi_n\|_{\ell,p}, 
\eeas
which serves to replace 
$p_n$ by $h^0_n$. 
Furthermore 
\beas 
\|1-\psi_n\|_{L^p}
\leq 
\|1-\psi_n\|_{\ell,p} 
&\leq& 
{\sf C_p}
\l(P\l[|\xi_n|> \half\r]^\theta+P\l[|\xi_n'|> \half\r]^\theta\r)
\eeas
for $\theta<1/p$. 
Here in order to obtain $P[|\xi_n'|>2^{-1}]=o(r_n^k)$ for any $k>0$, 
we can apply the truncation by $\xi_n$ for nondegeneracy of $\sigma_{Q_n}$. 
After all, we obtain the desired result by Lemma \ref{210215-2}. 
\begin{en-text}
\beas 
\sup_{(z,x)\in\bbR^{\check{d}}} |(z,x)|^k
\big|p_n(z,x)-h^0_n(z,x)\big|
&\leq&
\sup_{(z,x)\in\bbR^{\check{d}}} |(z,x)|^k
\big| 
\eeas
\end{en-text}
\begin{en-text}
\beas 
\big|\calf^{-1}_{(u,v)}
\partialbs^\alpha
E\big[\Psi_\infty(u,v)(1-\psi_n)\big]\big|
&\leq&
\sup_{(z,x)\in\bbR^{\check{d}}}
\|\calf^{-1}_{(u,v)}[
\partialbs^\alpha\Psi_\infty](z,x)\|_{p'}
\|1-\psi_n\|_p
\eeas 
; here 
the uniform local nondegeneracy of $(M_n,W_n,F_n)$ 
is working to validate the existence of $\calf^{-1}$, and 
the supremum in $(z,x)$ comes from ``$L^p$''-boundedness of 
variables. 
\end{en-text}
\qed\vspace{5mm}

{\it Proof of Theorem \ref{210926-6}.} 
It holds that 
$P\l[r_n^{-c}|C_n-C_\infty|>1\r]=O(r_n^{1+\kappa'})$ 
for any constant $c\in(2/3,{\coloro{1}})$ and {\coloro{any}} $\kappa'>0$, 
{\coloro{due to the boundedness of $\{\dotc_n\}$ in $L^{\infty-}{\coloraka=\cap_{p>1}L^p}$.}} 
%
We take constants $a$ and $c$ so that $2/3<1-a<c<1$. 
Let 
\beas 
\xi_n 
&=&
10^{-1}r_n^{-2c}|C_n-C_\infty|^2
+
2\Big[1+4\Delta_{(M_n+W_\infty,F_\infty)}s_n^{-1}\Big]^{-1}. 
\eeas
We only consider sufficiently large $n$. 
Then {\colb [B2]$_{\ell}$ and [B3] 
are} verified 
under {\colb [B2$'$]$_{\ell}$ and [B3$'$]}.
We take a sufficiently large $\theta\in(0,1)$ depending on $\kappa\wedge\kappa'$, 
and apply Theorem \ref{210215-1}. 
Note that 
\beas 
\Bigl\{|\xi_n|>\half\Bigr\} &\subset& 
\Bigl\{r_n^{-c}|C_n-C_\infty|>1\Bigr\}
\bigcup
\Bigl\{\Delta_{(M_n+W_\infty,F_\infty)}<s_n\Bigr\}.
\eeas
\qed\\

}}

\section{Expansion for the double stochastic integral}\label{220915-2} 

\begin{en-text}
{\coloro{To be removed. 
Define a $d$-dimensional process $M^n$ by 
\beas 
M^n_t&=& D^*(1_{[0,t]}D^*K^n), 
\eeas
where 
$K^n\in\bbD_{l,\infty-}(H\otimes H\otimes\bbR^d)$. 
Denote by $K^n(s,r)$ the $\bbR^d$-valued density function representing 
$K^n$ as a mapping in 
$L^2([0,1]^2,dsdr;\bbR^{{\sf r}}\otimes\bbR^{{\sf r}}\otimes\bbR^d)$. 
Also we denote by $K^n(s)$ the mapping 
$[0,1]\ni r\mapsto K^n(s,r)\in\bbR^{{\sf r}}\otimes\bbR^{{\sf r}}\otimes\bbR^d$. 

We assume that 
the processes $r\mapsto K^n(s,r)$ and $s\mapsto D^*K^n(s)$ 
are progressively measurable. 
Thus, 
\beas 
M^n_t
&=& 
\int_0^t D^*K^n(s) dw_s
=
\int_0^t \l(\int_0^1 K^n(s,r)dw_r\r)dw_s 
\eeas
is a double It\^o integral, where the element $w$ of $\bbW$ is also used for 
the canonical process on $\bbW$, thus it is 
an ${\sf r}$-dimensional Wiener process. 
\footnote{$K^n(s,t)\in
\bbR^d\otimes\bbR^{{\sf r}}\otimes\bbR^{{\sf r}}
\simeq \bbR^{{\sf r}}\otimes\bbR^{{\sf r}}\otimes\bbR^d$. 
Often $d=k(k+1)/2$ for some basic dimension $k$. }
}}
\end{en-text}

\begin{comment}
Define a $d$-dimensional process $M^n$ by 
\beas 
M^n_t&=&c_n  D^*(1_{[0,t]}A^n\odot D^*B^n), 
\eeas
where 
$c_n$ is a positive constant for scaling, 
$A^n\in\bbD_{{\coloro{\ell+1}},\infty}(H\otimes\bbR^d)$ and 
$B^n\in\bbD_{{\coloro{\ell+1}},\infty}(H\otimes H\otimes\bbR^d)$. 
%
{\coloro{
Here $\odot$ denotes the Hadamard product (entrywise product): 
$(A\odot B)^i_{\alpha,t}=A^i_{\alpha,t}B^i_{\alpha,t}$ 
for $i\in\{1,...,{\sf d}\}$, $\alpha\in\{1,...,{\sf r}\}$ and $t\in[0,1]$ 
for $A,B\in H\otimes\bbR^d$ with identification 
$H\cong L^2([0,1];\bbR^{{\sf r}})$. 
}}
%
%
For the martingale property of $M^n$, it is assumed that 
the process $s\to D^*B^n(s,\cdot)$ is progressively measurable. 
\end{comment}
%
\begin{comment}\rm 
Often $B^n(s,r)$ works for $r\leq s$ but up to $s$, and apparently it opposes to 
the adaptivity of $D^*B^n(s,\cdot)$, however this assumption 
places no restriction on the results. 
In practice, the ``$s$''-part of $B^n(s,r)$ is in most case deterministic 
or $\calf_r$-measurable. Moreover, in applications of our results, 
$M^n_1$ is the principal part of $Z_n$, and usually we can 
take $\calf_{t_{j-1}}$-measurable $D^*B^n(s,\cdot)$ as well as $A^n(s)$ 
for $s\in(t_{j-1},t_j]$. 
\end{comment}

\subsection{Kernel of the quadratic variation}

{\coloroy{Hereafter, let $\ell=\check{d}+8$. }}
In order to fix ideas, we will consider the kernel function
%
$K^n(s,r)$ defined by 
\bea\label{211003-1}
K^n(s,r) 
&=& 
r_n^{-1}\sum_j1_{(t_{j-1},t_j]}(s)\dot{K}^n(s)
\oslash
1_{(t_{j-1},s]}(r)\ddot{K}^n(r)
\eea
for 
{\coloro{
$\dot{K}^n\in\bbD_{\ell+1,\infty}(H\otimes\bbR^d)$ and 
$\ddot{K}^n\in\bbD_{\ell+1,\infty}(H\otimes\bbR^{{\sf r}}\otimes\bbR^d)$,
}} 
{\coloro{where $\oslash$ stands for the tensor product 
with partial Hadamard product given the index $(i,\alpha)\in\{1,...,d\}\times \{1,...,{\sf r}\}$. 
Here 
$H$ is identified with $L^2([0,1];\bbR^{{\sf r}})$, and 
the partial Hadamard product 
(entrywise product) 
given an index {\coloraka $\lambda\in\Lambda$} 
is defined as follows: 
{\colb 
for $a=(a_{\lambda j})_{\lambda\in\Lambda,j\in\calj}\in\bbR^{\Lambda}\otimes\bbR^{\calj}$ and 
$b=(b_{{\coloraka \lambda}k})_{\lambda\in\Lambda,k\in\calk}\in\bbR^{\Lambda}\otimes\bbR^{\calk}$, 
$a\oslash b\in\bbR^\Lambda\otimes\bbR^{\calj}\otimes\bbR^{\calk}$ is given by
\beas 
a\oslash b &=& (a_{\lambda j}b_{\lambda k})_{\lambda\in\Lambda, j\in\calj,  k\in\calk}. 
\eeas
When $\calj$ and $\calk$ are one-point sets, $\oslash$ is a usual Hadamard product $\odot$. 
}
}}
The sequence 
$\{t_j\}_{j=0,1,...,\bar{j}^n}$ $(n\in\bbN)$ {\coloro{is}} a triangular array 
of numbers such that $t_j=t^n_j$ depending on $n$, 
$0=t_0<t_1<\cdots<t_{\bar{j}^n}=1$, 
$\max_j|I_j|=o(r_n)$ and the sequence of measures 
\bea\label{220131-3}
\mu^n=r_n^{-2}\sum_j |I_j|^2 \delta_{t_{j-1}} 
&\to&
\mu
\eea
weakly for some 
{\coloraka measure} 
$\mu$ on $[0,1]$
{\coloraka with  a bounded derivative}, where 
$I_j=(t_{j-1},t_j]$.  
Suppose that $\dot{K}^n$ and $\ddot{K}^n$ are progressively measurable. 
More strongly, 
{\coloro{we assume the {\bf strong predictability 
condition} that $\dot{K}^n(s)$ is $\calf_{t_{j-1}}$-measurable 
for $s\in(t_{j-1},t_j]$.
\footnote{Obviously this condition is satisfied if one replaces 
$\dot{K}^n(s)$ by $\dot{K}^n(t_{j-1})$ 
in the representation of the kernel function $K^n(s,r)$. 
In examples, 
it will turn out that this replacement does not 
cause any practical difficulty. 
}
 }}
We write 
\beas 
\bar{K}^n(s,r)&=&\dot{K}^n(s)\oslash\ddot{K}^n(r). 
\eeas
{\coloro{

\begin{en-text}
then 
\beas 
K^n(s,r) 
&=&
r_n^{-1}\sum_j1_{(t_{j-1},t_j]}(s)
1_{(t_{j-1},s]}(r)\bar{K}^n(s,r).
\eeas
\end{en-text}
%
\begin{en-text} 
By definition, 
the stochastic integral of $K^n(s,\cdot)$ has 
a usual meaning and 
\beas 
\int_0^1 K^n(s,r)dw_r
&=&
r_n^{-1}\sum_j1_{(t_{j-1},t_j]}(s)
\dot{K}^n(s)
\oslash\int_{t_{j-1}}^s\ddot{K}^n(r)dw_r, 
\eeas
and
\end{en-text} 

{\colr{
Corresponding to the representation (\ref{211003-1}), we consider 
$M^n_t$ given by 
\beas 
M^n_t 
&=&
r_n^{-1}\sum_j 
\int_{t_{j-1}\wedge t}^{t_j\wedge t}
\dot{K}^n(s)
\oslash
\Bigl(\int_{t_{j-1}}^s \ddot{K}^n(r)dw_r\Bigr) dw_s,
\eeas
%
{\coloraka where $w$ is the canonical process on $\bbW$, extended naturally to $\bar{\Omega}$. } 

In this case, 
\begin{en-text}
\beas 
C_n &=& 
{\colorr{\mbox{Tr}^o}}
\int_0^1 \Bigl(\int_0^1 K^n(s,r)dw_r\Bigr)^{\otimes2} ds 
\eeas
\end{en-text}
\begin{en-text}
\beas 
C_n &=& 
{\colorr{\mbox{Tr}^o}}
\int_0^1 \bigg(\dot{K}^n(s)
\oslash
\Bigl(\int_{t_{j-1}}^s \ddot{K}^n(r)dw_r\Bigr) dw_s
\bigg)^{\otimes2} ds 
\eeas
\end{en-text}
\beas 
C_n &=& 
\mbox{Tr}^o\int_0^1
r_n^{-2}\sum_j1_{I_j}(s)\Big(\dot{K}^n(s)\oslash
\int_{t_{j-1}}^s\ddot{K}^n(r)dw_r\Big)^{\otimes2}\>ds
\eeas
and it will turn out that the in-p limit is 
\bea\label{210930-1} 
C_\infty &=& 
{\colorr{\half\>}}
\mbox{Tr}^*
\int_0^1 \bar{K}^\infty(t,t)^{\otimes2} \mu(dt), 
\eea 
where 
$\mbox{Tr}^o$ and 
$\mbox{Tr}^*$ denote the traces of 
the element of 
$\mbox{L}(\bbR^{\sf r};\bbR^{\sf r})$ and 
$\mbox{L}(\bbR^{\sf r}\otimes\bbR^{\sf r}
;\bbR^{\sf r}\otimes\bbR^{\sf r})$, 
respectively. 

Let $\Delta=\{(s,r);\> 0\leq r\leq s\leq 1\}$. 
Let $\Delta^n=\cup_j\{(s,r);\>t_{j-1}\leq r \leq s \leq t_j\}$. 
\begin{en-text}
\begin{description}
\item[[C$^\flat$\!\!]] $K^n$ satisfy the following conditions. 
\begin{description}
\item[(i)] 
$\sup_{n\in\bbN}\|K^n(0,0)\|_p<\infty$ for every $p>1$. 
\item[(ii)] 
There exists $c>0$ such that for every $p>1$, 
\beas 
\|\bar{K}^n(s,r)-\bar{K}^n(s',r')\|_p
\leq C_p(|s-s'|^c+|r-r'|^c)
\sskip ((s,r),(s',r')\in\Delta,\>n\in\bbN)
\eeas 
for some constant $C_p$ possibly dependeing on $p$. 
But it seems that $L^p$-continuity in $(s,r)$ is sufficient for our use. 
\item[(iii)] 
{\coloro{$\lim_{n\to\infty}\bar{K}^n(s,r)=\bar{K}^\infty(s,r)$ $P\times{\mathfrak l}\times{\mathfrak l}$-a.e. }} 
\footnote{${\mathfrak l}$ stands for the Lebesgue measure.}
\end{en-text}
%
\begin{comment}
\begin{description}
\item[[C$^\flat$\!\!]] $K^n$ satisfy the following conditions.
\footnote{Of course, in order to show this, the assumption is much stronger than necessary. }
\begin{description}
\item[(i)] 
$\sup_{n\in\bbN}\sup_{(s,r)\in\Delta}\|K^n(s,r)\|_p<\infty$ for every $p>1$. 
\item[(ii)] 
For every $p>1$, 
\beas 
\sup_{(s,r)\in\Delta^n}
\|\bar{K}^n(s,r)-\bar{K}^n(t_{j-1},t_{j-1})\|_p
&\to& 0
\eeas 
as $n\to\infty$. 
%
\item[(iii)] 
There is a continuous process $\bar{K}^\infty(t,t)$ such that 
\beas 
\sup_j\sup_{t\in I_j}\|\bar{K}^n(t_{j-1},t_{j-1})-\bar{K}^\infty(t,t)\|_1\to0
\eeas
as $n\to\infty$. 
\footnote{
The reader may use  ``$\bar{K}^\infty(t)$'' for $\bar{K}^\infty(t,t)$ while we prefer the latter notation.} 

\begin{en-text}
\item[(iii)] 
$\lim_{n\to\infty}\Big\|\sup_{(s,r)\in\Delta}|K^n(s,r)-K^\infty(s,r)|\Big\|_p=0$ for every $p\geq2$. 

\item[(iv)] 
Similar continuity for $Q_s=D_sW_\infty,\> D_sC_\infty,\> D_sF_\infty$. 

\item[(ii)] 
The mappings $Q^n:[0,1]\ni s\mapsto L^p(\Omega;\bbR^d)$ $(n\in\bbN)$ 
are equi-continuous 
for $Q^n=\dot{K}^n$ and $\ddot{K}^n$. 

\item[(iii)] 
$\lim_{n\to\infty}\dot{K}^n(t)=\dot{K}^\infty(t)$ $P\times{\mathfrak l}$-a.e.  
and 
$\lim_{n\to\infty}\ddot{K}^n(t)=\ddot{K}^\infty(t)$ $P\times{\mathfrak l}$-a.e.  
\footnote{${\mathfrak l}$ stands for the Lebesgue measure.} 
\end{en-text}
\end{description}
\end{description}
\end{comment}
%
%
Let 
${\boldsymbol{\delta}}_n=r_n^{-2}\sum_j|I_j|^2=O(1)$ as assumption. 
%
\begin{lemma}
Suppose that 
$K^n$ satisfy the following conditions.
\footnote{Clearly, in order to show the result, the assumption is much stronger than necessary. }
\begin{description}
\item[(i)] 
$\sup_{n\in\bbN}\sup_{(s,r)\in\Delta^{{\colb n}}}\|K^n(s,r)\|_p<\infty$ for every $p>1$. 
\item[(ii)] 
For every $p>1$, 
\beas 
{\colb \sup_j\sup_{(s,r)\in\Delta^n,\>s\in I_j}}
\|\bar{K}^n(s,r)-\bar{K}^n(t_{j-1},t_{j-1})\|_p
&\to& 0
\eeas 
as $n\to\infty$. 
\item[(iii)] 
There is a continuous process $\bar{K}^\infty(t,t)$ such that 
\beas 
\sup_j\sup_{t\in I_j}\|\bar{K}^n(t_{j-1},t_{j-1})-\bar{K}^\infty(t,t)\|_1\to0
\eeas
as $n\to\infty$. 
\footnote{
The reader may use  ``$\bar{K}^\infty(t)$'' for $\bar{K}^\infty(t,t)$ while we prefer the latter notation.} 
\ed
%
Then 
(\ref{210930-1}) holds true.  
\end{lemma}
\proof 
Condition (ii) 
together with the convergence of $\mu^n$ 
ensures
\beas 
C_n 
&=& 
\mbox{Tr}^o\int_0^1
r_n^{-2}\sum_j1_{I_j}(s)\Big(\dot{K}^n(s)\oslash
\int_{t_{j-1}}^s\ddot{K}^n(r)dw_r\Big)^{\otimes2}\>ds
\\&=&
{\mathfrak C}_n+o_p({\boldsymbol{\delta}}_n), 
\eeas
where 
\beas
{\mathfrak C}_n
&=&
\mbox{Tr}^o\int_0^1
r_n^{-2}\sum_j1_{I_j}(s)\Big(\dot{K}^n(t_{j-1})\oslash
\int_{t_{j-1}}^s\ddot{K}^n(t_{j-1})dw_r\Big)^{\otimes2}\>ds
\eeas
and $o_p({\boldsymbol{\delta}}_n)$ denotes a sequence of matrices of indicated order. 
We see 
\bea\label{221017-1} 
{\mathfrak C}_n
&=&
r_n^{-2}\sum_j \mbox{Tr}^o\mbox{Tr}^*\bigg\{
\bar{K}^n(t_{j-1},t_{j-1})^{\otimes2}\otimes\int_{t_{j-1}}^{t_j}
\Bigl\{ \Bigl( \int_{t_{j-1}}^s dw_r \Bigr)^{\otimes2} 
-\int_{t_{j-1}}^s \I_{{\sf r}}dr \Bigr\}\>ds\bigg\}
\nn\\&&
+
r_n^{-2}\sum_j \mbox{Tr}^o\mbox{Tr}^* \bigg\{
\bar{K}^n(t_{j-1},t_{j-1})^{\otimes2}\otimes\int_{t_{j-1}}^{t_j}
\int_{t_{j-1}}^s \I_{{\sf r}}dr \>ds\bigg\}. 
\eea
The square of $L^2$ norm of the first term on the right-hand side 
is of the order 
\beas 
r_n^{-4}\sum_j|I_j|^4
&\leq& 
(r_n^{-1}\max_j|I_j|)^2   r_n^{-2}\sum_j|I_j|^2
\to 0
\eeas
as $n\to\infty$. 
The second term converges in probability to the right-hand side of 
(\ref{210930-1}) by the $L^p$-equi-continuity of the random kernels. 
\qed \ \\

\begin{remark}\rm 
It may seem that we can specify the behavior of $\dotc_n$, however it is still related to 
the It\^o expansion of the process $\dot{K}^n(s)$, so we have left this procedure 
to each individual case. 
Later we will consider the situation where $\dot{K}^n(s)=\dot{K}^n(t_{j-1})$ 
for $s\in(t_{j-1},t_j]$. It is an easily tractable case in the above sense. 
We will do with it there. \\
\end{remark}

\begin{comment}
[For $\dotc\>\!\!^n_t:=\sqrt{n}\Big(C^n_t - C^\infty_t\Big)$, we have \koko 
\beas 
\dotc\>\!\!^n_t
&=& 
\sum_{j:t_j\leq t}\int_{t_{j-1}}^{t_j}1_{[0,t]}(s)\>4n\sqrt{n}a(w_{t_{j-1}})^2
\Big\{\Big(\int_{t_{j-1}}^sdw_r\Big)^2-(s-t_{j-1})\Big\}ds
\\&&
-
2\sqrt{n}\sum_{j:t_j\leq t} \int_{t_{j-1}}^{t_j} \Big(a(w_s)^2-a(w_{t_{j-1}})^2\Big)ds
+O_p\Big(\frac{1}{\sqrt{n}}\Big)
\\&=& 
\sum_{j:t_j\leq t}\int_{t_{j-1}}^{t_j}1_{[0,t]}(s)\>4n\sqrt{n}a(w_{t_{j-1}})^2
\Big\{\Big(\int_{t_{j-1}}^sdw_r\Big)^2-(s-t_{j-1})\Big\}ds
\\&&
-
2\sqrt{n}\sum_{j:t_j\leq t} \int_{t_{j-1}}^{t_j}
2a(w_{t_{j-1}})a'(w_{t_{j-1}})(w_s-t_{j-1})ds
+O_p\Big(\frac{1}{\sqrt{n}}\Big) 
\sskip\mbox{(this line can be omitted)}
\\&=& 
\sum_{j:t_j\leq t}\int_{t_{j-1}}^{t_j}1_{[0,t]}(s)\>4n\sqrt{n}a(w_{t_{j-1}})^2
\Big\{\Big(\int_{t_{j-1}}^sdw_r\Big)^2-(s-t_{j-1})\Big\}ds
\\&&
+O_p\Big(\frac{1}{\sqrt{n}}\Big). 
\eeas
%
Here we know that 
the supremum of ``$O_p(n^{-1/2})$'' in $t\in[0,1]$ is of $O_p(n^{-1/2})$. ]
%
\end{comment}

\begin{comment}
{\coloro{
Comment. Recall
\beas 
e^n_t(u)
&=&
\exp\l(iM^n_t[u]+\half C^n_t[u^{\otimes2}]\r)
\eeas

\beas
\Psi_n(u,v)
&=&
\exp\l\{ iW_n[u]-\half C_n[u^{\otimes2}]+iF_n[v] \r\} 
\eeas
}}
\end{comment}

%
{\coloro{
Let $\bbI_s=-\half(C^\infty_s-C^\infty_1)$, {\coloraka with} {\colb $C^\infty_s$ being the limit of $C^n_s$}. 

\bd
\item[[A1\!\!]] 
{\bf (i)} 
$\bar{K}^n\in\bbD_{\ell+1}(H\otimes H\otimes \bbR^d)$ and 
a representation density of each derivative admits 
\beas 
ess.\sup_{r_1,...,r_k,\in(0,1),\atop(s,r)\in\Delta^n,\>n\in\bbN}
\bigg\|D_{r_1,...,r_k}\bar{K}^n(s,r)\bigg\|_p<\infty
\eeas
for every $p\in(1,\infty)$ and $k\leq\ell+1$. 
\begin{description}
\item[(ii)] For every $\eta>0$ and $p\in(1,\infty)$,
\beas 
\sup_{s\in(0,1),e\in S^{d-1}}
\bigg\|\bigg[\frac{\bbI_s[e^{\otimes2}]}{(1-s)^{1+\eta}}\bigg]^{-1}\bigg\|_p<\infty. 
\eeas 

\item[(iii)] 
$\sup_{s\in[0,1],n\in\bbN}|C^n_s-C^\infty_s|\leq r_n^{2q}$ whenever $|\xi_n|\leq1$. 

\item[(iv)] 
For every $p>1$, 
\beas
\sup_j\sup_{(s,r)\in\Delta^n\atop s\in I_j} \|\bar{K}^n(s,r)-\bar{K}^n(t_{j-1},t_{j-1})\|_{\ell,p}
&=& O(r_n^{2q})
\eeas
and
\beas 
\sup_j\sup_{t\in I_j}
\|\bar{K}^{\colb n}(t_{j-1},t_{j-1})-\bar{K}^\infty(t,t)\|_{\ell,p}&=& O(r_n^{2q})
\eeas
as $n\to\infty$. 
\ed

\ed 

\begin{remark}\rm 
{\bf (i)}  
We note that $2q=1-a\in(\frac{2}{3},1)$ for (iii); in typical cases, $r_n^{2q}=n^{-q}>n^{-\half}$. 
{\coloraka It is essentially possible to remove (iii) by redefining $\xi_n$, as it will be done under another set of conditions later. 
However, we here keep (iii) because it shows a role of $\xi_n$ and making such $\xi_n$ in each case is rather routine. 
 }
\bd
\item[(ii)] 
Under [A1] (i), 
for every $p\in(1,\infty)$ and $k\leq\ell$, 
\beas 
ess.\sup_{r_1,...,r_k,s\in(0,1)}
\bigg\|\frac{|D_{r_1,...,r_k}\bbI_s|}{1-s}\bigg\|_p<\infty.
\eeas

\item[(iii)] 
As for [A1] (ii), 
in order to do with $\exp(\half(C^\infty_s-C^\infty_1)[u^{\otimes2}]))$ 
for $s$ near $1$, we use 
{\coloro{the nondegeneracy of the derivative of $C^\infty_s$ in $s$, or}} 
a large deviation argument. 

\item[(iv)] 
{\coloraka 
The nondegeneracy $\det C_\infty^{-1}\in\cap_{p\geq2} L^p$ follows from [A1] (ii). 
Indeed, it implies 
\beas 
\sup_{e\in\cals^{d-1}}P[C_\infty[e^{\otimes2}]<\ep] &\leq& \ep^p\sup_{e\in\cals^{d-1}}\|\bbI_0[e^{\otimes2}]^{-1}\|_p^p
\>\leq\>C_p\ep^p\sskip(\ep>0)
\eeas
for some constant $C_p$ for every $p>1$. 
Then the desired inequality is obtained; see e.g. Lemma 2.3.1 of Nualart \cite{Nualart2006}. 
}

\ed
\end{remark}

We assume: 

\bd
\item[{[A2]}$^\natural$] 
{\bf (i)}  
$F_\infty\in\bbD_{\ell+1,\infty}(\bbR^{d_1})$ 
and 
$W_\infty\in\bbD_{{\coloro{\ell+1}},\infty}(\bbR^d)$
. 

\bd
\item[(ii)]
$F_n\in\bbD_{\ell+1,\infty}(\bbR^{d_1})$, 
$W_n\in\bbD_{{\coloro{\ell+1}},\infty}(\bbR^d)$, 
$N_n\in\bbD_{{\coloro{\ell+1}},\infty}(\bbR^d)$  
and $\xi_n\in\bbD_{\ell,\infty}(\bbR)$. 
Moreover, 
\beas 
\sup_{{\coloraka n\in\bbN}}\Big\{ 
{\coloraka \|\dotc_n\|_{\ell,p}+}
\|\dotw_n\|_{{\coloro{\ell+1}},p}
+\|\dotf_n\|_{\ell+1,p}+\|N_n\|_{{\coloro{\ell+1}},p}+\|\xi_n\|_{\ell,p} 
\Big\} <\infty
\eeas
for every $p\geq2$. 

\begin{en-text}
\item[(iii)]
{\coloro{$\lim_{n\to\infty}P\l[|\xi_n|\leq\half\r]=1$. }}

\item[(iv)]
$|C_n-C_\infty|>r_n^{1-a}$ implies $|\xi_n|\geq 1$, 
where $a\in(0,1/3)$ is a constant.

\item[(v)] For every $p\geq2$, 
\beas 
\limsup_{n\to\infty}E\Big[1_{\{|\xi_n|\leq1\}} \Delta_{(M_n+W_\infty,F_\infty)}^{-p}\Big] 
<\infty
\eeas
and {\colb moreover} 
$C_\infty^{-1}\in\cap_{p\geq2} L^p$. 
\end{en-text}

\item[(iii)] 
%
{\coloraka $\underline{\sigma}\in \cals(\ell_*,2,1)$.} 
\item[(iv)] 
$({\coloro{M^n_\cdot}},N_n,\dotc_n,\dotw_n,\dotf_n)
\iku^{d_s(\calf)}
(M^\infty_\cdot,N_\infty,\dotc_\infty,\dotw_\infty,\dotf_\infty)$. 

\item[(v)] 
For $G=W_\infty$ and $F_\infty$, 
\beas 
ess. \sup_{r_1,...,r_k\in(0,1)} 
\|D_{r_1,...,r_k}G\|_p 
&<&\infty
\eeas
for every $p\in[2,\infty)$ and $k\leq {\coloraka \ell+1}$. 
%
Moreover, $r\mapsto D_rG$ and $(r,s)\mapsto D_{r,s}G$ ($r\leq s$) are 
continuous a.s. 
\ed
\ed

\begin{comment}
\begin{remark}\rm 
\bd

\item[(上から出る)]
$\bar{\sigma}\in\cals(\ell,{\mathfrak m},{\mathfrak n})$ and 
$\calf^{-1}_{(u,v)}[\tilde{\Phi}^{2,0}](z,x)$ 
has the random symbol $\bar{\sigma}$. 

\item[(出る)] $\call\{M^\infty_t|\calf\} = N_d(0,C^\infty_t)$. 
\ed
\end{remark}
\end{comment}

\begin{remark}\rm 
Under the assumptions, 
$\bar{\sigma}\in\cals(\ell,{\mathfrak m},{\mathfrak n})$ and 
$\calf^{-1}_{(u,v)}[\tilde{\Phi}^{2,0}](z,x)$ 
{\coloraka will have} the random symbol $\bar{\sigma}$, and 
$\call\{M^\infty_t|\calf\} = N_d(0,C^\infty_t)$. 
{\colb [A2]$^\natural$ (iv) is a condition for the joint convergence since we have not specified the random variables other than $M^n$.
We set (iii) by a similar reason. }
\end{remark}

The nondegeneracy of $(M^n_t+W_\infty,F_\infty)$ will be necessary. 
\bd
\item[{[A3]}$^\natural$] 
{\bf (i)} 
There exists a sequence $(t_n)_{n\in\bbN}$ in {\colorr{$[0,1]$ with $\sup_nt_n<1$}}  
such that 
%
the family $\{(M^n_t+W_\infty,F_\infty)\}_{t\geq t_n,n\in\bbN}$ is uniformly nondegenerate 
under the truncation by $\xi_n$, namely, 
\beas 
\sup_{t\geq t_n,n\in\bbN}
E\big[1_{\{|\xi_n|\leq1\}}
\big(\det\sigma_{(M^n_t+W_\infty,F_\infty)}\big)^{-p}\big] < \infty
\eeas
for every $p>1$, 
%


\begin{description}
\item[(ii)]
{\coloro{$\lim_{n\to\infty}P\l[|\xi_n|\leq\half\r]=1$. }} \\
\end{description}
\ed

\begin{remark}\rm 
We removed the condition 
``$|C_n-C_\infty|>r_n^{1-a}$ implies $|\xi_n|\geq 1$, 
where $a\in(0,1/3)$ is a constant'' 
because we can modify the definition of $\xi_n$ {\colb to satisfy this condition}.  
It is possible thanks to $L^p$-boundedness of $\dotc_n$. 
As for $t_n$,   typically $t_n=1/2$. 

\begin{comment}
Comment: We have removed ``(ii)  
There exists a constant $a_0\in[0,1)$ such that 
the condition that $C^\infty_{t_n}>a_0C_\infty$
implies $|\xi_n|>1$.\footnote{ Recall $C_\infty=C^\infty_1$.}  ''. 
It is verified with [A1](ii). 
\end{comment}

\end{remark}


\subsection{Anticipative random symbol {\coloraka and estimates of Fourier transforms}}

We are working with [A1], [A2]$^\natural$ and [A3]$^\natural$. 
We denote $A^n_j(s)=1_{(t_{j-1},t_j]}(s)\dot{K}^n(s)$ and 
$B^n_j(s,r)=1_{\{t_{j-1}<r<s\leq t_j\}}\ddot{K}^n(r)$. 
{\coloro{
Let $A^n_j=(A^{n,\lambda}_{j,\alpha})_{\lambda=1,...,{\sf d}\atop \alpha=1,...,{\sf r}}$ 
and $B^n_j=(B^{n,\lambda}_{j,\alpha,\beta})_{\lambda=1,...,{\sf d}\atop \alpha,\beta=1,...,{\sf r}}$. 
The $\alpha$-th entry of the density function $D_tF$ taking values in $\bbR^{{\sf r}}$ 
will be denoted by $D^{(\alpha)}_tF$. 
We apply the IBP formula to obtain  
\beas 
E\l[L^n_1(u)\Psi_\infty(u,v)\psi_n\r]
&=&
r_n^{-1}\sum_j{\coloroy i}E\Bigl[D^*
\Big(e^n_\cdot(u)A^n_j(\cdot){\colorro{\odot}} D^*B^n_j(\cdot)\Big)[u]
\Psi_\infty(u,v)\psi_n\Bigr]
\\&=&
r_n^{-1}\sum_j
iE\bigg[\bigg\langle \Big(e^n_\cdot(u)A^n_j(\cdot)
{\colorro{\odot}} D^*B^n_j(\cdot)\Big)[u],
D\l(\Psi_\infty(u,v)\psi_n\r)\bigg\rangle_H\bigg]
\\&=&
ir_n^{-1}\sum_j\sum_{\lambda,\alpha}\int_0^1 E\l[
{\colorro{
u^\lambda e^n_s(u)A^{n,\lambda}_{j,\alpha}(s) D^*B^{n,\lambda}_{j,\alpha,\cdot}(s)
D^{(\alpha)}_s\l(
\Psi_\infty(u,v)\psi_n\r)
}}
\r]\>ds
\\&=&
ir_n^{-1}\sum_j\sum_{\lambda,\alpha}\int_0^1
{\colorro{ E\l[
u^\lambda 
\bigg\langle B^{n,\lambda}_{j,\alpha,\cdot}(s,\cdot),
D\Big\{ e^n_s(u)A^{n,\lambda}_{j,\alpha}(s) D^{(\alpha)}_s\l(
\Psi_\infty(u,v)\psi_n\r)\Big\} \bigg\rangle_H\r]\>ds
}}
\\&=&
{\colorro{
ir_n^{-1}\sum_j\sum_{\lambda,\alpha,\beta} \int_0^1 \int_0^1 
E\bigg[
u^\lambda 
B^{n,\lambda}_{j,\alpha,\beta}(s,r)
D^{(\beta)}_r\Big\{ e^n_s(u)A^{n,\lambda}_{j,\alpha}(s) D^{(\alpha)}_s\l(
\Psi_\infty(u,v)\psi_n\r)\Big\} 
\bigg]\>dsdr.
}}
\eeas
}}

%
\begin{en-text}
\begin{remark}\rm 
The gap $\dot{K}(s)-\dot{K}_{t_{j-1}}$ does not  contribute 
even to the second order term of the asymptotic expansion. 
In one-dimensional and equi-spaced case for simplicity, 
\beas 
&&
E\bigg[ \bigg(\sqrt{n}\sum_j \int_{t_{j-1}}^{t_j}\dot{K}(s)\Big(\int_{t_{j-1}}^s \ddot{K}(r)dr\Big)
dw_s \bigg)^2 \bigg]
\\&=&
E\bigg[ n\sum_j\int_{t_{j-1}}^{t_j}\dot{K}(s)^2
\Big(\int_{t_{j-1}}^s\ddot{K}(r)dw_r\Big)^2 ds\bigg]
\\&=&
E\bigg[ n\sum_j\int_{t_{j-1}}^{t_j}\dot{K}(s)^2
\Big(\int_{t_{j-1}}^s\ddot{K}(r)dw_r\Big)^2 ds\bigg]
\eeas 
\end{remark}
\end{en-text}

%
\begin{comment}
Question. Does the gap $\dot{K}(s)-\dot{K}(t_{j-1})$ appear in the second-order limit? 
\end{comment}

\begin{en-text}
{\coloro{For a while, we assume the {\bf strong predictability 
condition} that $\dot{K}^n(s)$ is $\calf_{t_{j-1}}$-measurable 
for $s\in(t_{j-1},t_j]$. 
Obviously this condition is satisfied if one replaces 
$\dot{K}^n(s)$ by $\dot{K}^n(t_{j-1})$ 
in the representation of the kernel funciton $K^n(s,r)$. 
It will turn out that this replacement does not affect 
the expansion up to the second order. 
If $\dot{K}^n$ is strongly predictable, then 
\end{en-text}
{\coloro{Since}} $\dot{K}^n$ is strongly predictable, 
\beas 
E\l[L^n_1(u)\Psi_\infty(u,v)\psi_n\r]
&=&
i
\int_0^1 \int_0^1 
E\bigg[
\bigg\langle 
K^n(s,r)[u],
D_r\Big\{ e^n_s(u) D_s\l(
\Psi_\infty(u,v)\psi_n\r)\Big\} 
\bigg\rangle_{\bbR^{{\sf r}}\otimes\bbR^{{\sf r}}}
\bigg]\>dsdr. 
\eeas
}}
%
%
\begin{comment}
{\coloro{
Comment. 
Old version. 
Since 
\beas 
E\l[L^n_1(u)\Psi_\infty(u,v)\psi_n\r]
&=&
iE\l[\int_0^1 e^n_s(u) dM^n_s[u]\Psi_\infty(u,v)\psi_n\r]
\\&=&
iE\bigg[D^*\l(e^n_\cdot(u) D^*K^n(\cdot)\r)[u]\Psi_\infty(u,v)\psi_n\bigg]
\\&=&
iE\bigg[\bigg\langle e^n_\cdot(u) D^*K^n(\cdot)[u],
D\l(\Psi_\infty(u,v)\psi_n\r)\bigg\rangle_H\bigg]
\\&=&
i\int_0^1 E\l[e^n_s(u)D^*K^n(s)[u]\cdot D_s\l(
\Psi_\infty(u,v)\psi_n\r)\r]\>ds
\\&=&
i\int_0^1 E\l[\bigg\langle K^n(s,\cdot)[u],
D\Big\{ e^n_s(u)D_s\l(
\Psi_\infty(u,v)\psi_n\r)\Big\} \bigg\rangle_H\r]\>ds
\\&=&
i\int_0^1 \int_0^1 E\bigg[\bigg\langle K^n(s,r)[u],
D_r\Big\{ e^n_s(u)D_s\l(
\Psi_\infty(u,v)\psi_n\r)\Big\} \bigg\rangle_{\bbR^{{\sf r}}\otimes\bbR^{{\sf r}}}\bigg]\>dsdr
\\&=&
\eeas
}}
\end{comment}
%
%
or more generally 
{\colorr{ 
\beas 
&&
\partialbs^\alpha
E\l[L^n_1(u)\Psi_\infty(u,v)\psi_n\r]
\\&=&
\sum_{\alpha_0+\alpha_1=\alpha}
c_{\alpha_0,\alpha_1}
i\int_0^1 \int_0^1 E\bigg[\bigg\langle 
\partialbs^{\alpha_0}K^n(s,r)[u],
\partialbs^{\alpha_1}
D_r\Big\{ e^n_s(u)D_s\l(
\Psi_\infty(u,v)\psi_n\r)\Big\} \bigg\rangle_{\bbR^{{\sf r}}\otimes\bbR^{{\sf r}}}\bigg]\>dsdr
\eeas
with some constants $c_{\alpha_0,\alpha_1}$ appearing in Leibniz's rule 
for $\partialbs^\alpha$. 
We will find a representation of the random symbol $\bar{\sigma}$ 
associated to $\calf^{-1}_{(u,v)}[\tilde{\Phi}^{2,0}](z,x)$ 
and verify the convergence (\ref{220131-1}). 
}}



The density function appearing here is 
\beas &&
D_r\Big\{ e^n_s(u)D_s\l(\Psi_\infty(u,v)\psi_n\r)\Big\} 
\\&=&
D_r\Big\{ e^n_s(u)D_s\l(
\exp\l\{ iW_\infty[u]-\half C_\infty[u^{\otimes2}]+iF_\infty[v] \r\} \psi_n\r)\Big\}
\\&=&
e^n_s(u)
\Psi_\infty(u,v)\psi_n
\l(iD_rM^n_s[u]+\half D_rC^n_s[u^{\otimes2}]\r)\otimes 
\Big(iD_sW_\infty[u]-\half D_sC_\infty[u^{\otimes2}]+iD_sF_\infty[v] \Big)
\\&&
+e^n_s(u)\Psi_\infty(u,v)\psi_n
\Big(iD_rW_\infty[u]-\half D_rC_\infty[u^{\otimes2}]+iD_rF_\infty[v] \Big)
\otimes 
\Big(iD_sW_\infty[u]-\half D_sC_\infty[u^{\otimes2}]+iD_sF_\infty[v] \Big)
\\&&
+e^n_s(u)\Psi_\infty(u,v)\psi_n
\Big(iD_rD_sW_\infty[u]-\half D_rD_sC_\infty[u^{\otimes2}]+iD_rD_sF_\infty[v] \Big)
\\&&
+e^n_s(u)\Psi_\infty(u,v)D_r\psi_n\otimes 
\Big(iD_sW_\infty[u]-\half D_sC_\infty[u^{\otimes2}]+iD_sF_\infty[v] \Big)
\\&&
+
D_r\Big\{ e^n_s(u)\Psi_\infty(u,v)D_s\psi_n\Big\}. 
\eeas
{\colorr{Moreover, 
the Leibniz type formulas for 
$\partialbs^{\alpha_2}
D_r\Big\{ e^n_s(u)D_s\l(\Psi_\infty(u,v)\psi_n\r)\Big\}$ 
remain in force if one applies the Leibniz rule 
to each term on the right-hand side of the above equality. 
It will be observed that the function $|\partialbs^\alpha E\l[L^n_1(u)\Psi_\infty(u,v)\psi_n\r]|$ is 
dominated by a polynomial in $(u,v)$ at most fifth-order 
under a restricted range of $(u,v)$ by means of the truncation with $\psi_n$. 
}}

We introduce an $\bbR^{\sf r}\otimes\bbR^{\sf r}$-valued random symbol
\beas 
\sigma_{s,r}(iu,iv)
&=&
\half D_rC^{{\colorr{\infty}}}_s[u^{\otimes2}]\otimes 
\Big(iD_sW_\infty[u]-\half D_sC_\infty[u^{\otimes2}]+iD_sF_\infty[v] \Big)
\\&&
+
\Big(iD_rW_\infty[u]-\half D_rC_\infty[u^{\otimes2}]+iD_rF_\infty[v] \Big)
\otimes 
\Big(iD_sW_\infty[u]-\half D_sC_\infty[u^{\otimes2}]+iD_sF_\infty[v] \Big)
\\&&
+
\Big(iD_rD_sW_\infty[u]-\half D_rD_sC_\infty[u^{\otimes2}]+iD_rD_sF_\infty[v] 
\Big). 
\eeas

\subsubsection{The terms involving $D_rM^n_s$}\label{221024-4}

We have 
\beas 
D_rM^n_s
&=&
r_n^{-1}\sum_j\int_{t_{j-1}\wedge s}^{t_j\wedge s} 
D_r\dot{K}^n(s_1)
\oslash
\l(\int_{t_{j-1}}^{s_1}\ddot{K}^n(s_2)dw_{s_2}\r)
dw_{s_1}
\\&&
+r_n^{-1}\sum_j\int_{t_{j-1}\wedge s}^{t_j\wedge s}
\dot{K}^n(s_1)
\oslash
\l(\int_{t_{j-1}}^{s_1}D_r\ddot{K}^n(s_2)dw_{s_2}\r)dw_{s_1}
\\&&
+r_n^{-1}\sum_j\int_{t_{j-1}\wedge s}^{t_j\wedge s}
\dot{K}^n(s_1)
\oslash1_{(t_{j-1},s_1]}(r)\ddot{K}^n(r)
dw_{s_1}
\\&&
+r_n^{-1}\sum_j1_{(t_{j-1}\wedge s,t_j\wedge s]}(r)\dot{K}^n(r)
\oslash
\int_{t_{j-1}}^{r}\ddot{K}^n(s_1)dw_{s_1}.
\eeas
With this representation, we see 
\beas 
\int_0^1 K^n(s,r)\otimes D_rM^n_sdr
&=&
\int_0^1 
r_n^{-1}\sum_j1_{(t_{j-1},t_j]}(s)
\dot{K}^n(s)\oslash
1_{(t_{j-1},s]}(r)\ddot{K}^n(r)
\\&&
\otimes 
\bigg\{
r_n^{-1}\sum_{j'}\int_{t_{j'-1}\wedge s}^{t_{j'}\wedge s} 
D_r\dot{K}^n(s_1)
\oslash
\l(\int_{t_{j'-1}}^{s_1}\ddot{K}^n(s_2)dw_{s_2}\r)
dw_{s_1}
\\&&
+r_n^{-1}\sum_{j'}\int_{t_{j'-1}\wedge s}^{t_{j'}\wedge s}
\dot{K}^n(s_1)
\oslash
\l(\int_{t_{j'-1}}^{s_1}D_r\ddot{K}^n(s_2)dw_{s_2}\r)dw_{s_1}
\\&&
+r_n^{-1}\sum_{j'}\int_{t_{j'-1}\wedge s}^{t_{j'}\wedge s}
\dot{K}^n(s_1)1_{(t_{j'-1},s_1]}(r)\oslash\ddot{K}^n(r)
dw_{s_1}
\\&&
+r_n^{-1}\sum_{j'}1_{(t_{j'-1}\wedge s,t_{j'}\wedge s]}(r)
\dot{K}^n(r)\oslash
\int_{t_{j'-1}}^{r}\ddot{K}^n(s_1)dw_{s_1}
\bigg\}
dr.
\eeas
Thus we have 
\beas &&
\int_0^1 K^n(s,r)\otimes D_rM^n_sdr
\\&=&
r_n^{-1}\sum_{j}1_{(t_{j-1},t_{j}]}(s)
\dot{K}^n(s)\oslash
r_n^{-1}
\int_{t_{j-1}}^{s} 
\bigg\{
\int_{t_{j-1}\wedge s}^{t_{j}\wedge s}
\ddot{K}^n(r)
\otimes 
D_r\dot{K}^n(s_1)\oslash
\l(\int_{t_{j-1}}^{s_1}\ddot{K}^n(s_2)dw_{s_2}\r)
dw_{s_1}
\\&&
+
\ddot{K}^n(r)\otimes 
\int_{t_{j-1}\wedge s}^{t_{j}\wedge s}
\dot{K}^n(s_1)\oslash
\l(\int_{t_{j-1}}^{s_1}
D_r\ddot{K}^n(s_2)dw_{s_2}\r)dw_{s_1}
\\&&
+\ddot{K}^n(r)\otimes
\int_{t_{j-1}\wedge s}^{t_{j}\wedge s}
\dot{K}^n(s_1)\oslash1_{(t_{j-1},s_1]}(r)
\ddot{K}^n(r)
dw_{s_1}
\\&&
+1_{(
t_{j-1}\wedge s,t_{j}\wedge s]}(r)
\ddot{K}^n(r)\otimes 
\dot{K}^n(r)\oslash
\int_{t_{j-1}}^{{\sf r}}\ddot{K}^n(s_1)dw_{s_1}
\bigg\}dr. 
\eeas
Here we note that for terms not to vanish, it is necessary that 
\bea\label{200503-1}
r\leq s_1 \leq t_{j'}
\sskip \mbox{for } D_r\dot{K}^n(s_1)\not=0 \mbox{}
\eea
\bea\label{200503-5}
r\leq s_2 \leq s_1 \leq t_{j'}
\sskip \mbox{for } D_r\ddot{K}^n(s_2)\not=0 \mbox{}
\eea
\bea\label{200503-2}
t_{j-1}<r\leq s
\eea
\bea\label{200503-3}
t_{j-1}<s\leq t_j
\eea
\bea\label{200503-4}
t_{j'-1}< s
\eea
In particular, 
$t_{j-1}<t_{j'}$ 
from [(\ref{200503-1}) or (\ref{200503-5})] and (\ref{200503-2}), 
and 
$t_{j'-1}<t_j$ 
from (\ref{200503-3}) and (\ref{200503-4}), 
so that 
$j'=j$. 
Similar argument is valid for the last two terms to neglect off-diagonal 
elements for $j\not=j'$.

{\coloroy{
We will assume that 
\bea\label{220131-2}
r_n^{-8}\sum_j|I_j|^5&=&O(1)
\eea
as $n\to\infty$. 
Then, applying Jensen's inequality, we have 
\bea\label{220211-1}
r^{(1)}_n:=r_n^{-2}\sum_j|I_j|^{\frac{5}{2}}=O(r_n)
\eea
and 
\bea\label{220211-2}
r^{(2)}_n:=r_n^{-2}\sum_j|I_j|^3=O(r_n^2),
\eea
moreover, 
$
r_n^{(2)}\leq O(r_n^\half r_n^{(1)})
$ 
since $\max_j|I_j|=O(r_n)$ 
by assumption. 
}}

Let 
\beas 
Q_s &=& iD_sW_\infty[u]-\half D_sC_\infty[u^{\otimes2}]+iD_sF_\infty[v] . 
\eeas
Without scaling by $r_n^{-1}$, 
\bea
\cali_n
&:=&
i\int_0^1 \int_0^1 E\bigg[\bigg\langle K^n(s,r)[u],
e^n_s(u)\Psi_\infty(u,v)\psi_n
iD_rM^n_s[u]\otimes Q_s
\bigg\rangle_{\bbR^{{\sf r}}\otimes\bbR^{{\sf r}}}\bigg]\>dsdr
\label{220504-1}
\\&=&
-{\mbox{Tr}^*}
E\Bigg[
\int_0^1 ds
e^n_s(u)\Psi_\infty(u,v)\psi_n
r_n^{-1}\sum_{j}1_{(t_{j-1},t_{j}]}(s)
\dot{K}^n(s)
\nn\\&&
\oslash 
r_n^{-1}
\int_{t_{j-1}}^{s} \ddot{K}^n(r)
\otimes 
\bigg\{
\int_{t_{j-1}\wedge s}^{t_{j}\wedge s}
\dot{K}^n(s_1)1_{(t_{j-1},s_1]}(r)\oslash 
\ddot{K}^n(r)
dw_{s_1}
\nn\\&&
+1_{(t_{j-1}\wedge s,t_{j}\wedge s]}(r)
\dot{K}^n(r)\oslash
\int_{t_{j-1}}^{r}\ddot{K}^n(s_1)dw_{s_1}
\bigg\}dr\otimes Q_s
\Bigg]
[u^{\otimes2}]
+O({\coloroy{r_n^{(2)}}})
\nn\\&=&
-{\mbox{Tr}^*}
\int_0^1
E\Bigg[
e^n_s(u)\Psi_\infty(u,v)\psi_n
r_n^{-1}\sum_{j}1_{(t_{j-1},t_{j}]}(s)
\nn\\&&
\cdot r_n^{-1}
\int_{t_{j-1}}^{s} \bar{K}^n(s,r)
\otimes
\bigg\{
\int_{t_{j-1}\wedge s}^{t_{j}\wedge s}
1_{(t_{j-1},s_1]}(r)
\bar{K}^n(s_1,r)
dw_{s_1}
\nn\\&&
+1_{(t_{j-1}\wedge s,t_{j}\wedge s]}(r)\dot{K}^{n}(r)
\oslash\int_{t_{j-1}}^{r}\ddot{K}^{n}(s_1)dw_{s_1}
\bigg\}\otimes Q_s dr
\Bigg][u^{\otimes2}]
\> ds
+O({\coloroy{r_n^{(2)}}})
\nn\\&=&
-{\mbox{Tr}^*}E\Bigg[
\Psi_\infty(u,v)
\psi_n
\sum_j
\xi^n_j(u)
\otimes Q_{t_{j-1}}
\Bigg][u^{\otimes2}] 
+o({\coloroy{r_n^{(1)}}}), 
\nn
\eea
where 
\beas 
\xi^n_j(u)
&=&
e^n_{t_{j-1}}(u)
\bar{K}^n(t_{j-1},t_{j-1})
\otimes 
\bar{K}^n(t_{j-1},t_{j-1})
\> {\coloroy{r_n^{-2}}}
\bigg\{
\int_{t_{j-1}}^{t_{j}}
(t_j-s_1)(s_1-t_{j-1})dw_{s_1}
\\&&
+\int_{t_{j-1}}^{t_j}\frac{(t_j-s_1)^2}{2}dw_{s_1}
\bigg\}.
\eeas 
The last equality 
in the expression of $\cali_n$ 
is by straightforward $L^p$-estimate 
with triangular inequality and $L^p$-continuity of the process $Q_s$; 
note that $Q_s$ may be anticipative and we do not assume the 
predictability of $Q_{t_{j-1}}$, which has no more meaning than 
approximation to $Q_s$. 

Let 
\bea
I^n(u,v)&=&
r_n^{-1}
i\int_0^1 \int_0^1 E\bigg[\bigg\langle 
\partialbs^{\alpha_0}K^n(s,r)[u],
\partialbs^{\alpha_1}\bigg\{
e^n_s(u)
\Psi_\infty(u,v)\psi_n
iD_rM^n_s[u]
\nn\\&&
\otimes 
\Big(iD_sW_\infty[u]-\half D_sC_\infty[u^{\otimes2}]+iD_sF_\infty[v] \Big)
\bigg\}
%
%
\bigg\rangle_{\bbR^{{\sf r}}\otimes\bbR^{{\sf r}}}\bigg]\>dsdr
\label{220322-1}
\eea
%
We will apply the integration-by-parts formula 
{\coloroy{$\check{d}+6$}} 
times 
for ${\colb s}\geq t_n$, 
by taking advantage of the nondegeneracy of 
$\{(M^n_s+W_\infty,F_\infty)\}_{s\geq t_n,n\in\bbN}$ when $|\xi_n|\leq1$
, and that 
of the decay of  {\coloraka $\exp\{\half(C^n_s-C_\infty)[u^{\otimes2}]\}$} 
for $s<t_n$ 
{\coloroy{as well as the integration-by-parts formula for $F_\infty$}} 
in order to obtain integrability in $(u,v)$ of 
$I^n(u,v)$ of \eqref{220322-1}. 
{\coloroy{We need several steps to achieve this plan.}}
\ \\

%
\bd
\item[(a)]  
$(M^n_s+W_\infty,F_\infty)$ and $C^n_s-C_\infty$ 
appear in the factorization
\beas 
e^n_s(u)\Psi_\infty(u,v)
&=&
\exp\bigg(iM^n_s[u]+\half(C^n_s-C_\infty)[u^{\otimes2}]\bigg)
\exp\bigg(iW_\infty[u]+iF_\infty[v]\bigg)
\\&=&
\bbF^n_s\bbG_s\bbH^n_s,
\eeas
where
\beas
\bbF^n_s &=&\exp\bigg(i(M^n_s+W_\infty)[u]+iF_\infty[v]\bigg),
\\
\bbG_s &=& 
\exp\bigg(\half(C^\infty_s-C^\infty_1)[u^{\otimes2}]\bigg),
\\
\bbH^n_s &=&
\exp\bigg(\half(C^n_s-C^\infty_s)[u^{\otimes2}]\bigg). 
\eeas

\item[(b)] 
By the Leibniz rule, the density function 
$D_{r_1,...,r_k}\exp(-\bbI_s[u^{\otimes2}])$ of the $k$th Malliavin derivative 
of $\exp(-\bbI_s[u^{\otimes2}])$ is a linear combination of the terms of the form 
\beas 
\bbA
&=&
D_{r_1,...,r_{k_1}}\bbI_s[u^{\otimes 2}]\otimes
\cdots \otimes D_{r_{k_{m-1}+1},...,r_{k_m}}\bbI_s
[u^{\otimes 2}]
\>
\exp(-\bbI_s[u^{\otimes2}]), 
\eeas
where 
$1\leq k_1<\cdots<k_{m-1}<k_m=k$ and $ 1\leq m\leq k$. 
We denote by $\ep$ a positive number and we will make its value as small as we want 
in the context. It is possible to choose such an $\ep$ because we only change 
its values finitely often. 
Let $u\in\bbR^d\setminus\{0\}$ and let $e_u=|u|^{-1}u$. 
\beas 
|\bbA| 
&\leq& 
\frac{|D_{r_1,...,r_{k_1}}\bbI_s|}{1-s}\cdots
\frac{|D_{r_{k_{m-1}+1},...,r_{k_m}}\bbI_s|}{1-s}
\bigg[\frac{\bbI_s[e_u^{\otimes2}]}{(1-s)^{1+{\colb \eta(\ep)}}}\bigg]^{-m(1-\ep)}
\\&&
\cdot\>\bbI_s[u^{\otimes2}]^{m(1-\ep)}\exp(-\bbI_s[u^{\otimes2}])
|u|^{2m\ep}
\eeas
for 
{\colb $\eta(\ep)=\ep(1-\ep)^{-1}$. }

\item[(c)] 
For every $m$ and $\ep\in(0,1)$, 
$\bbI_s[u^{\otimes2}]^{m(1-\ep)}\exp(-\bbI_s[u^{\otimes2}])$ is bounded 
since $\bbI_s$ is nonnegative-definite and 
$\sup_{x\geq0}|x|^m e^{-x}<\infty$. 
This and the inequality in (b) together with [A1] (i)-(ii) imply 
\beas 
ess. \sup_{r_1,...,r_k,s\in(0,1)} \|\bbA\|_p 
&\leq&
C(p) 
|u|^{2m\ep}
\eeas
for some constant $C(p)$.

\item[(d)] 
$(C^n_s-C^\infty_s)[u,u]$ is bounded 
in $u\in\Lambda^0_n$ 
whenever $\xi_n\leq1$, due to [A1] (iii). 
Therefore, 
{\coloroy{by [A1] (iv), 
$\big\{\|D^k\bbH^n_s1_{\{\xi_n\leq1\}}\|_p;\> s\in[0,1], \>u\in\Lambda(d,q),\>n\in\bbN\big\}$ 
is bounded for every $p>1$}} 
%
for $k\leq\ell$. 
We notice that 
due to (\ref{220131-2}), 
the first term on the right-hand side of (\ref{221017-1}) is 
$O(r_n)$ 
in norms. 

\item[(e)] 
The Sobolev norms of 
$\{D_rM^n_s\}_{(r,s)\in[0,1]^2,n\in\bbN}$ are bounded. 

\item[(f)] 
$\sup_{n\in\bbN} r_n^{-2}\int_0^1\int_0^1
\sum_j1_{(t_{j-1},t_j]}(s)1_{(t_{j-1},s]}(r)\>dsdr<\infty$. 

\item[(g)] 
\begin{en-text}
\item[(h)] In the estimation of $\bbA$ in (b), we can re-estimate the factor 
$\bbI_s[u^{\otimes2}]^{m(1-\ep)}\exp(-\bbI_s[u^{\otimes2}])$ as 
\beas 
\bbI_s[u^{\otimes2}]^{m(1-\ep)}\exp(-\bbI_s[u^{\otimes2}])
&\leq&
\bbI_s[u^{\otimes2}]^{m(1-\ep)}\exp(-\half\bbI_s[u^{\otimes2}])
\cdot \exp(-\half a_0C_\infty[u^{\otimes2}])
\\&\leq& 
{\sf C}\ \exp(-\half a_0C_\infty[u^{\otimes2}])
\\&\leq& 
{\sf C'}\ |u|^{-2m}(C_\infty[e_u^{\otimes2}])^{-m}
\eeas
whenever $\xi_n\leq1$. 
\end{en-text}
{\coloroy{
In the estimation of $\bbA$ in (b), we can re-estimate the factor 
$\bbI_s[u^{\otimes2}]^{m(1-\ep)}\exp(-\bbI_s[u^{\otimes2}])$ as 
\beas 
\bbI_s[u^{\otimes2}]^{m(1-\ep)}\exp(-\bbI_s[u^{\otimes2}])
&\leq& 
{\sf C(k)}\ |u|^{-2k}(\bbI_{t_n}[e_u^{\otimes2}])^{-k}
\eeas
whenever $\xi_n\leq1$, for any $k\in\bbN$. 
It will be applied for $s<t_n$. 
}}

\item[(h)]
With (a), (c), (d), (e) and (f), we repeatedly apply the integration-by-parts formula 
based on $(M^n_s+W_\infty,F_\infty)$ 
at (\ref{220504-1}) for $s\geq t_n$. 
{\coloroy{For $s<t_n$,}}  we apply 
the integration-by-parts formula based on $F_\infty$ (for all $v\in\bbR^{d_1}$) 
with the help of (d) and (g), 
to obtain 
\bea\label{220505-1}
\sup_n\sup_{(u,v)\in\Lambda^0_n(\check{d},q)} |(u,v)|^{\check{d}+2-{\colb \ep'}}
|I^n(u,v)|<\infty
\eea
{\colb for some $\ep'\in(0,1)$. }
We note 
{\colb 
that the differentiation with respect to $(u,v)$ in (\ref{220322-1}) does not change the estimate essentially.
}
Note also 
that $I^n(u,v)$ is like a $4th$-order polynomial in $(u,v)$ in growth rate, and that  
$\ell-4-{\colb \ep'}=\check{d}+2-{\colb \ep'}$.\footnote{Here we have an index larger than necessary. 
But it becomes $\check{d}+1-{\colb \ep'}$ for some other terms. }

\ed
}}

\begin{comment}
, 
also that the Malliavin covariance of $M^n_s$ 
is not less than the predictable blacket $C^n_s$ essentially; see 
memo.
\end{comment}

%
\begin{en-text}
The last equality is by the Schwarz inequality and the $L^2$-estimate 
for the martingale transform, the inside of $E[...]$ without the factor 
$\Psi_\infty(u,v)$. 
\end{en-text}
%
%
Let 
\beas 
\mathscr{M}^n_t 
&=&
r_n^{-1}
\sum_j
e^n_{t_{j-1}}(u)
\bar{K}^n(t_{j-1},t_{j-1})
\\&&
\otimes 
\bar{K}^n(t_{j-1},t_{j-1})
\>  {\coloroy{r_n^{-2}}}
\bigg\{
\int_{t_{j-1}\wedge t}^{t_{j}\wedge t}
(t_j-s_1)(s_1-t_{j-1})dw_{s_1}
+
\int_{t_{j-1}\wedge t}^{t_{j}\wedge t}
\frac{(t_j-s_1)^2}{2}dw_{s_1}
\bigg\}.
\eeas

\begin{en-text}

Then a stable convergence occurs for 
the complex-valued martingale $\mathscr{M}^n$, that is, 
\bea\label{200529-1} 
\mathscr{M}^n &\iku^d \mathscr{M}^\infty
\eea 
in $C([0,1])$ as $n\iku\infty$, and 
$\mathscr{M}^\infty$ can be expressed as 
\beas 
\mathscr{M}^\infty &=& 
\int_0^\cdot 
e^\infty_t(u)\bar{K}^\infty(t,t)\otimes\bar{K}^\infty(t,t)
d\mathscr{B}_t,
\eeas
where 
\beas 
e^\infty_t(u) &=& 
\exp\bigg(iM^\infty_t[u]+\half C^\infty_t[u^{\otimes2}]\bigg)
\eeas
and 
$(M^\infty,\mathscr{B}^\infty)$ is a continuous martingale 
indepdentend of $\calf$; recall that $\calf=\sigma[\calf^w]$. 
(\UTF{00E8}\UTF{00E2}\UTF{00E4}\UTF{02D9}\UTF{00ED}l\UTF{00C7}\UTF{00E6}\UTF{00C7}\UTF{00D8}\UTF{00EB}\UTF{00C2}\UTF{00C7}°≠\UTF{00C7}°\UTF{201A}\UTF{00C7}\UTF{00B5}\UTF{00C7}¶\UTF{220F}\UTF{00C7}\UTF{0178}\UTF{00C7}°\UTF{00AF}\UTF{00C7}\UTF{2122}\UTF{00C7}\UTF{00CA}\UTF{00C7}°\UTF{00D2}\UTF{00C7}\UTF{00A9}\UTF{FFE0}\UTF{00DA}D) 
To obtain the above convergence, we first observe that 
$(M^n,\mathscr{B}^n)$, $\mathscr{B}^n$ corresponding 
to the sum of $dw$ integrals, converges stably to 
$(M^\infty,\mathscr{B}^n)$ that is a continuous martingale 
orthogonal to $w$. Next, the convergence (\ref{200529-1}) 
follows thanks to the theory of weak convergence of 
stochastic integrals. 

\end{en-text}
%

Obviously, the Burkholder-Davis-Gundy inequality 
{\coloroy{applied under (\ref{220131-2})}} implies 
\beas 
\|\sup_{t\in[0,1]}|\mathscr{M}^n_t|\|_p=O(r_n)
\eeas 
for every $p>1$.

For $p>1$ and $\ep>0$, 
{\coloroy{thanks to (\ref{220211-1}),}} 
we can find an increasing sequence 
$T_k$ ($k=0,1,...,{\sf K}$) in $[0,1]$ such that 
$0=T_0<\cdots<T_{{\sf K}}=1$ and 
\beas 
&&
\Bigg| r_n^{-1}
E\Bigg[
\Psi_\infty(u,v)\psi_n\sum_j
\xi^n_j(u)\otimes \big(Q_{T_{k(t_{j-1})}}-Q_{t_{j-1}}\big)
\Bigg]
\Bigg| 
< \ep
\eeas
uniformly for large $n$, where 
$k(t_{j-1})=\max\{k;\> T_k\leq t_{j-1}\}$ depending on $j$ and $n$. 
Moreover, 
\beas
&&
r_n^{-1}
E\Bigg[
\Psi_\infty(u,v)\psi_n\sum_j
\xi^n_j(u)
\otimes Q_{T_{k(t_{j-1})}}
\bigg]
\\&=&
\sum_{k=1}^{{\sf K}} 
E\bigg[ \Psi_\infty(u,v) \psi_n
(\mathscr{M}^n_{T_k}-\mathscr{M}^n_{T_{k-1}})\otimes Q_{T_{k-1}}
\bigg]
+o(1)
\\&=&
o(1)
\eeas
since the stable limit is centered. 
Here the first equality was due to C-tightness 
of the sequence $\mathscr{M}^n$.

\begin{en-text}
Because of the uniform integrability of the variables inside of 
the expectation, we obtain 
\beas 
&&
\lim_{n\iku\infty}
\sum_{k=1}^{{\sf K}} 
E\bigg[ \Psi_\infty(u,v) 
(\mathscr{M}^n_{T_k}-\mathscr{M}^n_{T_{k-1}})\otimes Q_{T_{k(t_{j-1})}}
\bigg]
\\&=&
\sum_{k=1}^{{\sf K}} 
E\bigg[ \Psi_\infty(u,v) 
(\mathscr{M}^\infty_{T_k}-\mathscr{M}^\infty_{T_{k-1}})
\otimes Q_{T_{k(t_{j-1})}}
\bigg]
\\&=&0
\eeas 
by the fact that $\mathscr{M}^\infty$ is a 
conditionally zero-mean martingale given $\calf$. 
Since $\ep$ is arbitrary, 
\beas 
&&
\lim_{n\iku\infty}
\Bigg| r_n^{-1}
E\Bigg[
\Psi_\infty(u,v)\sum_j
e^n_{t_{j-1}}(u)
\bar{K}^n(t_{j-1},t_{j-1})
\\&&
\otimes 
\bar{K}^n(t_{j-1},t_{j-1})
\> n
\bigg\{
\int_{t_{j-1}}^{t_{j}}
(t_j-s_1)(s_1-t_{j-1})dw_{s_1}
\\&&
+\int_{t_{j-1}}^{t_j}\frac{(t_j-s_1)^2}{2}dw_{s_1}
\bigg\}\otimes Q_{t_{j-1}}
\Bigg]
\Bigg| =0. 
\eeas

\end{en-text}

Since $\ep$ is arbitrary, 
as a consequence of the above estimates, we obtain 
\beas
i\int_0^1 \int_0^1 E\bigg[\bigg\langle K^n(s,r)[u],
e^n_s(u)\Psi_\infty(u,v)\psi_n
iD_rM^n_s[u]\otimes Q_s
\bigg\rangle_{\bbR^{{\sf r}}\otimes\bbR^{{\sf r}}}\bigg]\>dsdr
&=&
o(r_n)
\eeas
as $n\iku\infty$. 
{\coloroy{
Even in the case involving derivatives in $(u,v)$, 
following the same argument, we obtain 
\bea\label{220505-2}
I^n(u,v) &=& o(1)
\eea
as $n\to\infty$ for every $(u,v)\in\bbR^{\check{d}}$. 

 }}
\subsubsection{The terms involving $D_rC^n_s$ {\coloraka and others}}\label{221024-5}

By definition, 
\beas 
\int_0^1 K^n(s,r)\otimes D_rC^n_sdr
&=&
r_n^{-1}\sum_j1_{(t_{j-1},t_j]}(s)
\dot{K}^n(s)\oslash
\int_0^1 
1_{(t_{j-1},s]}(r)\ddot{K}^n(r)\otimes D_rC^n_sdr.
\eeas

With scaling by $r_n^{-1}$, 
\bea\label{220505-3}
&&
r_n^{-1}i\int_0^1 \int_0^1 E\bigg[\bigg\langle K^n(s,r)[u],
e^n_s(u)\Psi_\infty(u,v)\psi_n
\half D_rC^n_s[u^{\otimes2}]
\otimes Q_s
\bigg\rangle_{\bbR^{{\sf r}}\otimes\bbR^{{\sf r}}}\bigg]\>dsdr
\\&=&
r_n^{-1}i \>\mbox{Tr}^* \int_0^1 
E\bigg[e^n_s(u)\Psi_\infty(u,v)\psi_n
r_n^{-1}\sum_j1_{(t_{j-1},t_j]}(s)\dot{K}^n(s)\oslash
\int_0^1 1_{(t_{j-1},s]}(r)\ddot{K}^n(r)\otimes 
{\coloro{\half}}D_rC^n_sdr\otimes Q_s
\bigg][u^{\otimes3}]\>ds
\nn\\&=&
{\coloroy{r_n^{-2}i}} \>\mbox{Tr}^* \int_0^1 \sum_j1_{(t_{j-1},t_j]}(s)
{\coloro{(s-t_{j-1})}}
E\bigg[e^n_s(u)\Psi_\infty(u,v)\psi_n
\dot{K}^n(s)\oslash
\ddot{K}^n(t_{j-1})\otimes 
{\coloro{\half}}D_{t_{j-1}}C^n_s\otimes Q_s
\bigg][u^{\otimes3}]\>ds+o(1)
\nn\\&=&
{\coloroy{r_n^{-2}i}} \>\mbox{Tr}^*  \sum_j{\coloro{\half|I_j|^2}}
E\bigg[e^n_{t_{j-1}}(u)\Psi_\infty(u,v)\psi_n
\dot{K}^n({t_{j-1}})\oslash
\ddot{K}^n(t_{j-1})\otimes 
{\coloro{\half}}
D_{t_{j-1}}C^n_{t_{j-1}}\otimes Q_{t_{j-1}}
\bigg][u^{\otimes3}]+o(1)
\nn\\&=&
\frac{i}{2} \>\mbox{Tr}^*  \int_0^1 
E\bigg[e^\infty_t(u)\Psi_\infty(u,v)
\bar{K}^\infty(t,t)\otimes 
{\coloro{\half}}D_tC^\infty_t\otimes Q_t
\bigg][u^{\otimes3}]\>{\coloro{\mu(dt)}} +o(1),
\nn
\eea
where 
$D_tC^\infty_t=\lim_{s\up t}D_sC^\infty_t$. 
For the last equality, the first integration in the last line can be approximated by 
that with respect to $\mu^n$. 
If we realize the weak convergence of random variables including 
$e^n_\cdot(u)$ taking values in the space of continuous functions with uniform norm 
as a.s. convergence, then the integrand of the last line 
can be approximated by $E[e^n_t(u)\Psi_\infty(u,v)\psi_n
\bar{K}^n(t,t)\otimes {\coloro{\half}}D_tC^n_t\otimes Q_t]$ uniformly in $t$ 
with the help of the truncation. 

By conditioning, 
\beas 
&&
E\bigg[e^\infty_t(u)\Psi_\infty(u,v)
\bar{K}^\infty(t,t)\otimes 
D_tC^\infty_t\otimes Q_t
\bigg]
\\&=&
E\bigg[\Psi_\infty(u,v)
\bar{K}^\infty(t,t)\otimes 
D_tC^\infty_t\otimes Q_t
\bigg]
\\&=&
E\bigg[
\exp\l\{ iW_\infty[u]-\half C_\infty[u^{\otimes2}]+iF_\infty[v] \r\} 
\bar{K}^\infty(t,t)\otimes 
D_tC^\infty_t\otimes Q_t
\bigg]
\\&=& 
E\bigg[ \int_{\bbR^d}\exp\big(iu\cdot z+iW_\infty[u]+iF_\infty[v] \big) 
\phi(z;0,C_\infty)
dz \>
\bar{K}^\infty(t,t)\otimes 
D_tC^\infty_t\otimes Q_t
\bigg]
\eeas
Therefore 
\beas
&&
r_n^{-1}i\int_0^1 \int_0^1 E\bigg[\bigg\langle K^n(s,r)[u],
e^n_s(u)\Psi_\infty(u,v)\psi_n
\half D_rC^n_s[u^{\otimes2}]
\otimes Q_s
\bigg\rangle_{\bbR^{{\sf r}}\otimes\bbR^{{\sf r}}}\bigg]\>dsdr
\\&=&
\frac{i}{2} 
E\bigg[ \int_{\bbR^d}\exp\big(iu\cdot z+iW_\infty[u]+iF_\infty[v] \big) 
\phi(z;0,C_\infty)dz 
\>
\mbox{Tr}^*  \int_0^1 
\bar{K}^\infty(t,t)\otimes 
{\coloro{\half}}{\coloroy{D}} _tC^\infty_t\otimes Q_t\>{\coloro{\mu(dt)}}
\bigg][u^{\otimes3}]
 +o(1).
\eeas
{\coloraka Moreover, $D_tC^\infty_t=0$ in this case. }

We can do with 
the terms involving either $D_rW_\infty$, $D_rF_\infty$, 
$D_rD_sW_\infty$, $D_rD_sC_\infty$, 
or $D_rD_sF_\infty$ 
in the same way to obtain
\bea\label{220505-9}
\tilde{\Phi}^{2,0}(u,v)
&:=&
\nn\lim_{n\iku\infty} r_n^{-1}\Phi^{2,0}_n(u,v)
\\&=&
\lim_{n\iku\infty} 
r_n^{-1}
E\l[L^n_1(u)\Psi_\infty(u,v)\psi_n\r]
\nn\\&=&
\frac{i}{2} 
E\bigg[ \int_{\bbR^d}\exp\big(iu\cdot z+iW_\infty[u]+iF_\infty[v] \big) 
\phi(z;0,C_\infty)
\>dz 
\>
\mbox{Tr}^*  \int_0^1 
\bar{K}^\infty(t,t)[u]\otimes 
\sigma_{t,t}(iu,iv)\>{\coloroy{\mu (dt)}}
\bigg]
\nn\\&=&
E\bigg[ \int_{\bbR^d}\exp\big(iu\cdot z+iF_\infty[v] \big) 
\phi(z;W_\infty,C_\infty)\>dz 
\>\bar{\sigma}(iu,iv)
\bigg], 
\eea
where $\sigma_{t,t}(iu,iv)=\lim_{s\up t}\sigma_{t,s}(iu,iv)$ and 
\bea\label{241003-1} 
\bar{\sigma}(iu,iv)
&=&
\half\mbox{Tr}^*  \int_0^1 
\bar{K}^\infty(t,t)[iu]\otimes 
\sigma_{t,t}(iu,iv)\>{\coloroy{\mu (dt)}}. 
\eea
{\colb The random symbol $\sigma_{t,t,}(iu,iv)$ can be expressed formally as 
\beas 
\sigma_{t,t}(iu,iv)
&=&
\Big(iD_tW_\infty[u]-\half D_tC_\infty[u^{\otimes2}]+iD_tF_\infty[v] \Big)
\otimes 
\Big(iD_tW_\infty[u]-\half D_tC_\infty[u^{\otimes2}]+iD_tF_\infty[v] \Big)
\\&&
+
\Big(iD_tD_tW_\infty[u]-\half D_tD_tC_\infty[u^{\otimes2}]+iD_tD_tF_\infty[v] 
\Big). 
\eeas
}

{\colorr{
By tracing the derivation of the above limit 
with the Leibniz rule, 
we obtain in a similar manner 
\beas 
\tilde{\Phi}^{2,\alpha}(u,v)
&:=&
\lim_{n\iku\infty} r_n^{-1}
\partialbs^\alpha\Phi^{2,0}_n(u,v)
\\&=&
\lim_{n\iku\infty} 
r_n^{-1}
\partialbs^\alpha
E\l[L^n_1(u)\Psi_\infty(u,v)\psi_n\r]
\\&=&
\frac{i}{2} 
\partialbs^\alpha
E\bigg[ \int_{\bbR^d}\exp\big(iu\cdot z+iW_\infty[u]+iF_\infty[v] \big) 
\phi(z;0,C_\infty)
\>dz 
\>
\mbox{Tr}^*  \int_0^1 
\bar{K}^\infty(t,t)[u]\otimes 
\sigma_{t,t}(iu,iv)\>dt
\bigg]
\\&=&
\partialbs^\alpha
E\bigg[ \int_{\bbR^d}\exp\big(iu\cdot z+iF_\infty[v] \big) 
\phi(z;W_\infty,C_\infty)\>dz 
\>\bar{\sigma}(iu,iv)
\bigg]. 
\eeas

{\coloroy{
Applying the integration-by-parts formula at (\ref{220505-3})  {\colb or its derivatives}
in the same way  
as we reached (\ref{220505-1}), we obtain 
\bea\label{220505-10}
\sup_n\sup_{(u,v)\in\Lambda^0_n(\check{d},q)}|(u,v)|^{\check{d}+1-{\colb \ep'}}
r_n^{-1}|\Phi^{2,\alpha}_n(u,v)|
&<&\infty. 
\eea
}}

\subsection{Asymptotic expansion of the double stochastic integral}\label{221003-1}
{\colory{
We are now on the point of presenting our results with the aid of  
the preceding subsection. 
{\coloraka The density $p_n$ is given by (\ref{220829-2}), (\ref{220829-1}), (\ref{240218-1}) and (\ref{241003-1}). }
%
\begin{theorem}\label{210530-1}
Suppose that Conditions {\rm [A1], [A2]$^\natural$} and {\rm [A3]$^\natural$} are fulfilled. 
Let $M,\gamma\in(0,\infty)$ and $\theta\in(0,1)$ be arbitrary numbers. 
Then, for some constant $C_1=C(M,\gamma,\theta)$,  {\rm (\ref{210926-5})} holds as $n\to\infty$. 
\end{theorem}
\proof 
Inequality (\ref{210314-1}) has been verified in the present situation by (\ref{220505-10}). 
We can verify 
{\colb [B1], [B2]$_\ell$, [B3] and [B4]$_{\ell,{\mathfrak m},{\mathfrak n}}$ 
}
for $({\mathfrak m},{\mathfrak n})=(5,2)$  
to apply 
Theorem \ref{210215-1} (b). 
[Recall that now $\ell$ 
is different from  ``$\ell$'' in 
 Theorem \ref{210215-1}.]  
We apply \cite{Jacod1997} for 
[B1] (ii){\coloraka, however} 
we still need the joint convergence assumption [A2]$^\natural$ (iv). 
\qed 
\ \\

We shall present a version of Theorem \ref{210530-1}. 

\bd
\item[{[A1]}$^\flat$] 
Conditions in [A1] hold except for {\rm (iii)}. 
\ed

$s_n:\Omega\to\bbR$ is a positive functional.

\bd
\item[{[A2]}] 
Condition [A2]$^\natural$ holds, replacing 
its (ii) by 
\bd
\item[(ii)]
$F_n\in\bbD_{\ell+1,\infty}(\bbR^{d_1})$, 
$W_n\in\bbD_{{\coloro{\ell+1}},\infty}(\bbR^d)$, 
$N_n\in\bbD_{{\coloro{\ell+1}},\infty}(\bbR^d)$  
and $s_n\in\bbD_{\ell,\infty}(\bbR)$. 
Moreover, 
\beas 
\sup_{{\coloraka n\in\bbN}}\Big\{ 
{\coloraka \|\dotc_n\|_{\ell,p}+}
\|\dotw_n\|_{{\coloro{\ell+1}},p}
+\|\dotf_n\|_{\ell+1,p}+\|N_n\|_{{\coloro{\ell+1}},p}+\|s_n\|_{\ell,p} 
\Big\} <\infty
\eeas
for every $p\geq2$. 
\ed
\ed

\begin{en-text}
\bd
\item[{[A2]}] 
{\bf (i)}  
$F_\infty\in\bbD_{\ell+1,\infty}(\bbR^{d_1})$ 
and 
$W_\infty\in\bbD_{{\coloro{\ell+1}},\infty}(\bbR^d)$
. 

\bd
\item[(ii)]
$F_n\in\bbD_{\ell+1,\infty}(\bbR^{d_1})$, 
$W_n\in\bbD_{{\coloro{\ell+1}},\infty}(\bbR^d)$, 
$N_n\in\bbD_{{\coloro{\ell+1}},\infty}(\bbR^d)$  
and $s_n\in\bbD_{\ell,\infty}(\bbR)$. 
Moreover, 
\beas 
\sup\Big\{ 
\|\dotw_n\|_{{\coloro{\ell+1}},p}
+\|\dotf_n\|_{\ell+1,p}+\|N_n\|_{{\coloro{\ell+1}},p}+\|s_n\|_{\ell,p} 
\Big\} <\infty.
\eeas
for every $p\geq2$. 

\item[(iii)] 
$\tilde{C}_\infty^{j,k}$, $\tilde{W}_\infty^j$, $\tilde{N}_\infty^j$ $(j,k=1,...,d)$ 
and $\tilde{F}_\infty^l$ $(l=1,...,d_1)$ are in $\calc(2[(d_1+3)/2],2)$. 

\item[(iv)] 
$({\coloro{M^n_\cdot}},N_n,\dotc_n,\dotw_n,\dotf_n)
\iku^{d_s(\calf)}
(M^\infty_\cdot,N_\infty,\dotc_\infty,\dotw_\infty,\dotf_\infty)$. 

\item[(v)] 
For $G=W_\infty$ and $F_\infty$, 
\beas 
ess. \sup_{r_1,...,r_k\in(0,1)} 
\|D_{r_1,...,r_k}G\|_p 
&<&\infty
\eeas
for every $p\in[2,\infty)$ and $k\leq d_1+6$. 
\footnote{$2[(2(={\mathfrak n}+d_1+2)/2]+2$ (of $D_rD_sW_\infty$ etc. )}
Moreover, $r\mapsto D_rG$ and $(r,s)\mapsto D_{r,s}G$ are 
continuous a.e. with respect to the Lebesgue measures 
[or in $L^{2+}$ and $L^{1+}$-sense, respectively]. 
\ed
\ed
\end{en-text}

The nondegeneracy of $(M^n_t+W_\infty,F_\infty)$ will be necessary. 
\bd
\item[{[A3]}] 
{\bf (i)} 
There exist a sequence $(t_n)_{n\in\bbN}$ in {\colorr{$[0,1]$ with $\sup_nt_n<1$}}  
such that 
%
%

$
\sup_{t\geq t_n}
P\big[
\det\sigma_{(M^n_t+W_\infty,F_\infty)}< s_n\big] =O(r_n^{\colb \nu})
$
as $n\to\infty$ for some ${\colb \nu>\ell/3}$. 

\begin{description}
\item[(ii)] 
For every $p\geq2$, 
$
\limsup_{n\to\infty} E[s_n^{-p}]<\infty  
$
{\coloraka .}
%
%



\end{description}
\ed

\ \\

\begin{theorem}
Suppose that Conditions {\rm [A1]}$^\flat$, {\rm[A2]} and {\rm [A3]} are fulfilled. 
Then for any positive numbers $M$ and $\gamma$, 
\beas
\sup_{f\in\cale(M,\gamma)}\Delta_n(f)
&=&
o(r_n)
\eeas
as $n\to\infty$. 
\end{theorem}
\proof 
As in the proof of Theorem \ref{210530-1}, 
we will verify 
{\colb [B1], [B2]$_{{\coloraka \ell'}}$, [B3] and [B4]$_{{\coloraka \ell'},{\mathfrak m},{\mathfrak n}}$} 
for $(\ell',{\mathfrak m},{\mathfrak n})=(\check{d}+3,5,2)$  
to apply 
Theorem \ref{210215-1} (b). 
Define $\xi_n$ by 
\beas 
\xi_n &=& 
10^{-1}r_n^{-2c}|C_n-C_\infty|^2
+
2\Big[1+4\Delta_{(M_n+W_\infty,F_\infty)}s_n^{-1}\Big]^{-1}
\\&&
+L^*\int_{[0,1]^2}
\bigg(
\frac{|C^n_t-C^\infty_t-C^n_s+C^\infty_s| r_n^{-2q}}{|t-s|^{3/8}}\bigg)^8 dtds, 
\eeas
where $c$ is a constant given in the proof of Theorem \ref{210926-6} and 
$L^*$ is a sufficiently large constant. 
{\coloraka Here we can choose a number $q$ that is smaller than the given $q$.} 
 Now the rest is to verify (\ref{210314-1}) for ``$\ell$''$=\check{d}+3$, 
 that is, (\ref{220505-10}). 
 There were two steps to reach (\ref{220505-10}): 
 (\ref{220505-1}) and the argument just before (\ref{220505-10}) concerning (\ref{220505-3}). 
 The reasoning is quite the same in those cases, so we will show (\ref{220505-1}). 
 However we need to do it under [A1]$^\flat$, [A2] and [A3] this time, 
 not under $[A1]$, [A2]$^\natural$ and [A3]$^\natural$. 
 Obviously we have [A2]$^\natural$. 
 Condition [A1] (iii) is satisfied due to the definition of the above $\xi_n$ with suitable $L^*$. 

In order to estimate (\ref{220504-1}) once again under the present assumptions, 
we can follow (a)-(g) in Section \ref{221024-4}. 
Let 
\beas 
\psi_{n,s}&=&
\psi\bigg(2\Big[1+4\Delta_{(M^n_s+W_\infty,F_\infty)}s_n^{-1}\Big]^{-1}\bigg)
\eeas
for $s\geq t_n$. 
Then the integrand of (\ref{220504-1}) is decomposed as 
\beas 
&&
E\bigg[\bigg\langle K^n(s,r)[u],
e^n_s(u)\Psi_\infty(u,v)\psi_n
iD_rM^n_s[u]\otimes Q_s
\bigg\rangle_{\bbR^{{\sf r}}\otimes\bbR^{{\sf r}}}\bigg]
\\&=&
E\bigg[\bigg\langle K^n(s,r)[u],
e^n_s(u)\Psi_\infty(u,v)\psi_n\psi_{n,s}
iD_rM^n_s[u]\otimes Q_s
\bigg\rangle_{\bbR^{{\sf r}}\otimes\bbR^{{\sf r}}}\bigg]
+
R_{n}(s,r)
\eeas
with 
\beas 
|R_{n}(s,r)|
&\leq& 
C(p)r_n^{-{\mathfrak k}q}
\>\sup_{s'}\|1-\psi_{n,s'}\|_p\>
\>r_n^{-1}\sum_j1_{(t_{j-1},t_j]}(s)1_{(t_{j-1},s]}(r)
\eeas
for all $n$, $s$ and restricted $(u,v)$, where ${\mathfrak k}=4$ and 
$C(p)$ is a constant independent of them 
but on $p\in(1,\infty)$. 
Taking a small $q>1/3$, we have 
{\colb $\nu-(\check{d}+1+{\mathfrak k})q>0$} 
(even for ${\mathfrak k}=5$). 
We can apply 
the integration-by-parts formula repeatedly as in (h) of Section \ref{221024-4} 
with truncation by $\psi_n\psi_{n,s}$ instead of $\psi_n$, and then 
(\ref{220505-1}) follows from [A3](i) by choosing a small $p>1$. 
It was the estimate for the terms concerning $D_rM^n_s$. 
We can do the same kind estimate for the terms involving $D_rC^n_s$ {\coloraka and others} 
to finally obtain (\ref{220505-10}) while ${\mathfrak k}=5$ in this case. 
\qed 
}}

\section{Quadratic form of a Wiener process}\label{220826-1}
\subsection{Asymptotic expansion}\label{240324-2} 

The quadratic form of the increments of a diffusion process 
with a strongly predictable kernel 
plays a central role in the inference for diffusion coefficients. 
The case of the Wiener process shows 
what is most essential in the analysis. 
Let $a\in C^\infty_\up(\bbR)$, the set of smooth functions with all derivatives at most polynomial growth. 
We write 
$1_j=1_{(t_{j-1},t_j]}$, $t_j=j/n$. 
Let 
\beas 
M^n_t &=& \sqrt{n}\sum_j
2a(w_{t_{j-1}})\int_{t_{j-1}\wedge t}^{t_j\wedge t}\int_{t_{j-1}}^s dw_rdw_s.
\eeas
In this section, $w=(w_t)_{t\in[0,1]}$ denotes a standard Wiener process starting at $w_0$. 
Then 
\beas 
C^n_t &=& n\int_0^t \sum_j 1_j(s) \>4a(w_{t_{j-1}})^2
\Big(\int_{t_{j-1}}^sdw_r\Big)^2 ds
\eeas
and we see
\beas 
C^n_t &=& 
\sum_j\int_{t_{j-1}}^{t_j}1_{[0,t]}(s)\>4a(w_{t_{j-1}})^2
n\Big\{\Big(\int_{t_{j-1}}^sdw_r\Big)^2-(s-t_{j-1})\Big\}ds
\\&&
+
\frac{2}{n}\sum_{j:t_j\leq t}a(w_{t_{j-1}})^2
+O_p\Big(\frac{1}{n}\Big)
\\&\to^p&
2\int_0^t a(w_s)^2ds=C^\infty_t. 
\eeas
In this example,we will consider 
\begin{comment}
$W_\infty=0$ and 
\beas 
F_n &=& F_\infty =2\int_0^1 a(X_s)^2ds \equiv C^\infty_1.
\eeas
\end{comment}
%
\beas 
F_n &=& \frac{2}{n}\sum_j a(w_{t_{j-1}})^2. 
\eeas
Obviously  
\beas 
F_n &\to^p& 2\int_0^1 a(w_s)^2ds=C^\infty_1\equiv C_\infty. 
\eeas
This setting is natural in estimation of the diffusion coefficient (volatility), where 
$M^n$ becomes the ``deviation'' of the estimator for the cumulative variance 
of the system and $F_n$ is an estimator of the ``asymptotic variance''. 

{\coloroy{Let $\alpha(x)=a(x)^2$. 
In order to obtain asymptotic expansion, 
we will assume the nondegeneracy condition
\bd
\item[[HW1\!\!]]\hspace{5mm}
$\inf_{x\in\bbR} |a(x)| > 0.$
\ed
In application to statistical estimation of volatility, 
this condition is translated as the uniform ellipticity of the diffusion process. 
For the nondegeneracy of $F_\infty$, we need the following condition. 
\bd 
\item[[HW2\!\!]]\hspace{5mm}
$\sum_{i=1}^\infty|\partial^i\alpha(w_0)|>0$, 
where $w_0$ is the initial value of $w$. 
\ed
The initial value $w_0$ can be random but has a compact support.

Set 
\beas 
{\mathfrak C}_0 &=& 2\int_0^1 a(w_s)^2ds 
= C^\infty_1\equiv C_\infty = F_\infty,
\\
{\mathfrak C}_1 &=& 
{\coloro{\frac{2}{3}\int_0^1a(w_s)^3ds\Big( \int_0^1a(w_s)^2ds\Big)^{-1}}}
\\
{\mathfrak C}_2 &=& \int_0^1a(w_t)\Big(\int_t^1aa'(w_\nu)d\nu\Big)^2dt
\\
{\mathfrak C}_3 &=& \int_0^1a(w_t)\int_t^1
\{aa''+(a')^2\}(w_\nu)d\nu\> dt. 
\eeas
Then the random symbol $\sigma$ is given by 
\beas 
\sigma(z,iu,iv) 
&=& 
{\coloro{{\mathfrak C}_1 z (iu)^2 }}
+ {\mathfrak C}_2 iu\Big(-2u^2+4iv\Big)^2
+ {\mathfrak C}_3 iu\Big(-2u^2+4iv\Big). 
\eeas
The asymptotic expansion is done by
\beas 
p_n(z,x) 
&=&
\phi(z;0,x)E[\delta_x({\mathfrak C}_0)] 
+\frac{1}{\sqrt{n}} \dot{p}(z,x),
\eeas
where 
\beas 
\dot{p}(z,x)
&=&
E\Big[\sigma(z,\partial_z,\partial_x)^*
\Big\{\phi(z;0,C_\infty)\delta_x(C_\infty)\Big\}\Big]
\\&=&
{\coloro{
E[{\mathfrak C}_1\delta_x({\mathfrak C}_0)]  \partial_z^2\Big\{z\phi(z;0,x)\Big\}
}}
\\&&
-\partial_z\Big(2\partial_z^2-4\partial_x\Big)^2 
\Big\{E[{\mathfrak C}_2\delta_x({\mathfrak C}_0)]\phi(z;0,x)\Big\}
\\&&
-\partial_z\Big(2\partial_z^2-4\partial_x\Big) 
\Big\{E[{\mathfrak C}_3\delta_x({\mathfrak C}_0)]\phi(z;0,x)\Big\}
\\&=&
{\colb 
E[{\mathfrak C}_1\delta_x({\mathfrak C}_0)]\partial_z^2\big\{z\phi(z;0,x)\big\}
-16\big\{\partial_x^2E[{\mathfrak C}_2\delta_x({\mathfrak C}_0)]\big\}\partial_z\phi(z;0,x)
}
\\&&
{\colb 
+4\big\{\partial_xE[{\mathfrak C}_3\delta_x({\mathfrak C}_0)]\big\}\partial_z\phi(z;0,x)
}
\\&=:&
p_1(z,x)+p_2(z,x)+p_3(z,x). 
\eeas

\begin{remark}\rm 
In order to obtain rough evaluation of the terms involving $\delta_x$, 
we may apply the IBP formula, the kernel method, or other methods. 
\end{remark}

\begin{theorem}\label{220820-3}
Suppose that {\rm [HW1]} and {\rm[HW2]} are satisfied. 
Then for any positive numbers $M$ and $\gamma$, 
\beas 
\sup_{f\in\cale(M,\gamma)} \bigg|
E\big[f(M^n_1,F_n)\big]-\int_{\bbR^2}f(z,x)p_n(z,x)dzdx\bigg|
&=&
o\bigg(\frac{1}{\sqrt{n}}\bigg)
\eeas
as $n\to\infty$, 
where $\cale(M,\gamma)$ is the set of measurable functions $f:\bbR^2\to\bbR$ 
satisfying $|f(z,x)|\leq M(1+|z|+|x|)^\gamma$ for all $z,x\in\bbR$.  
\end{theorem}

\subsection{Proof}
\subsubsection{Representation of the limit variables}\label{220826-5}

}}

The variable 
$M^n_1$ admits the representation 
\beas 
M^n_1 &=& 
n^{-\half} \sum_{j:t_j\leq 1}{\coloro{a}}(w_{t_{j-1}})
\Big((\sqrt{n}\Delta_jw)^2-1\Big),
\eeas
where $\Delta_jw=w_{t_j}-w_{t_{j-1}}$. 
%
For $\dotc\>\!\!^n_t:=\sqrt{n}\Big(C^n_t - C^\infty_t\Big)$, we have 
\beas 
\dotc\>\!\!^n_t
&=& 
\sum_{j:t_j\leq t}\int_{t_{j-1}}^{t_j}1_{[0,t]}(s)\>4n\sqrt{n}a(w_{t_{j-1}})^2
\Big\{\Big(\int_{t_{j-1}}^sdw_r\Big)^2-(s-t_{j-1})\Big\}ds
\\&&
-
2\sqrt{n}\sum_{j:t_j\leq t} \int_{t_{j-1}}^{t_j} \Big(a(w_s)^2-a(w_{t_{j-1}})^2\Big)ds
+O_p\Big(\frac{1}{\sqrt{n}}\Big)
\\&=& 
\sum_{j:t_j\leq t}\int_{t_{j-1}}^{t_j}1_{[0,t]}(s)\>4n\sqrt{n}a(w_{t_{j-1}})^2
\Big\{\Big(\int_{t_{j-1}}^sdw_r\Big)^2-(s-t_{j-1})\Big\}ds
\\&&
-
2\sqrt{n}\sum_{j:t_j\leq t} \int_{t_{j-1}}^{t_j}
2a(w_{t_{j-1}})a'(w_{t_{j-1}})(w_s-w_{t_{j-1}})ds
+O_p\Big(\frac{1}{\sqrt{n}}\Big) 
\\&=& 
\sum_{j:t_j\leq t}\int_{t_{j-1}}^{t_j}1_{[0,t]}(s)\>4n\sqrt{n}a(w_{t_{j-1}})^2
\Big\{\Big(\int_{t_{j-1}}^sdw_r\Big)^2-(s-t_{j-1})\Big\}ds
\\&&
+O_p\Big(\frac{1}{\sqrt{n}}\Big). 
\eeas
Here we know that 
the supremum of ``$O_p(n^{-1/2})$'' in $t\in[0,1]$ is of $O_p(n^{-1/2})$. 
With the same argument, we have 
\beas 
\sqrt{n}(F_n-F_\infty) &\to^p& 0=\dotf_\infty. 
\eeas 

We should specify the distribution of 
\beas 
(M_\infty, \dotc_\infty)
\eeas
For computations, we will use the discrete filtration 
$\F^n=(\bar{\calf}^n_t)_{t\in  [0,1]}$ with 
$\bar{\calf}^n_t=\calf_{[nt]/n}$. 
\def\dla{\langle\!\langle}
\def\dra{\rangle\!\rangle}
The bracket for $\F^n$ is denoted by 
$\dla\cdot\dra$; though it depends on $n$, we suppress it from notation. 
Let 
$H_1(x)=x$ and $H_2(x)=(x^2-1)/\sqrt{2}$. 
Denote $\Delta_j w=w_{t_j}-w_{t_{j-1}}$, which depends on $n$ as well as $j$. 
The discrete version of $M^n$ is given by 
\beas 
\bar{M}^{2,n}_t 
&=&
\frac{1}{\sqrt{n}}\sum_{j:t_j\leq t}\sqrt{2}a(w_{t_{j-1}}) H_2(\sqrt{n}\Delta_j w). 
\eeas
The principal part of $\dotc\>\!\!^n$ is $\F^n$-martingale 
\beas 
\bar{M}^{\xi,n}_t 
&=&
\sum_{j:{\colb t_j}\leq t}\int_{t_{j-1}}^{t_j} 4n\sqrt{n}a(w_{t_{j-1}})^2
\Big\{\Big(\int_{t_{j-1}}^sdw_r\Big)^2-(s-t_{j-1})\Big\}ds
\eeas
The discrete version of $w$ is denoted by $\bar{w}^n_t=w_{[nt]/n}$.  
Then we have 
\beas 
\dla \bar{w}^n,\bar{w}^n \dra_t &=& \frac{[nt]}{n}\to t
\\
\dla \bar{M}^{2,n},\bar{M}^{2,n} \dra_t 
&=& 
\frac{2}{n}\sum_{j:{\colb t_j}\leq t} a(w_{t_{j-1}})^2 
\to^{ucp} 
2\int_0^t a(w_s)^2 ds
\\
\dla \bar{M}^{\xi,n},\bar{M}^{\xi,n} \dra_t 
&=&
\frac{16}{3n}\sum_{j:{\colb t_j}\leq t} a(w_{t_{j-1}})^4
\to^{ucp} 
\frac{16}{3}\int_0^t a(w_s)^4ds 
\\
\dla \bar{w}^n,\bar{M}^{k,n} \dra_t &=& 0
\sskip(k=2,\xi)
\\
\dla \bar{M}^{2,n},\bar{M}^{\xi,n} \dra_t &=& 
\frac{{\coloro{8}}}{3n}\sum_{j:{\colb t_j}\leq t}a_{t_{j-1}}^3
\to^{ucp} 
\frac{{\coloro{8}}}{3}\int_0^t a(w_s)^3 ds. 
\eeas
Obviously, 
the orthogonality between those martingales and any bounded martingales orthogonal to $w$ holds. 
Therefore, we obtain the following stable convergence with a representation of the limit: 
\beas 
(\bar{M}^{2,n}, \bar{M}^{\xi,n})
&\to^{d_s(\calf)}& 
\Big( \int_0^\cdot \sqrt{2}a(w_s)dB_s,
\int_0^\cdot \frac{{\coloro{4}}\sqrt{2}}{3}a(w_s)^2dB_s
+\int_0^\cdot \frac{{\coloro{4}}}{3}a(w_s)^2dB'_s\Big),
\eeas
where $(B,B')$ is a two-dimensional standard Wiener process, 
independent of $\calf$,  
defined on the extension $\bar{\Omega}$. 
In particular, 
\beas 
(M_\infty,\dotc_\infty)
&=^d&
\Big( \int_0^1 \sqrt{2}a(w_s)dB_s,
\int_0^1 \frac{{\coloro{4}}\sqrt{2}}{3}a(w_s)^2dB_s
+\int_0^1 \frac{{\coloro{4}}}{3}a(w_s)^2dB'_s\Big), 
\eeas
so that 
\beas 
\check{C}_\infty(\omega,M_\infty)
=\bbE[\dotc_\infty|\check{\calf}]
&=&
{\coloro{\frac{4}{3}\int_0^1a(w_s)^3ds\Big( \int_0^1a(w_s)^2ds\Big)^{-1}}}\> M_\infty. 
\eeas
Thus we have 
\beas 
\tilde{C}_\infty(z)
&=&
\check{C}_\infty(\omega,z-W_\infty)
=
\check{C}_\infty(\omega,z)
\\&=&
{\coloro{\frac{4z}{3}\int_0^1a(w_s)^3ds\Big( \int_0^1a(w_s)^2ds\Big)^{-1}}}. 
\eeas
By definition of the variable in this example, 
we have $\tilde{F}_\infty(\omega,z)=0$ and $\tilde{N}_\infty(\omega,z)=0$. 
{From} above computations, 
the random symbol $\underline{\sigma}(z,iu,iv)$ is given by 
\beas 
\underline{\sigma}(z,iu,iv)
&=&
{\coloro{\frac{2z}{3}\int_0^1a(w_s)^3ds\Big( \int_0^1a(w_s)^2ds\Big)^{-1}}} \>(iu)^2. 
\eeas

In order to obtain the random symbol $\sigma_{s,r}(iu,iv)$, we note 
\beas 
D_rC^\infty_s 
&=& 
4\int_r^s aa'(w_\nu)d\nu\>1_{\{r\leq s\}}, 
\\
D_sC_\infty
&=&
D_sC^\infty_1
=
D_sF_\infty
=
4\int_s^1 aa'(w_\nu)d\nu\>1_{\{s\leq1\}}, 
\\
D_rD_sC_\infty
&=&
D_rD_sF_\infty
=
4\int_{s\vee r}^1 \{aa''+(a')^2\}(w_\nu)d\nu. 
\eeas
Therefore, 
\beas 
\sigma_{s,r}
&=&
2 u^2\Big(-2u^2+4iv\Big)
\int_r^s aa'(w_\nu)d\nu 
\int_s^1 aa'(w_\nu)d\nu
\\&&
+
\Big(-2 u^2+4iv\Big)^2\int_r^1 aa'(w_\nu)d\nu
\int_s^1 aa'(w_\nu)d\nu
\\&&
+
\Big(-2u^2+4iv\Big)
\int_{s\vee r}^1 \{aa''+(a')^2\}(w_\nu)d\nu
\eeas
and 
\beas
\sigma_{t,t}(iu,iv)
&=&
\Big(-2 u^2+4iv\Big)^2 \Big(\int_t^1 aa'(w_\nu)d\nu\Big)^2
\\&&
+
\Big(-2u^2+4iv\Big)
\int_t^1 \{aa''+(a')^2\}(w_\nu)d\nu. 
\eeas
Thus
\beas 
\bar{\sigma}(iu,iv)
&=&
\int_0^1 a(w_t)iu\>\sigma_{t,t}(iu,iv)\>dt
\\&=&
iu\Big(-2 u^2+4iv\Big)^2 
\int_0^1 a(w_t) \Big(\int_t^1 aa'(w_\nu)d\nu\Big)^2\>dt
\\&&
+
iu\Big(-2u^2+4iv\Big)
\int_0^1 a(w_t)
\int_t^1 \{aa''+(a')^2\}(w_\nu)d\nu\>dt . 
\eeas

\begin{en-text}
\begin{remark}\rm 
{\coloroy{When the function $a$ is odd,}}
the conditional expectation 
of $\bar{\sigma}(iu,iv)$ given $F_\infty$ vanishes by symmetry. 
\end{remark}
\end{en-text}

\begin{remark}\rm 
From computational point of view, 
only rough simulation for  the conditional expectation 
of $\bar{\sigma}(iu,iv)$ given $F_\infty$ will be necessary 
to obtain the second-order term. 
This is called the {\it Hybrid I method}. 
\end{remark}

{\coloroy{Thanks to the nondegeneracy of $a$, 
it is easy to verify [A1](ii); otherwise, we would be involved in a tedious large deviation argument 
which we did not want to pursuit here. 
Other conditions in [A1] are also easy to prove.  }}
\begin{en-text}
We will sketch a way to verify [A1](ii); other conditions in [A1] are easy to prove. 
Suppose that $\liminf_{|x|\to\infty}|a(x)|>0$. 
[Otherwise, there is an example for which $a(w_t)=0$ for all $t\in[s,1]$ with positive probability. 
It is the case when $a(x)=0$ outside of the unit disk. ]
For $n\in\bbN$ and $s\in(0,1-n^{-1/2}]$, 
\beas 
P\bigg[
(1-s)^{-(1+\eta)}\int_s^1a(w_t)^2dt<n^{-c_1}\bigg]
&\leq&
P\bigg[
\int_s^{s+n^{-1}}a(w_t)^2dt<n^{-c_1}\bigg]
\\
&\leq&
P\bigg[\int_s^{s+n^{-1}}a(w_t)^2dt<n^{-c_1},\> |a(w_s)|>n^{-\half c_1}\bigg]
\\&&
+
P\bigg[\int_s^{s+n^{-1}}\tilde{a}(w_t)^2dt<n^{-c_1}\bigg],
\eeas
where in the last term we replaced $a$ by a bounded $\tilde{a}$ with all derivatives bounded 
due to the compactness of the set $\{x\in\bbR;|a(x)|\leq n^{-\half c_1}\}$. 
The first term on the right-hand side is exponentially small, and 
the second term turns out to be so if we apply the Key lemma by Kusuoka and Stroock 
(cf. Ikeda and Watanabe p.398). 
For example, 
it is sufficient if 
\bea\label{H-condition} 
\inf_{x\in\bbR}\sum_{i=0}^k |\partial_x^i\alpha(x)|>0
\eea
for some $k\in\bbN$, where $\alpha(x)=a(x)^2$. 
The derivatives of $\alpha$ come from a sequential It\^o expressions: 
\beas 
\xi(t) &=& \int_s^t \alpha(w_r) dr,
\\
\alpha(w_t)&=&\alpha(w_s)+\int_s^t\alpha'(w_r)\circ dw_r
=\alpha(w_s)+\int_s^t\alpha'(w_r)dw_r+\int_s^t\half \alpha''(w_r)dr,
\\
\alpha'(w_t)&=&\alpha'(w_s)+\int_s^t\alpha''(w_r)\circ dw_r
=\alpha'(w_s)+\int_s^t\alpha''(w_r)dw_r+\int_s^t\half \alpha'''(w_r)dr,
\\
&&\cdots.
\eeas
\end{en-text}

\subsubsection{Nondegeneracy}

Here we will briefly discuss the nondegeneracy in 
Malliavin's sense. We have 
\beas 
D_rM^n_t &=& 
\sqrt{n}\sum_{j=1}^n
2{\coloroy{a(w_{t_{j-1}})}}
\bigg(\int_{t_{j-1}\wedge t}^{t_j\wedge t}dw_s \bigg) 1_{(t_{j-1}\wedge t,t_j\wedge t)}(r)
\\&&
+\sqrt{n}\sum_{j=1}^na'(w_{t_{j-1}}){\coloroy{1_{(0,t_{j-1}\wedge t]}(r)}}
\bigg((\int_{t_{j-1}\wedge t}^{t_j\wedge t}dw_s)^2-(t_{j-1}\wedge t-t_{j-1}\wedge t)\bigg)
\\&=:&
D_1(n,t)_r+D_2(n,t)_r. 
\eeas

{\coloroy{
Let 
$\eta_j(t) = \sqrt{n}(w(t_j\wedge t)-w(t_{j-1}\wedge t))$ and 
\beas 
\xi_j(t) &=& 
n\bigg((w(t_j\wedge t)-w(t_{j-1}\wedge t))^2-(t_j\wedge t - t_{j-1}\wedge t)\bigg). 
\eeas
Then 
\beas 
D_2(n,t)_r 
&=& 
n^{-\half}\sum_{j=1}^na'(w_{t_{j-1}}){\coloroy{1_{(0,t_{j-1}\wedge t]}(r)}}\xi_j(t)
\\&=&
n^{-\half}\sum_{j=1}^{n-1}\bigg(\sum_{k=j+1}^n a'(w_{t_{k-1}})\xi_k(t)\bigg)1_{I_j(t)}(r),
\eeas
where $I_j(t)=(t_{j-1}\wedge t,t_j\wedge t]$.  
{\coloraka For a while, we assume $t\in\{t_j\}_j$.} 
Hence 
\beas 
\sigma_{11}(n,t):=\sigma_{M^n_t} 
&=&
\frac{1}{n}
\sum_{j=1}^n \bigg[2a({\coloraka w(t_{j-1})})\eta_j(t)
+\frac{1}{\sqrt{n}}\sum_{k=j+1}^na'(w(t_{k-1}))\xi_k(t)\bigg]^2,
\eeas
where we read $\sum_{k=n+1}^n...=0$. 
Moreover, 
\beas 
\sigma_{12}(n,t)
:=\langle M^n_t,F_\infty\rangle_H
&=& 
\sum_{j=1}^n \bigg[2a(w_{t_{j-1}})\eta_j(t)
+\frac{1}{\sqrt{n}}\sum_{k=j+1}^na'(w(t_{k-1}))\xi_k(t)\bigg]
\cdot \int_{I_j(t)}4\int_r^1aa'(w_s)ds \>dr
\eeas
and 
\beas 
\sigma_{22}(t)
&=&
\int_0^t\bigg[4\int_r^1aa'(w_s)ds\bigg]^2dr.
\eeas
The Malliavin covariance matrix of $(M^n_t,F_\infty)$ is then given by 
\beas
\sigma_{(M^n_t,F_\infty)}
=\l[
\begin{array}{cc} 
\sigma_{11}(n,t) & \sigma_{12}(n,t) \\ \sigma_{12}(n,t) & \sigma_{22}(1)
\end{array} 
\r]. 
\eeas
Let 
\beas
\sigma(n,t)&:=&
\l[
\begin{array}{cc} 
\sigma_{11}(n,t) & \sigma_{12}(n,t) \\ \sigma_{12}(n,t) & \sigma_{22}(t)
\end{array} 
\r]. 
\eeas
Let 
\beas 
\tilde{\sigma}_{11}(n,t)
&=&
\frac{1}{n}\sum_{j=1}^n \bigg[2a(w_{t_{j-1}})\eta_j(t)\bigg]^2
+\frac{1}{n}\sum_{j=1}^n \bigg[\frac{1}{\sqrt{n}}\sum_{k=j+1}^na'(w(t_{k-1}))\xi_k(t)\bigg]^2,
\eeas
\beas 
\tilde{\sigma}_{12}(n,t)
&=& 
\sum_{j=1}^n \bigg[
\frac{1}{\sqrt{n}}\sum_{k=j+1}^na'(w(t_{k-1}))\xi_k(t)\bigg]
\cdot \int_{I_j(t)}4\int_r^1aa'(w_s)ds \>dr
\eeas
and 
\beas 
\tilde{\sigma}(n,t)
=\l[
\begin{array}{cc} 
\tilde{\sigma}_{11}(n,t) & \tilde{\sigma}_{12}(n,t) \\ \tilde{\sigma}_{12}(n,t) & \sigma_{22}(t)
\end{array} 
\r]. 
\eeas
We shall show 
\bea\label{220820-1} 
\big\| \sigma(n,t) - \tilde{\sigma}(n,t) \big\|_p &=& O(n^{-\half})
\eea
for every $p>1$ and 
{\coloraka uniformly in $t$.} 

Let $\cali$ denote the set of sequences $J^{(\nu)}=(J^{(\nu)}_{n,j})$ of 
multiple It\^o stochastic integrals  
taking the form
\beas 
J^{(\nu)}_{n,j} 
&=& n^{\frac{\nu}{2}}
\int_{t_{j-1}}^{t_j}dw_{s_1}a^{(\nu)}_{n,j,1}(s_1)
\int_{t_{j-1}}^{s_1}dw_{s_2}a^{(\nu)}_{n,j,2}(s_2)
\int_{t_{j-1}}^{s_2}
\cdots
\int_{t_{j-1}}^{s_\nu}dw_{s_\nu}a^{(\nu)}_{n,j,\nu}(s_\nu),
\eeas
where 
$\{a^{(\nu)}_{n,j,i};\>i=1,...,\nu$, $j=1,...,n$, $n\in\bbN\}$ is a family of progressively measurable processes.   
In the following lemma, 
$J^{(\nu_1)}_{n,j_1}\cdots J^{(\nu_m)}_{n,j_m}$, 
$J^{(\mu_1)}_{n,k_1}\cdots J^{(\mu_m)}_{n,k_m}$ 
are in $\cali$, and each of them has $a^{(*)}_{n,j,i}$ which may possibly differ from 
those of other indices $\nu$'s and $\mu$'s even if 
the values of indices coincide each other.  
\begin{lemma}\label{220816-1}
Suppose that 
\beas 
&&
{\coloraka 
\sup_{j\in\{1,...,n\},\>n\in\bbN,\>\gamma\in\{0,1,...,m\},\atop\coloraka r_1,....,r_\gamma,s\in[0,1]}
\big\| |D_{r_1,....,r_\gamma}{\coloraka \calo}|\big\|_p
<\infty
}
\eeas
{\coloraka for all $\calo=a^{(*)}_{n,j}$, $b^{(*)}_{n,j}$ and $a^{(*)}_{n,j,i}(s)$, and } 
for every $p>1$.\footnote{$\gamma=0$ denotes the case with no derivative.} 
Then 
\bd
\item[(a)] 
Suppose that $b^{(d)}_{n,k}$ are $\calf_{t_{k-1}}$-measurable. 
Then for $\nu_1,...,\nu_m,\mu_1,...,\mu_q\in\bbN$, 
\beas 
\frac{1}{n^m}\sum_{j_1,....,j_m}^n E\bigg[
a^{(1)}_{n,j_1}J^{(\nu_1)}_{n,j_1}\cdots a^{(m)}_{n,j_m}J^{(\nu_m)}_{n,j_m}
\bigg(\frac{1}{\sqrt{n}}\sum_{k_1=j_1+1}^nb^{(1)}_{n,k_1}J^{(\mu_1)}_{n,k_1}\bigg)
\cdots
\bigg(\frac{1}{\sqrt{n}}\sum_{k_q=j_m+1}^nb^{(q)}_{n,k_q}J^{(\mu_q)}_{n,k_q}\bigg)
\bigg]
&=&
O\bigg(\frac{1}{n^{m/2}}\bigg). 
\eeas
\item[(b)] 
For $\nu_1,...,\nu_m\in\bbN$, 
\beas 
\frac{1}{n^m}\sum_{j_1,....,j_m}^n E\bigg[
a^{(1)}_{n,j_1}J^{(\nu_1)}_{n,j_1}\cdots a^{(m)}_{n,j_m}J^{(\nu_m)}_{n,j_m}\bigg]
&=&
O\bigg(\frac{1}{n^{m/2}}\bigg). 
\eeas
\ed
{\coloraka The constants in the above estimates depend only on the given supremums. } 
\end{lemma}
\proof 
First we will show (a). 
We use the $L^2([0,T])$-orthogonality between $1_{(t_{j-1},t_j]}$ and $1_{(t_{k-1},t_k]}$ for $j\not=k$. 
If the number of single $j_*$'s is $\alpha$, then the outside summation has 
at most $n^\alpha\times n^{\frac{m-\alpha}{2}}$ terms of such type. 
The $\alpha$ times IBP-formula for those single $j_*$'s deduces the order $n^{-m/2}$ 
if $k_1,...,k_q$ are different from any of $j_*$'s; otherwise, we also get $n^{-1/2}$ 
in each  IBP-formula. 
Note also that the derivative of $b$'s do not change the form of ``martingale'', besides 
it gives $n^{-1/2}$. 
After all, total order becomes 
\beas 
n^{-m}\times n^\alpha\times n^{\frac{m-\alpha}{2}}\times n^{-\half \alpha}
=n^{-\frac{m}{2}}. 
\eeas
In a similar way, we can obtain (b). 
\qed\\

For example, we apply (b) for 
\beas
a^{(c)}_{n,j}
&=& 
2a(w_{t_{j-1}})\times n\int_{I_j(t)}\bigg(4\int_r^1 aa'(w_s)ds\bigg)dr. 
\eeas
By Lemma \ref{220816-1}, we see 
\beas 
\bigg\|
\frac{1}{n}\sum_{j=1}^n 2a(w_{t_{j-1}})\eta_j(t)
\bigg(\frac{1}{\sqrt{n}}\sum_{k=j+1}^na'(w(t_{k-1}))\xi_k(t)\bigg)
\bigg\|_p
&=&
\bigg(\frac{1}{\sqrt{n}}\bigg)
\eeas
and 
\beas
\bigg\|
\sum_{j=1}^n 2a(w_{t_{j-1}})\eta_j(t)\int_{I_j(t)}4\int_r^1aa'(w_s)ds \>dr
\bigg\|_p
&=&
\bigg(\frac{1}{\sqrt{n}}\bigg)
\eeas
as $n\to\infty$ for every $p>1$. 
Consequently, we obtain (\ref{220820-1}).

Now we have 
\bea\label{220820-2} 
\det\tilde{\sigma}(n,t)
&=&
\frac{1}{n}\sum_{j=1}^n \bigg[2a(w_{t_{j-1}})\eta_j(t)\bigg]^2 \int_0^t\bigg[4\int_r^1aa'(w_s)ds\bigg]^2dr
\nn\\&&
+\sum_{j=1}^n \bigg[\frac{1}{\sqrt{n}}\sum_{k=j+1}^na'(w(t_{k-1}))\xi_k(t)\bigg]^2 
\times\frac{1}{n}\sum_{j=1}^n\int_{I_j(t)}\bigg[4\int_r^1aa'(w_s)ds\bigg]^2dr
\nn\\&&
-\bigg\{
\sum_{j=1}^n \bigg[
\frac{1}{\sqrt{n}}\sum_{k=j+1}^na'(w(t_{k-1}))\xi_k(t)\bigg]
\cdot \int_{I_j(t)}4\int_r^1aa'(w_s)ds \>dr
\bigg\}^2
\nn\\&\geq&
\frac{1}{n}\sum_{j=1}^n \bigg[2a(w_{t_{j-1}})\eta_j(t)\bigg]^2 \int_0^t\bigg[4\int_r^1aa'(w_s)ds\bigg]^2dr,
\eea
where we used the Schwarz inequality 
as well as $|I_j(t)|\leq 1/n$.

Let $t_n=1/2$ and $c_0=\inf_{x\in\bbR}|a(x)|$. 
\begin{en-text}
Define $\xi_n$ and $s_n$ by 
\beas 
\xi_n &=& 
\int_{[0,1]^2}
\bigg(
\frac{|C^n_t-C^\infty_t-C^n_s+C^\infty_s|n^q}{|t-s|^{3/8}}\bigg)^8 dtds. 
\eeas
\end{en-text}
Define $s_n$ by 
%
\beas 
s_0:=s_n &=& \half c_0^2 \int_0^\half\bigg[4\int_r^1aa'(w_s)ds\bigg]^2dr. 
\eeas
In particular, 
$s_n$ does not depend on $n$ in this case. 

We consider the system of stochastic differential equations: 
\beas 
\l\{\begin{array}{rcl}
d\>w_t &=& d\>w_t,\\
d\>f_t &=& 2\alpha(w_t)\>d\>t . 
\end{array}\r.
\eeas
For this system, 
we have $V_0=2\alpha(w)\partial_2$, 
$V_1=\partial_1$, 
where $\partial_1$ and $\partial_2$ correspond to 
$w$ and $f$, respectively. 
We see that Condition [HW2] together with [HW1] applied at $w(0)$  
implies the H\"ormander condition for this system; 
see for example Ikeda and Watanabe. 
The boundedness of derivatives of the coefficients 
assumed there can be removed in our case by means of 
a large deviation argument. 
In particular, $F_\infty=f_1$ is nondegenerate, 
even up to $t=1/2$, that is, 
\bea\label{220820-4}
s_0^{-1}\in L^{\infty-}. 
\eea
\begin{en-text}
Thus it is not difficult to show the uniform nondegeracy of 
$F_n$, if necessary, under suitable truncation which is 
exponentially small. 
\end{en-text}
Since $M^n_t$ and $F_n$ are asymptotically orthogonal 
in the $H$-space in the sense of (\ref{220820-2}), 
the convergence is sufficiently fast as shown by (\ref{220820-1}), 
and $s_0$ is nondegenerate as (\ref{220820-4}), 
we conclude the uniform nondegeneracy of $(M^n_t,F_\infty)$ 
for $t\geq 1/2$, as follows. 
Due to $\sigma_{(M^n_t,F_\infty)}\geq \sigma(n,t)$ by definition, for every $K>0$ 
{\coloraka and $t^\dagger_n=\min\{t_j;t_j\geq 1/2\}$}, 
\beas 
\sup_{t\geq1/2}P\bigg[\det\sigma_{(M^n_t,F_\infty)}<s_0\bigg]
&\leq&
P\bigg[\det\sigma(n,{\coloraka t^\dagger_n})<{\coloraka 1.5}s_0\bigg] {\coloraka +O(n^{-K})}
\\&\leq&
P\bigg[\det\tilde{\sigma}(n,{\coloraka t^\dagger_n})<2s_0\bigg]+O(n^{-K})
\\&=&
O(n^{-K})
\eeas
as $n\to\infty$.

\begin{en-text}
{\bf $\bigg[$} 
The following lemma does not hold and we do not use it!! 
\begin{lemma}\label{220816-2}
\beas 
\frac{1}{n^m}\sum_{j_1,....,j_m}^n E\bigg[
a^{(1)}_{n,j_1}\cdots a^{(m)}_{n,j_m}
\bigg(\frac{1}{\sqrt{n}}\sum_{k_1=j_1+1}^nb^{(1)}_{n,k_1}J^{(\mu_1)}_{n,k_1}\bigg)
\cdots
\bigg(\frac{1}{\sqrt{n}}\sum_{k_q=j_m+1}^nb^{(m)}_{n,k_q}J^{(\mu_m)}_{n,k_m}\bigg)
\bigg]
&=&
O\bigg(\frac{1}{n^{m/2}}\bigg). 
\eeas
\end{lemma}
{\bf $\bigg]$}

\koko
}}

}\coloroy{
By elementary argument, 
\beas 
\bigg\|\|D_1(n,1)\|_H^2-4\int_0^1a(w_t)^2dt\bigg\|_p &=& O(n^{-\half})
\eeas
for every $p>1$. 
[The first term on the right-hand side is dominating. $\leftarrow$ not true!]
On an extension of the original stochastic basis, 
there exists a Wiener process $B=(B_t)_{t\in[0,1]}$ independent of $\calf_1$ such that 
\beas 
\bigg(w,n^{-\half}\sum_{j=2}^{[n\cdot]} a'(w_{t_{j-1}})\xi_j\bigg)
&\to^d&
\bigg(w,\int_0^\cdot a'(w_t)\sqrt{2}dB_t\bigg)
\eeas
in $C([0,1];\bbR^2)$ as $n\to\infty$; 
this is a consequence of the theory of convergence of stochastic integrals 
by Jakbowskii et al. \cite{jak} or Kurtz and Protter \cite{Kur-Pro90}. 
Therefore, 
\beas 
\bigg\| \|D_2(n,1)\|_H^2 - \int_0^1 
\big(\int_r^1 a'(w_t)\sqrt{2}dB_t\big)^2 dr \bigg\|_p 
\koko
\eeas

}}
\koko

\beas 
D_rF_n &=& 
\frac{4}{n}\sum_{j=1}^n a'(w_{t_{j-1}})1_{(0,t_{j-1}]}(r).
\eeas

\beas 
|\langle DM^n_t,F_n\rangle |
&\leq&
\langle D_2(n,t) \rangle^\half \langle DF_n\rangle^\half
\\&\to&
0
\eeas
[An estimate of the convergence can be obtained. ]
\end{en-text}

\subsubsection{Proof of Theorem \ref{220820-3}}

Theorem \ref{220820-3} now follows from the results of the preceding subsections.

\subsection{Studentization}
We shall consider the expansion of the expectation $E[g(F_n^{-\half}M^n_1)]$. 
This form corresponds to a studentized statistic in the statistical context. 
The contribution of the principal part is given by 
\beas 
\int g\Big(\frac{z}{\sqrt{x}}\Big) \phi(z;0,x)E[\delta_x({\mathfrak C}_0)]\>dzdx
&=&
\int g(z)\phi(z;0,1)dz\> \int p^{C_\infty}(x)dx
=
\int g(z)\phi(z;0,1)dz\
\eeas
%

Let $g\in\cals(\bbR)$. 
For the second-order terms, we have 
{\coloro{
\beas 
\int g\Big(\frac{z}{\sqrt{x}}\Big) p_1(z,x)\>dzdx
&=&
\int g\Big(\frac{z}{\sqrt{x}}\Big) 
E[{\mathfrak C}_1\delta_x({\mathfrak C}_0)]  \partial_z^2\Big\{{\coloro{z}}\phi(z;0,x)\Big\}
\>dzdx
\\&=&
\int \frac{1}{x}g''\Big(\frac{z}{\sqrt{x}}\Big) 
E[{\mathfrak C}_1\delta_x({\mathfrak C}_0)]  z\phi(z;0,x)
\>dzdx
\\&=&
\int g''(z)
E[{\mathfrak C}_1\delta_x({\mathfrak C}_0)]  \frac{z}{\sqrt{x}}\phi(z;0,1)\>dzdx
\\&=&
\int 
E[{\mathfrak C}_1 \delta_x({\mathfrak C}_0)] \frac{1}{\sqrt{x}}  \>dx \cdot 
\int g(z) (z^3-3z)\phi(z;0,1)\>dz
\\&=&
E\Big[ \frac{{\mathfrak C}_1}{\sqrt{{\mathfrak C}_0}} \Big]\cdot 
\int g(z) (z^3-3z)\phi(z;0,1)\>dz
\eeas
}}
Here we used 
\beas 
\int E[{\mathfrak C}_1 \delta_x({\mathfrak C}_0)] \frac{1}{\sqrt{x}}  \>dx
&=&
\int 
E\Big[\frac{{\mathfrak C}_1}{\sqrt{{\mathfrak C}_0}}
\delta_x({\mathfrak C}_0)\Big]  \>dx
\\&=&
-\int \partial_x E\Big[ \frac{{\mathfrak C}_1}{\sqrt{{\mathfrak C}_0}}
1_{\bbR_+}({\mathfrak C}_0-x)\Big] \>dx
\\&=&
E\Big[ \frac{{\mathfrak C}_1}{\sqrt{{\mathfrak C}_0}} \Big]
\eeas

Define polynomials $P_{\beta,\nu}(z,x)$ by 
\beas 
(-\partial_x)^\beta g\Big(\frac{z}{\sqrt{x}}\Big)
&=&
\sum_{\nu\leq\beta}P_{\beta,\nu}\Big(\frac{z}{\sqrt{x}},\frac{1}{\sqrt{x}}\Big)
g^{(\nu)}\Big(\frac{z}{\sqrt{x}}\Big) 
\eeas
for a $\beta$-times differentiable function $g$. 
Set 
\beas 
Q_{\alpha,\beta,\nu}(z,x)&=& x^\alpha P_{\beta,\nu}(z,x). 
\eeas
For $g\in\cals(\bbR)$ and any smooth functional $\mathfrak D$. 
we have 
%
\beas &&
\int g\Big(\frac{z}{\sqrt{x}}\Big)
\partial_z^\alpha\partial_x^\beta 
\Big\{E\Big[{\mathfrak D}\delta_x({\mathfrak C}_0)\Big]\phi(z;0,x)\Big\}\>dzdx
\\&=&
\int g\Big(\frac{z}{\sqrt{x}}\Big)
\partial_x^\beta 
\Big\{E\Big[{\mathfrak D}\delta_x({\mathfrak C}_0)\Big]
x^{-\half(1+\alpha)}\partial_y^\alpha\phi(y;0,1)|_{y=zx^{-1/2}}\Big\}\>dzdx
\\&=&
\int \sum_{\nu\leq\beta}P_{\beta,\nu}\Big(\frac{z}{\sqrt{x}},\frac{1}{\sqrt{x}}\Big)
g^{(\nu)}\Big(\frac{z}{\sqrt{x}}\Big)
\Big\{E\Big[{\mathfrak D}\delta_x({\mathfrak C}_0)\Big]
x^{-\half(1+\alpha)}\partial_y^\alpha\phi(y;0,1)|_{y=zx^{-1/2}}\Big\}\>dzdx
\\&=&
\int \sum_{\nu\leq\beta}P_{\beta,\nu}\Big(y,\frac{1}{\sqrt{x}}\Big)
g^{(\nu)}(y)
E\Big[{\mathfrak D}\delta_x({\mathfrak C}_0)\Big]
x^{-\half\alpha}\partial_y^\alpha\phi(y;0,1)\>dydx
\\&=&
\int g(y)
\sum_{\nu\leq\beta}(-\partial_y)^{(\nu)}\Big\{
\partial_y^\alpha\phi(y;0,1)
\int P_{\beta,\nu}\Big(y,\frac{1}{\sqrt{x}}\Big)
E\Big[{\mathfrak D}\delta_x({\mathfrak C}_0)\Big]
x^{-\half\alpha}\>dx\Big\}
\>dy
\\&=&
\int g(y)
\sum_{\nu\leq\beta}(-\partial_y)^{(\nu)}\Big\{
\partial_y^\alpha\phi(y;0,1)
\int 
E\Big[Q_{\alpha,\beta,\nu}\Big(y,\frac{1}{\sqrt{{\mathfrak C}_0}}\Big)
{\mathfrak D}\delta_x({\mathfrak C}_0)\Big]
\>dx\Big\}
\>dy,
\eeas
therefore we obtain the formula
\beas &&
\int g\Big(\frac{z}{\sqrt{x}}\Big)
\partial_z^\alpha\partial_x^\beta 
\Big\{E\Big[{\mathfrak D}\delta_x({\mathfrak C}_0)\Big]\phi(z;0,x)\Big\}\>dzdx
\\&=&
\int g(y)
\sum_{\nu\leq\beta}(-\partial_y)^{(\nu)}\Big\{
\partial_y^\alpha\phi(y;0,1)
E\Big[Q_{\alpha,\beta,\nu}\Big(y,\frac{1}{\sqrt{{\mathfrak C}_0}}\Big)
{\mathfrak D}\Big]\Big\}
\>dy. 
\eeas
%
%
\begin{en-text}
\beas &&
\int g\Big(\frac{z}{\sqrt{x}}\Big)
\partial_z^\alpha\partial_x^\beta 
\Big\{E\Big[{\mathfrak D}\delta_x({\mathfrak C}_0)\Big]\phi(z;0,x)\Big\}\>dzdx
\\&=&
\int 
(-1)^{\alpha}x^{-\half\alpha}
g^{(\alpha)}\Big(\frac{z}{\sqrt{x}}\Big)
\partial_x^\beta 
\Big\{E\Big[{\mathfrak D}\delta_x({\mathfrak C}_0)\Big]\phi(z;0,x)\Big\}\>dzdx
\\&=&
\int g(y)
\partial_x^\alpha
\Big\{E\Big[{\mathfrak D}\delta_x({\mathfrak C}_0)\Big]
x^{-\half\beta}\partial_y^\alpha\phi(y;0,1)\Big\}\>dydx
\eeas
\end{en-text}
By definition, 
\beas &&
Q_{\alpha,0,0}(y,x)=x^\alpha, 
\\&&
Q_{\alpha,1,0}(y,x)=0, \sskip 
Q_{\alpha,1,1}(y,x)=\half yx^{\alpha+2}
\\&&
Q_{\alpha,2,0}(y,x)=0,\sskip
Q_{\alpha,2,1}(y,x)=\frac{3}{4}yx^{\alpha+{\colorr{4}}},\sskip 
Q_{\alpha,2,2}(y,x)=\frac{1}{4}y^2x^{\alpha+{\colorr{4}}}. 
\eeas

Applying the formulas, we have 
\beas 
\int g\Big(\frac{z}{\sqrt{x}}\Big) p_2(z,x)\>dzdx
&=&
-\int g\Big(\frac{z}{\sqrt{x}}\Big)
\partial_z\Big(2\partial_z^2-4\partial_x\Big)^2 
\Big\{E[{\mathfrak C}_2\delta_x({\mathfrak C}_0)]\phi(z;0,x)\Big\}
\>dzdx
\\&=&
\int g(y) \Big[\!\!\Big[
-4\partial_y^5\phi(y;0,1)E[{\mathfrak C}_0^{-\frac{5}{2}}{\mathfrak C}_2]
\\&&
-8\partial_y\Big(y\partial_y^3\phi(y;0,1)\Big)
E[{\mathfrak C}_0^{-\frac{5}{2}}{\mathfrak C}_2]
\\&&
+12\partial_y
\Big(y\partial_y\phi(y;0,1)
E[{\mathfrak C}_0^{-\frac{{\colorr{5}}}{2}}{\mathfrak C}_2]\Big)
\\&&
-4\partial_y^2
\Big(y^2\partial_y\phi(y;0,1)
E[{\mathfrak C}_0^{-\frac{{\colorr{5}}}{2}}{\mathfrak C}_2]\Big)
\Big]\!\!\Big]\>dy
\\&=&
\int g(y) 
{\colorr{
E[{\mathfrak C}_0^{-\frac{5}{2}}{\mathfrak C}_2](12y)
}}
\phi(y;0,1)\>dy
\eeas
and 
\beas 
\int g\Big(\frac{z}{\sqrt{x}}\Big) p_3(z,x)\>dzdx
&=&
-\int g\Big(\frac{z}{\sqrt{x}}\Big)
\partial_z\Big(2\partial_z^2-4\partial_x\Big)
\Big\{E[{\mathfrak C}_3\delta_x({\mathfrak C}_0)]\phi(z;0,x)\Big\}
\>dzdx
\\&=&
\int g(y) \Big[\!\!\Big[
-2\partial_y^3\phi(y;0,1)E[{\mathfrak C}_0^{-\frac{3}{2}}{\mathfrak C}_3]
-2\partial_y\Big(y\partial_y\phi(y;0,1)E[{\mathfrak C}_0^{-\frac{3}{2}}{\mathfrak C}_3]\Big)
\Big]\!\!\Big]\>dy
\\&=&
\int g(y) 
E[{\mathfrak C}_0^{-\frac{3}{2}}{\mathfrak C}_3](-2y)
\phi(y;0,1)\>dy. 
\eeas

After all, we obtain 
\beas 
q_n(z) 
&=&
\phi(z;0,1)
+\frac{1}{\sqrt{n}}\Big\{
{\coloro{E[ {\mathfrak C}_0^{-\half}{\mathfrak C}_1]}}
(z^3-3z)
\\&&
+E[{\mathfrak C}_0^{-\frac{5}{2}}{\mathfrak C}_2]
{\colorr{(12z)}}
\\&&
+E[{\mathfrak C}_0^{-\frac{3}{2}}{\mathfrak C}_3](-2z)
\Big\}\phi(z;0,1)
\eeas
as the second-order approximate density to the distribution of $F_n^{-1/2}M^n_1$.

\begin{comment}\rm 
21.11.07 checked thus far. 
\end{comment}

\section{Quadratic form of a diffusion process}\label{220915-3}

We shall apply and extend the result in Section \ref{220826-1} 
to the quadratic variation of a diffusion process. 
It is called a realized volatility in financial context recently. 
We consider a diffusion process 
satisfying the It\^o integral equation (\ref{240221-1}). 
Here
$b$ and $\sigma$ are assumed to be smooth with 
bounded derivatives {\coloraka of positive order}.

For simplicity we only treat the one-dimensional case; 
multivariate analogue is straightforward. 
Even extension to It\^o processes is also possible 
but the descriptions would be involved. 
We write $b_t$ for $b(X_t)$ and $\sigma_t$ for $\sigma(X_t)$. 
The It\^o decomposition of $\sigma_t=\sigma(X_t)$ is 
denoted by 
\beas 
\sigma_t &=& \sigma_0 + \int_0^t \sigma^{[1]}_s dw_s
+\int_0^t \sigma^{[0]}_sds. 
\eeas
Though $\sigma^{[1]}_s$ and $\sigma^{[0]}_s$ have a simple expression with $b$, $\sigma$ and $X_s$, 
those symbols are convenient to simplify the notation.  
This rule will be applied for other functionals.

We consider the quadratic form (\ref{240221-10}) of the increments of $X$ 
with strongly predictable kernel. 
\begin{en-text}
$t_0\leq t_1\leq \cdots \leq t_n=1$ gives 
a partition of $[0,1]$. 
Note that we admit multiple of times, so 
the total number of subintervals 
may be less than $n$. Also note that 
$t_j$ can depend on $n$. 
\end{en-text}
Here we are interested in the asymptotic expansion of the 
normalized error 
\beas
Z_n&=& 
\sqrt{n}(U_n-U_\infty)
\eeas
for $U_\infty$ in (\ref{240221-12}). 
We are assuming $c\in C^\infty_\up(\bbR)$.

\subsection{Stochastic expansion} 

We will need a stochastic expansion of $Z_n$. 

\begin{lemma} $Z_n$ admits the following stochastic expansion:
\beas 
Z_n &=& 
M^n_1+\frac{1}{\sqrt{n}}N_n, 
\eeas
where 
\beas 
M^n_t &=& 
\sqrt{n}\sum_{j=1}^n 2c_{t_{j-1}}\sigma_{t_{j-1}}^2
{\coloraka \int_{t_{j-1}\wedge t}^{t_j\wedge t}}
\int_{t_{j-1}}^sdw_rdw_s,
\eeas
and 
\beas 
N_n &=& 
%
6n\sum_{j=1}^n c_{t_{j-1}}\sigma_{t_{j-1}}\sigma_{t_{j-1}}^{[1]}
\int_{t_{j-1}}^{t_j} \int_{t_{j-1}}^t \int_{t_{j-1}}^s
dw_udw_sdw_t
\\&&
+2\sum_{j=1}^n c_{t_{j-1}}b_{t_{j-1}}\sigma_{t_{j-1}}
\int_{t_{j-1}}^{t_j}dw_t
+
2n\sum_{j=1}^n c_{t_{j-1}}\sigma_{t_{j-1}} \sigma_{t_{j-1}}^{[1]} 
\int_{t_{j-1}}^{t_j}(t-t_{j-1}) dw_t
\\&&
{\coloraka
+
n^{-1}\sum_{j=1}^n  c_{t_{j-1}} b_{t_{j-1}}^2 
+n^{-1}\sum_{j=1}^n c_{t_{j-1}} \sigma_{t_{j-1}} b^{[1]}_{t_{j-1}} 
}
\\&&
{\coloro{-}}n\sum_{j=1}^n  c^{[1]}_{t_{j-1}} \sigma_{t_{j-1}}^2 
\int_{t_{j-1}}^{t_j} \int_{t_{j-1}}^t dw_sdt
\\&&
{\coloro{
-\frac{1}{2n}\sum_{j=1}^n c^{[0]}_{t_{j-1}}\sigma_{t_{j-1}}^2
-\frac{1}{n}\sum_{j=1}^n c^{[1]}_{t_{j-1}}\sigma_{t_{j-1}}\sigma^{[1]}_{t_{j-1}}
}}
+o_M(1). 
\eeas
Here $o_M(1)$ denotes a term of $o(1)$ as $n\to\infty$ 
with respect to $\bbD_{s,p}$-norms of any order. 
{\coloraka The families $\{M^n_t\}_{t\in[0,1],n\in\bbN}$ and $\{N_n\}_{n\in\bbN}$ are bounded in every $\bbD_{s,p}$-norm. }
\end{lemma}
By somewhat long computations, it is possible to obtain the above lemma. 
We omit details.

\subsection{Asymptotic expansion}

For a reference variable, we will consider 
\beas 
F_n = 
{\coloro{
\frac{1}{n}\sum_{j=1}^n \beta(X_{t_{j-1}})
}}
\hspace{5mm}\mbox{or} \hspace{5mm}
F_n = 
F_\infty:=
{\coloro{
\int_0^1\beta(X_t) dt,
}}
\eeas
where $\beta\in C^\infty_\up(\bbR,\bbR^{d_1})$. 
The results will be the same in these cases. 
It is statistically natural to consider those functionals because, for example, $F_\infty$ gives 
the conditional asymptotic variance of the estimation error $Z^n_1$ 
as in Section \ref{220826-1}.  
We will derive asymptotic expansion of the joint distribution of $(M^n_1,F_n)$. 

Let $a(x)=c(x)\sigma(x)^2$. 
Let 
\beas 
V_0(x_1,x_2) = 
\l[\begin{array}{c} b(x_1)-\half \sigma(x_1)\partial_{x_1}\sigma(x_1)\y
\beta(x_1) \end{array} \r]
&\mbox{  and  }&
V_1(x_1,x_2)=
\l[\begin{array}{c} \sigma(x_1)\y
0 \end{array} \r]
\eeas
for $x_1\in\bbR$ and $x_2\in\bbR^{d_1}$. 
The Lie algebra generated by 
\beas 
V_1, \>[V_i,V_j]\>(i,j=0,1),\>[V_i,[V_j,V_k]]\>(i,j,k=0,1),....
\eeas
at $(x_1,x_2)$ is denoted by $\mbox{Lie}[V_0;V_1](x_1,x_2)$. 

Assume that $\mbox{supp}(X_0)$ is compact. 
Moreover, for nondegeneracy, we assume
\bd
\item[[H1\!\!]]\hspace{5mm}
$\inf_{x\in\bbR} |a(x)| > 0.$
\ed
\bd
\item[[H2\!\!]]\hspace{5mm}
$\mbox{Lie}[V_0;V_1](X_0,0)=\bbR^{1+d_1}$ a.s. 
\ed

\begin{remark}\rm 
Under [H1], both $\mbox{ess.}\inf|\sigma(X_0)|>0$ and $\mbox{ess.}\inf|c(X_0)|>0$. 
Then [H2] is equivalent to the linear hull 
$\mbox{L}[\partial_{x_1}^i\beta(X_0);\>i\in\bbN]=\bbR^{d_1}$ a.s. 
It is rather simple but we prefer to keep [H2], which is suitable for more general form of $F_\infty$.  
\end{remark}
\begin{remark}\rm 
Condition [H1] is usually from the uniform ellipticity of the diffusion process $X_t$ 
and a reasonable choice of the estimator for the quadratic variation. 
In this sense, it is a natural assumption in statistical context. 
\end{remark}
\begin{remark}\rm 
Consider a $(1+d_1)$-dimensional stochastic integral equation
\beas 
\check{X}_t &=& \check{X}_0 +\int_0^tV_0(\check{X}_s)ds
+\int_0^tV_1(\check{X}_s)\circ dw_s , \sskip  t\in[0,1]. 
\eeas
Then $\check{X}_1=(X_1,F_\infty)$, and the nondegeneracy condition entails 
the nondegeneracy of $F_\infty$ in particular. 
\end{remark}

We see $\tilde{W}_\infty(z)=0$ and $\tilde{F}_\infty(z)=0$. 
It is necessary to specify the limit $(M_\infty,\dotc_\infty,N_\infty)$. 
The ``martingale part'' of $N_n$ with respect to $\F^n$ is given by 
\beas 
\dot{N}^n_t
&=&
6n\sum_{j:t_j\leq t} c_{t_{j-1}}\sigma_{t_{j-1}}\sigma_{t_{j-1}}^{[1]}
\int_{t_{j-1}}^{t_j} \int_{t_{j-1}}^t \int_{t_{j-1}}^s
dw_udw_sdw_t
\\&&
+2\sum_{j:t_j\leq t}  c_{t_{j-1}}b_{t_{j-1}}\sigma_{t_{j-1}}
\int_{t_{j-1}}^{t_j}dw_t
+
2n\sum_{j:t_j\leq t}  c_{t_{j-1}}\sigma_{t_{j-1}} \sigma_{t_{j-1}}^{[1]} 
\int_{t_{j-1}}^{t_j}(t-t_{j-1}) dw_t
\\&&
{\coloraka
-n\sum_{j:t_j\leq t}   c^{[1]}_{t_{j-1}} \sigma_{t_{j-1}}^2 
\int_{t_{j-1}}^{t_j} \int_{t_{j-1}}^t dw_sdt. 
}
\eeas
We redefine $\bar{M}^{2,n}$ and $\bar{M}^{\xi,n}$ by the same equations 
in Section \ref{220826-5} but with $a(X_{t_{j-1}})$ in place of $a(w_{t_{j-1}})$. 
Then
\beas 
\dla \bar{w}^n,\bar{w}^n \dra_t &=& \frac{[nt]}{n}\to t
\\
\dla \bar{M}^{2,n},\bar{M}^{2,n} \dra_t 
&=& 
\frac{2}{n}\sum_{j:t_{j-1}\leq t} a(X_{t_{j-1}})^2 
\to^p 
2\int_0^t a(X_s)^2 ds
\\
\dla \bar{M}^{\xi,n},\bar{M}^{\xi,n} \dra_t 
&=&
\frac{16}{3n}\sum_{j:t_{j-1}\leq t} a(X_{t_{j-1}})^4
\to^p 
\frac{16}{3}\int_0^t a(X_s)^4ds 
\\ 
{\coloraka \dla\bar{w}^n,\dot{N}\dra_t} &{\coloraka \to^p}&{\coloraka \int_0^t k_sds } 
\\
\dla \bar{w}^n,\bar{M}^{k,n} \dra_t &=& 0
\sskip(k=2,\xi)
\\
\dla \bar{M}^{2,n},\bar{M}^{\xi,n} \dra_t &=& 
\frac{{\coloro{8}}}{3n}\sum_{j:t_{j-1}\leq t} a(X_{t_{j-1}})^3
\to^p 
\frac{{\coloro{8}}}{3}\int_0^t a(X_s)^3 ds,
\\
\dla \bar{M}^{2,n},\dot{N}^n \dra_t &\to^p&0,
\\
\dla \bar{M}^{\xi,n},\dot{N}^n \dra_t &\to^p&0,
\\
\dla \dot{N}^n,\dot{N}^n \dra_t &\to^p& \int_0^t q_s^2 ds
\eeas
as $n\to\infty$ for each $t\in[0,1]$, 
where $\bbR_+$-valued process $q_t$ 
takes the form 
\beas 
q_t^2 &=& p(c_t,c^{[1]}_t,b_t,\sigma_t,\sigma^{[1]}_t)
\eeas
for some polynomial $p$; 
\begin{en-text}
is defined by 
\beas 
q_t^2 &=&
\frac{22}{3}c_t^2\sigma_t^2(\sigma^{[1]}_t)^2
+\frac{1}{3}(c^{[1]}_t)^2\sigma_t^4
+\frac{2}{3}c_tc^{[1]}_t\sigma_t^3\sigma^{[1]}_t
+\frac{28}{3}c_t^2b_t^2\sigma_t^2
+\frac{16}{3}c_t^2b_t\sigma_t^2\sigma^{[1]}_t
+\frac{10}{3}c_tb_t\sigma_t^3c^{[1]};
\eeas
\end{en-text}
it is possible to give an explicit expression of $p$, 
however we do not need the precise form of $q_t$ later. 
The orthogonality of $\bar{M}^{2,n}$, $\bar{M}^{\xi,n}$ and $\dot{N}^n$ 
to any bounded martingale orthogonal to $w$ is obvious, thus 
with a representation of $\dotc\>\!\!^n_t$ in Section \ref{220826-5} 
with $a(X_{t_{j-1}})$ for $a(w_{t_{j-1}})$, 
and 
those of $M^n_t$ and $N_n$, we obtain
\beas 
(M_\infty,\dotc_\infty,N_\infty)
&=^d&
\Big( \int_0^1 \sqrt{2}a(X_s)dB_s,
\int_0^1 \frac{{\coloro{4}}\sqrt{2}}{3}a(X_s)^2dB_s
+\int_0^1 \frac{{\coloro{4}}}{3}a(X_s)^2dB'_s,
\\&&
{\coloraka \int_0^1 k_sdw_s+}\int_0^1 {\coloraka \sqrt{q_s^2-k_s^2}}dB''_s+\int_0^1 h_sds
\Big), 
\eeas
where $(B,B', B'')$ is a three-dimensional standard Wiener process, 
independent of $\calf$,  
defined on the extension $\bar{\Omega}$, 
and 
\beas 
h_t&=&
c_tb_t^2+c_tb^{[1]}_t\sigma_t-\half c^{[0]}_t\sigma_t^2-c^{[1]}_t\sigma_t\sigma^{[1]}_t.
\eeas 
and 
{\coloraka 
\beas 
k_t &=& 2c_tb_t\sigma_t+c_t\sigma_t\sigma^{[1]}_t-\half c^{[1]}_t\sigma_t^2.
\eeas
}
Since 
\beas 
\tilde{N}_\infty(z) &=& {\coloraka \int_0^1 k_tdw_t+}\int_0^1 h_tdt, 
\eeas
the random symbol $\underline{\sigma}(z,iu,iv)$ is given by 
\beas 
\underline{\sigma}(z,iu,iv)
&=&
{\coloro{\frac{2z}{3}\int_0^1a(X_s)^3ds\Big( \int_0^1a(X_s)^2ds\Big)^{-1}}} \>(iu)^2 
{\coloraka +iu\int_0^1 k_t dw_t}
+iu\int_0^1 h_tdt.
\eeas

Let us find the anticipative random symbol $\overline{\sigma}(iu,iv)$. 
Recall that $\alpha(x)=a(x)^2$, 
\beas 
C_s^\infty=2\int_0^s\alpha(X_t)dt, \sskip
C_\infty=2\int_0^1 \alpha(X_t)dt, \sskip
F_\infty=\int_0^1\beta(X_t)dt\sskip\mbox{and}\sskip
W_\infty=0. 
\eeas 
The random symbol $\sigma_{s,r}(iu,iv)$ admits the expression
\beas 
&&
\sigma_{s,r}(iu,iv)
\\&=&
u^2\int_r^s\alpha'(X_t)D_rX_tdt
\Big(- u^2\int_s^1\alpha'(X_t)D_sX_tdt+i\int_s^1\beta'(X_t)[v]D_sX_tdt \Big)
\\&&
+
\Big(-u^2\int_r^1\alpha'(X_t)D_rX_tdt+i\int_r^1\beta'(X_t)[v]D_rX_tdt \Big)
\Big(-u^2\int_s^1\alpha'(X_t)D_sX_tdt+i\int_s^1\beta'(X_t)[v]D_sX_tdt \Big)
\\&&
+
\Big(-u^2\int_{r\vee s}^1\{\alpha''(X_t)D_rX_tD_sX_t
+\alpha'(X_t)D_rD_sX_t\}dt
+i\int_{r\vee s}^1\{\beta''(X_t)[v]D_rX_tD_sX_t
+\beta'(X_t)[v]D_rD_sX_t\}dt
\Big)
\eeas
for $r\leq s$, 
where the prime $'$ stands for the derivative in $x_1\in\bbR$. 
The processes $D_sX_t$ and $D_rD_sX_t$ are determined according to routine; 
$D_sX_t$ satisfies the equation
\beas 
D_sX_t &=&
\sigma(X_s)+\int_s^t {\coloraka b}'(X_{t_1})D_sX_{t_1}d{t_1}+\int_s^t\sigma'(X_{t_1})D_sX_{t_1}dw_{t_1}
\eeas
for $s\leq t$
{\coloraka, and  
\beas
D_rD_sX_t &=&
\sigma'(X_s)D_rX_s
+\int_s^t{\coloraka b}''(X_{t_1})D_rX_{t_1}D_sX_{t_1}d{t_1}
+\int_s^t{\coloraka b}'(X_{t_1})D_rD_sX_{t_1}d{t_1}
\\&&
+\int_s^t\sigma''(X_{t_1})D_rX_{t_1}D_sX_{t_1}dw_{t_1}
+\int_s^t\sigma'(X_{t_1})D_rD_sX_{t_1}dw_{t_1}
\eeas
for $r{\coloraka <} s\leq t$. 
}
Those equations form a graded system of 
partially linear equations; therefore the $L^p$-estimates of the solution are at hand. 
%
\begin{comment}
\beas 
\sigma_{s,r}(iu,iv)
&=&
\half D_rC^{{\colorr{\infty}}}_s[u^{\otimes2}]\otimes 
\Big(iD_sW_\infty[u]-\half D_sC_\infty[u^{\otimes2}]+iD_sF_\infty[v] \Big)
\\&&
+
\Big(iD_rW_\infty[u]-\half D_rC_\infty[u^{\otimes2}]+iD_rF_\infty[v] \Big)
\otimes 
\Big(iD_sW_\infty[u]-\half D_sC_\infty[u^{\otimes2}]+iD_sF_\infty[v] \Big)
\\&&
+
\Big(iD_rD_sW_\infty[u]-\half D_rD_sC_\infty[u^{\otimes2}]+iD_rD_sF_\infty[v] 
\Big). 
\eeas
\end{comment}
%
Now we obtain the anticipative random symbol 
\beas 
\bar{\sigma}(iu,iv)
&=&
 \int_0^1 
iu\>a(X_s)
\sigma_{s,s}(iu,iv)\>ds
\eeas
with
\beas 
\sigma_{s,s}(iu,iv)
&=&
\Big(-u^2\int_s^1\alpha'(X_t)D_sX_tdt+i\int_s^1\beta'(X_t)[v]D_sX_tdt \Big)^2
\\&&
-u^2\int_s^1\{\alpha''(X_t)(D_sX_t)^2
+\alpha'(X_t)D_sD_sX_t\}dt
\\&&
+i\int_s^1\{\beta''(X_t)[v](D_sX_t)^2
+\beta'(X_t)[v]D_sD_sX_t\}dt
\eeas
%
\begin{comment}
\beas 
\bar{\sigma}(iu,iv)
&=&
\half\mbox{Tr}^*  \int_0^1 
\bar{K}^\infty(t,t)[iu]\otimes 
\sigma_{t,t}(iu,iv)\>{\coloroy{\mu (dt)}}. 
\eeas
\end{comment}
%
As before, define the total random symbol $\sigma$ by (\ref{220829-1}) and 
the density function $p_n(z,x)\in C^\infty(\bbR^{1+d_1})$ by (\ref{220829-2}). 
{\coloraka By a quite similar} 
argument as we proved Theorem \ref{220820-3}, 
it is easy to obtain the following theorem. 
{\coloraka See \cite{Yoshida2012arXiv} for details. }

\begin{theorem}\label{220829-3}
Suppose that {\rm [H1]} and {\rm[H2]} are satisfied. 
Then for any positive numbers $M$ and $\gamma$, 
\beas 
\sup_{f\in\cale(M,\gamma)} \bigg|
E\big[f(Z_n,F_n)\big]-\int_{\bbR^{1+d_1}}f(z,x)p_n(z,x)dzdx\bigg|
&=&
o\bigg(\frac{1}{\sqrt{n}}\bigg)
\eeas
as $n\to\infty$, 
where $\cale(M,\gamma)$ is the set of measurable functions $f:\bbR^{1+d_1}\to\bbR$ 
satisfying $|f(z,x)|\leq M(1+|z|+|x|)^\gamma$ for all $(z,x)\in\bbR\times\bbR^{d_1}$.  \\
\end{theorem}

\begin{remark}\rm
The hybrid I method with a rough Monte-Carlo method in the second order term 
is useful in the application of the expansion formula to numerical approximation. 
Obviously, the asymptotic expansion applies to statistical hypothesis testing. 
We will show its applications to prediction and option pricing in other papers.  
\end{remark}

\begin{remark}\rm 
While our method working in the mixed normal limit case 
enabled us to introduce a {\it conditioning} variable as $F_n$, 
it is possible to consider versions of our results without $F_n$. 
It will reduce the regularity condition of smoothness, that is, indices of 
differentiability 
of other variables. Of course, in that case, we only obtain a single (not joint) expansion. 
\end{remark}

\begin{remark}\rm 
Conditional limit theorems are in our scope. 
Indeed, it is also possible to obtain asymptotic expansion of the conditional distribution, 
as it was reported 
at the meetings mentioned in the footnote of p.1. 
\end{remark}

\begin{remark}\rm 
In $[A2^\natural]$ and $[A2]$, the index of differentiability of $N_n$ etc. could be reduced 
by more sophisticated multiple truncations. 
\end{remark}



\begin{comment}
The following versions are to be written: 
\begin{enumerate}
\item asymptotic expansion without $F_n$ to sharpen indices
\item limit theorems for the conditional distribution (in the first order)
\item asymptotic expansion of the conditional distribution (as it has been already done by myself)
\end{enumerate}
\end{comment}


\bibliographystyle{spmpsci}      
\bibliography{bibtex-20101014-20120906}   

\def\cprime{$'$} \def\cprime{$'$}
\begin{thebibliography}{10}
\providecommand{\url}[1]{{#1}}
\providecommand{\urlprefix}{URL }
\expandafter\ifx\csname urlstyle\endcsname\relax
  \providecommand{\doi}[1]{DOI~\discretionary{}{}{}#1}\else
  \providecommand{\doi}{DOI~\discretionary{}{}{}\begingroup
  \urlstyle{rm}\Url}\fi

\bibitem{BGJPS2006}
Barndorff-Nielsen, O.E., Graversen, S.E., Jacod, J., Podolskij, M., Shephard,
  N.: A central limit theorem for realised power and bipower variations of
  continuous semimartingales.
\newblock In: From stochastic calculus to mathematical finance, pp. 33--68.
  Springer, Berlin (2006)

\bibitem{BhattacharyaRanga1986}
Bhattacharya, R.N., Ranga~Rao, R.: Normal approximation and asymptotic
  expansions.
\newblock Robert E. Krieger Publishing Co. Inc., Melbourne, FL (1986).
\newblock Reprint of the 1976 original

\bibitem{Dohnal1987}
Dohnal, G.: On estimating the diffusion coefficient.
\newblock J. Appl. Probab. \textbf{24}(1), 105--114 (1987)

\bibitem{Genon-CatalotJacod1993}
Genon-Catalot, V., Jacod, J.: On the estimation of the diffusion coefficient
  for multi-dimensional diffusion processes.
\newblock Ann. Inst. H. Poincar\'e Probab. Statist. \textbf{29}(1), 119--151
  (1993)

\bibitem{GotzeHipp1983}
G{\"o}tze, F., Hipp, C.: Asymptotic expansions for sums of weakly dependent
  random vectors.
\newblock Z. Wahrsch. Verw. Gebiete \textbf{64}(2), 211--239 (1983)

\bibitem{GotzeHipp1994}
G{\"o}tze, F., Hipp, C.: Asymptotic distribution of statistics in time series.
\newblock Ann. Statist. \textbf{22}(4), 2062--2088 (1994)

\bibitem{HayashiYoshida2008a}
Hayashi, T., Yoshida, N.: Nonsynchronous covariance estimator and limit theorem
  ii.
\newblock Institute of Statistical Mathematics, Research Memorandum
  \textbf{1067} (2008)

\bibitem{IacusUchidaYoshida2009}
Iacus, S.M., Uchida, M., Yoshida, N.: Parametric estimation for partially
  hidden diffusion processes sampled at discrete times.
\newblock Stochastic Process. Appl. \textbf{119}(5), 1580--1600 (2009).
\newblock \doi{10.1016/j.spa.2008.08.004}.
\newblock \urlprefix\url{http://dx.doi.org/10.1016/j.spa.2008.08.004}

\bibitem{IkedaWatanabe1989}
Ikeda, N., Watanabe, S.: Stochastic differential equations and diffusion
  processes, \emph{North-Holland Mathematical Library}, vol.~24, second edn.
\newblock North-Holland Publishing Co., Amsterdam (1989)

\bibitem{Jacod1997}
Jacod, J.: On continuous conditional {G}aussian martingales and stable
  convergence in law.
\newblock In: S\'eminaire de Probabilit\'es, XXXI, \emph{Lecture Notes in
  Math.}, vol. 1655, pp. 232--246. Springer, Berlin (1997)

\bibitem{Kessler1997}
Kessler, M.: Estimation of an ergodic diffusion from discrete observations.
\newblock Scand. J. Statist. \textbf{24}(2), 211--229 (1997)

\bibitem{KusuokaYoshida2000}
Kusuoka, S., Yoshida, N.: Malliavin calculus, geometric mixing, and expansion
  of diffusion functionals.
\newblock Probab. Theory Related Fields \textbf{116}(4), 457--484 (2000)

\bibitem{KutoyantsYoshida2007}
Kutoyants, Y.A., Yoshida, N.: Moment estimation for ergodic diffusion
  processes.
\newblock Bernoulli \textbf{13}(4), 933--951 (2007).
\newblock \doi{10.3150/07-BEJ1040}.
\newblock \urlprefix\url{http://dx.doi.org/10.3150/07-BEJ1040}

\bibitem{MasudaYoshida2005}
Masuda, H., Yoshida, N.: Asymptotic expansion for {B}arndorff-{N}ielsen and
  {S}hephard's stochastic volatility model.
\newblock Stochastic Process. Appl. \textbf{115}(7), 1167--1186 (2005)

\bibitem{Nualart2006}
Nualart, D.: The {M}alliavin calculus and related topics, second edn.
\newblock Probability and its Applications (New York). Springer-Verlag, Berlin
  (2006)

\bibitem{PodolskijVetter2009}
Podolskij, M., Vetter, M.: Estimation of volatility functionals in the
  simultaneous presence of microstructure noise and jumps.
\newblock Bernoulli \textbf{15}(3), 634--658 (2009).
\newblock \doi{10.3150/08-BEJ167}.
\newblock \urlprefix\url{http://dx.doi.org/10.3150/08-BEJ167}

\bibitem{PrakasaRao1983}
Prakasa~Rao, B.L.S.: Asymptotic theory for nonlinear least squares estimator
  for diffusion processes.
\newblock Math. Operationsforsch. Statist. Ser. Statist. \textbf{14}(2),
  195--209 (1983)

\bibitem{PrakasaRao1988}
Prakasa~Rao, B.L.S.: Statistical inference from sampled data for stochastic
  processes.
\newblock In: Statistical inference from stochastic processes (Ithaca, NY,
  1987), \emph{Contemp. Math.}, vol.~80, pp. 249--284. Amer. Math. Soc.,
  Providence, RI (1988)

\bibitem{SakamotoYoshida1998a}
Sakamoto, Y., Yoshida, N.: Asymptotic expansion of {$M$}-estimator over
  {W}iener space.
\newblock Stat. Inference Stoch. Process. \textbf{1}(1), 85--103 (1998)

\bibitem{SakamotoYoshida2004}
Sakamoto, Y., Yoshida, N.: Asymptotic expansion formulas for functionals of
  {$\epsilon$}-{M}arkov processes with a mixing property.
\newblock Ann. Inst. Statist. Math. \textbf{56}(3), 545--597 (2004)

\bibitem{ShimizuYoshida2006}
Shimizu, Y., Yoshida, N.: Estimation of parameters for diffusion processes with
  jumps from discrete observations.
\newblock Stat. Inference Stoch. Process. \textbf{9}(3), 227--277 (2006).
\newblock \doi{10.1007/s11203-005-8114-x}.
\newblock \urlprefix\url{http://dx.doi.org/10.1007/s11203-005-8114-x}

\bibitem{SoerensenUchida2003}
S{\o}rensen, M., Uchida, M.: Small-diffusion asymptotics for discretely sampled
  stochastic differential equations.
\newblock Bernoulli \textbf{9}(6), 1051--1069 (2003).
\newblock \doi{10.3150/bj/1072215200}.
\newblock \urlprefix\url{http://dx.doi.org/10.3150/bj/1072215200}

\bibitem{Uchida2010}
Uchida, M.: Contrast-based information criterion for ergodic diffusion
  processes from discrete observations.
\newblock Ann. Inst. Statist. Math. \textbf{62}(1), 161--187 (2010).
\newblock \doi{10.1007/s10463-009-0245-1}.
\newblock \urlprefix\url{http://dx.doi.org/10.1007/s10463-009-0245-1}

\bibitem{UchidaYoshida2001}
Uchida, M., Yoshida, N.: Information criteria in model selection for mixing
  processes.
\newblock Stat. Inference Stoch. Process. \textbf{4}(1), 73--98 (2001)

\bibitem{UchidaYoshida2008}
Uchida, M., Yoshida, N.: Estimation of the volatility for a stochastic
  differential equation  (2008)

\bibitem{Watanabe1983}
Watanabe, S.: Malliavin's calculus in terms of generalized {W}iener
  functionals.
\newblock In: Theory and application of random fields (Bangalore, 1982),
  \emph{Lecture Notes in Control and Inform. Sci.}, vol.~49, pp. 284--290.
  Springer, Berlin (1983)

\bibitem{Yoshida1992b}
Yoshida, N.: Estimation for diffusion processes from discrete observation.
\newblock J. Multivariate Anal. \textbf{41}(2), 220--242 (1992)

\bibitem{Yoshida1997}
Yoshida, N.: Malliavin calculus and asymptotic expansion for martingales.
\newblock Probab. Theory Related Fields \textbf{109}(3), 301--342 (1997)

\bibitem{Yoshida2001c}
Yoshida, N.: Malliavin calculus and martingale expansion.
\newblock Bull. Sci. Math. \textbf{125}(6-7), 431--456 (2001).
\newblock Rencontre Franco-Japonaise de Probabilit{\'e}s (Paris, 2000)

\bibitem{Yoshida2001b}
Yoshida, N.: Partial mixing and conditional edgeworth expansion for diffusions
  with jumps.
\newblock Probab. Theory Related Fields \textbf{129}, 559--624 (2004)

\bibitem{Yoshida2010}
Yoshida, N.: Expansion of the asymptotically conditionally normal law.
\newblock Research Memorandum \textbf{1125, The Institute of Mathematical
  Statistics} (2010)

\bibitem{Yoshida2012arXiv}
Yoshida, N.: Asymptotic expansion for the quadratic form of the diffusion
  process.
\newblock Preprint  (2012)

\end{thebibliography}

\end{document}